\documentclass[EJP]{ejpecp} % replace ECP by EJP if needed.
\usepackage{stmaryrd}
\usepackage[normalem]{ulem}
\usepackage{graphicx,xcolor}
\usepackage{amsfonts}
\usepackage{amssymb,mathrsfs,amsthm}
\usepackage[all]{xy} 
\usepackage{color}
\usepackage{enumerate,vmargin}
\usepackage{verbatim}
\usepackage{paralist}
\usepackage[labelfont=bf,size=small,font=it]{caption}
\usepackage[numbers,sort&compress]{natbib}
\usepackage{tikz}
\usetikzlibrary{shapes,arrows}

\theoremstyle{plain}
\newtheorem{lem}{Lemma}[section]
\newtheorem{thm}[lem]{Theorem}
\newtheorem{prop}[lem]{Proposition}
\newtheorem{cor}[lem]{Corollary}
\newtheorem{alg}[lem]{Algorithm}

\theoremstyle{remark}

\hypersetup{
    unicode      = false,     % non-Latin characters in Acrobatâs bookmarks
    pdftoolbar   = true,      % show Acrobatâs toolbar?
    pdfmenubar   = true,      % show Acrobatâs menu?
    pdffitwindow = true,      % page fit to window when opened
    pdfnewwindow = true,      % links in new window
    colorlinks   = true,      % false: boxed links; true: colored links
    linkcolor    = blue,      % color of internal links
    citecolor    = blue,      % color of links to bibliography
    filecolor    = black,      % color of file links
    urlcolor     = black       % color of external links
}

\SHORTTITLE{Critical random intersection graphs}

\TITLE{Large random intersection graphs inside the critical window and triangle counts\support{supported by EPSRC grant EP/W033585/1.}\
      } % \thanks is optional. Insert line breaks with \\

\AUTHORS{%
  Minmin Wang\footnote{University of Sussex, United Kingdom.
    \EMAIL{Minmin.Wang@sussex.ac.uk}}}
\KEYWORDS{random intersection graph; random bipartite graph; scaling limit; clustering} % Separate items with ;

\AMSSUBJ{60C05, 05C80, 60F05, 60G52} % Edit. Separate items with ;
\SUBMITTED{8th Oct.~2023} % Edit.
\ACCEPTED{8th Jan.~2025} % Edit.

\ARXIVID{2309.13694} % Edit.
\VOLUME{0}
\YEAR{2020}
\PAPERNUM{0}
\DOI{10.1214/YY-TN}

\ABSTRACT{We identify the scaling limit of random intersection graphs inside their critical windows. The limit graphs vary according to the clustering regimes, and coincide with the continuum Erd\H{o}s--R\'enyi graph in two out of the three regimes.  Our approach to the scaling limit relies upon the close connection of random intersection graphs with binomial bipartite graphs, as well as a graph exploration algorithm on the latter. This further allows us to prove limit theorems for the number of triangles in the large connected components of the random intersection graphs. 
}

\newcommand{\cR}{\mathcal R}
\newcommand{\noi}{\noindent}

\newcommand{\rS}{\mathrm S}

\newcommand{\rX}{\mathrm X}

\newcommand{\rB}{\mathrm B}

\newcommand{\sL}{\mathscr L}

\newcommand{\n}{^{n, m}}
\newcommand{\m}{^{m, n}}

\newcommand{\lf}{\lfloor}
\newcommand{\rf}{\rfloor}

\newcommand{\dhaus}{\operatorname{d^{\mathrm{Haus}}_{E}}}

\newcommand{\dpr}{\operatorname{d^{\mathrm{Pr}}_{E}}}

\newcommand{\dghp}{\operatorname{\delta_{GHP}}}
\newcommand{\dTV}{\operatorname{d_{\mathrm{TV}}}}

\newcommand{\gr}{\operatorname{gr}}
\newcommand{\dgr}{d_{\gr}}

\newcommand{\CL}{\operatorname{CL}}

\newcommand{\cT}{\mathcal T}
\newcommand{\cF}{\mathcal F}
\newcommand{\cC}{\mathcal C}
\newcommand{\cO}{\mathcal O}
\newcommand{\cP}{\mathcal P}

\newcommand{\cS}{\mathcal S}
\newcommand{\cB}{\mathcal B}
\newcommand{\cL}{\mathcal L}
\newcommand{\cM}{\mathcal M}
\newcommand{\cA}{\mathcal A}
\newcommand{\cI}{\mathcal I}
\newcommand{\cV}{\mathcal V}
\newcommand{\cU}{\mathcal U}
\newcommand{\cN}{\mathcal N}
\newcommand{\cD}{\mathcal D}
\newcommand{\cE}{\mathcal E}
\newcommand{\cG}{\mathcal G}
\newcommand{\cH}{\mathcal H}

\newcommand{\rL}{\mathrm L}

\newcommand{\R}{\mathbb R}
\newcommand{\N}{\mathbb N}
\newcommand{\Z}{\mathbb Z}

\newcommand{\bi}{\mathrm{bi}}
\newcommand{\RIG}{\mathrm{RIG}}
\newcommand{\ER}{\mathrm{ER}}

\newcommand{\Var}{\operatorname{Var}}

\newcommand{\eqd}{\stackrel{(d)}=}

\begin{document}

\setcounter{tocdepth}{2}
\tableofcontents

\section{Background}
\label{sec: intro}

In this work, we investigate the asymptotic behaviours of a critical random graph with non trivial clustering properties in its large-size limit. We carry out this investigation by identifying the scaling limit of the graph. By scaling limit, we are referring to a phenomenon in which the graph of size $n\in \N$, with its edge length rescaled to $a_{n}\in (0, \infty)$, converges in a suitable sense to a non trivial object as $n\to\infty$. Studying scaling limits of random graphs provides us with a panoramic perspective of the large-sized graphs, as well as direct access to global functionals including the  diameters of the graphs. Meanwhile, as the current work shows, we can also gain insight into certain local functionals such as the numbers of triangles. Addario-Berry, Broutin and Goldschmidt \cite{ABBrGo12} first identified the scaling limit of the Erd\H{o}s--R\'enyi graph $G(n, p_{e})$, where each pair of vertices is independently connected by an edge with probability $p_{e}$. Since then, the focus has been on models of random graph with inhomogeneous degree sequences, a feature often shared by real-world networks. These include various versions of the inhomogeneous random graphs \cite{BHS18, BrDuWa21, BrDuWa22} and the configuration model \cite{BhDhvdHSe20a, DhvHvL21, CG23}. On the other hand, clustering,  another important aspect in the modelling of complex networks, has so far  received little attention regarding its role in the scaling limit of random graphs. As a first attempt to fill this gap, our work investigates the question of existence of a scaling limit for the random intersection graph model introduced in Karo\'{n}ski et al.~\cite{KaScSC99}.  

\begin{figure}[htp]
\centering
\includegraphics[height = 4.5cm]{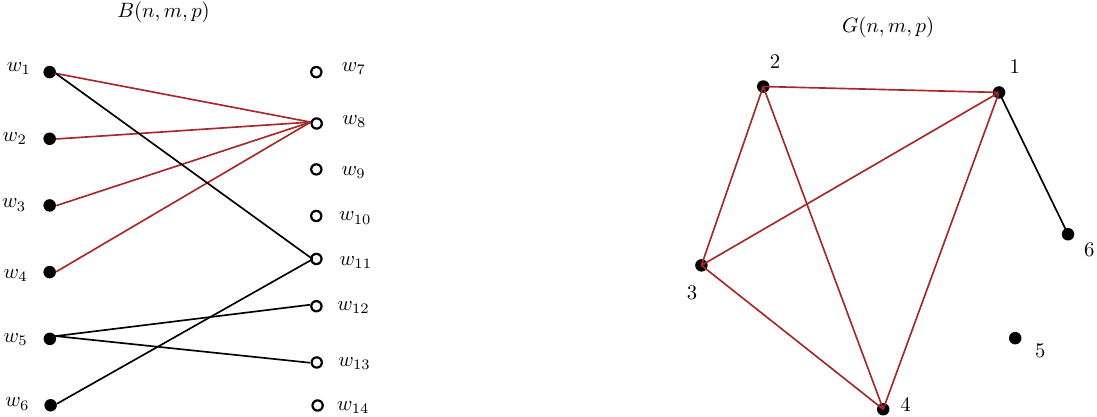}
\caption{\label{fig:1} On the left, an example of $B(n, m, p)$ with $n=6$, $m=8$. On the right, the induced intersection graph $G(n, m, p)$. Observe that the red edges and the adjacent vertices on the left, which depict a community with 4 members, induce a $K_{4}$ as subgraph (in red) on the right.}
\end{figure}

Let $n, m\in \mathbb N$ and $p\in [0, 1]$.  A random intersection graph with parameters $n, m, p$ can be obtained by first sampling a bipartite graph $B(n, m, p)$. More precisely, the vertex sets of this bipartite graph are denoted as  $\cV=\{w_{1}, w_{2}, \dots, w_{n}\}$ and $\cU=\{w_{n+1}, w_{n+2}, \dots, w_{n+m}\}$, of respective sizes $n$ and $m$. 
Among the $nm$ possible edges between $\cV$ and $\cU$, each is present in $B(n,m,p)$ independently with probability $p$. 
It is often useful to think of $\cV$-vertices (i.e.~elements of $\cV$) as individuals and $\cU$-vertices as communities. The graph $B(n, m, p)$  can be then interpreted as a network model where each individual joins a community independently with probability $p$. 
The {\it random intersection graph} $G(n, m, p)$ induced by $B(n,m,p)$ is the graph on the vertex set 
$[n]:=\{1, 2, 3, \dots, n\}$ where two vertices $i, j$ are adjacent to each other in $G(n,m,p)$ if and only if the individuals they represent, i.e.~$w_{i}$ and $w_{j}$, have joined a common community; see Fig.~\ref{fig:1} for an example. 
In particular, we note that individuals from the same community are always adjacent to each other in $G(n, m, p)$; in other words, a community in $B(n, m, p)$ translates into a ``clique'', i.e.~a copy of complete graph, inside $G(n, m, p)$. This suggests that the local structure of $G(n, m, p)$ is not ``tree-like'', unlike  many random graphs studied previously.  

Our main tool for studying $G(n, m, p)$ is a depth-first exploration of the bipartite graph $B(n, m, p)$, which allows us to extract a spanning forest from $B(n, m, p)$. Thanks to the criticality assumption, the question concerning the scaling limit of the graph largely reduces to identifying the scaling limit of this spanning forest. This, however, is still not a trivial question, as spanning forests that arise from the depth-first explorations of random graphs do not have ``nice'' distributional properties in general. 
Our approach, inspired by Conchon-Kerjan and Goldschmidt~\cite{CG23}, relies on representing the spanning forest of $B(n,m,p)$ as a sequence of tilted Bienaym\'e trees. 

The rest of the article is organised as follows. In Section~\ref{sec: param}, we introduce the relevant asymptotic regimes for the parameters of the model. In Section~\ref{sec: ER-intro}, we briefly review the result from \cite{ABBrGo12} concerning the convergence of critical Erd\H{o}s--R\'enyi graphs. The Erd\H{o}s--R\'enyi case serves as a valuable benchmark in our study.  We announce the main results in Section~\ref{sec: results}, whose proofs are found in Section~\ref{sec: proof}. 

\subsection{Phase transition and clustering regimes} 
\label{sec: param}
We are interested in how $G(n, m, p)$ behaves when both $n, m\to \infty$ and $p \sim \frac{c}{\sqrt{mn}}$ for some $c\in (0, \infty)$. Note in this case, the expected degree of any fixed vertex is given by 
\[
(n-1)\cdot \big(1-(1-p^2)^m\big) \to c^{2}, \quad \text{ as } n, m\to\infty. 
\]
It is shown in \cite{Beh2007, BST14} that the size of the largest connected component of $G(n, m, p)$ experiences a phase transition at $c=1$.  
In this work, we study the asymptotic behaviours of $G(n, m, p)$ in a regime around the criticality threshold $c=1$ called {\it critical window}.   
The analytic expression for the critical window of $G(n, m, p)$ depends, however, on the {\it clustering regimes}, a notion we are about to introduce. %we are in. Let us now explain this. 
We refer to the following quantity as the {\it clustering coefficient} of $G(n, m, p)$: 
\[
\CL(G(n, m, p))=\mathbb P\big(\text{vertex $2$ adjacent to $3$} \,|\, \text{both $2, 3$ adjacent to vertex $1$}\big).
\]
Clearly in above, $(1, 2, 3)$ can be replaced by any distinct triplet $(i, j, k)\in [n]^{3}$. In words, this is the probability of any three vertices forming a triangle knowing that they are already connected by edges. 
One can easily extend this definition of clustering coefficient to other graph models. In particular, for the Erd\H{o}s--R\'enyi model $G(n, p_{e})$, 
we note its clustering coefficient is always equal to $p_{e}$, since edges there are drawn independently. 
In contrast, for the random intersection graph $G(n, m, p)$, assuming $p\sim \frac{1}{\sqrt{mn}}$, we can distinguish three different outcomes:
\begin{itemize}
\item
If $\frac{m}{n}\to \infty$, then $\CL(G(n, m, p))\to 0$; this will be called the {\it light clustering regime}.
\item
If $\frac{m}{n}\to 0$, then $\CL(G(n, m, p))\to 1$; this will be called the {\it heavy clustering regime}.
\item
If $\frac{m}{n}\to \theta\in (0, \infty)$, then $\CL(G(n, m, p))\to (1+\sqrt{\theta})^{-1}\in (0, 1)$; this will be called the {\it moderate clustering regime}.
\end{itemize}
Critical window in this article refers to the following sets of assumptions on $(m, n, p)$:
%{\leqnomode
\begin{align}
\label{hyp: L}\tag{L}
&\text{either}\qquad \frac{m}{n}\to\infty\quad \text{and}\quad p = \frac{1}{\sqrt{mn}}(1+\lambda n^{-\frac13}) \text{ for some $\lambda\in \R$}; \\  \label{hyp: M} \tag{M}
&\text{or}\qquad\quad\;\,\exists\, \theta\in (0, \infty)\  \text{ s.t. } \ n^{\frac13}\big(\tfrac{m}{n}-\theta\big) \to 0\ \text{ and }\  p = \frac{1}{\sqrt{mn}}(1+\lambda n^{-\frac13})\text{ for some $\lambda\in \R$};\\ \label{hyp: H} \tag{H}
&\text{or}\qquad\quad\;\,\frac{m}{n}\to 0\quad \text{and} \quad p = \frac{1}{\sqrt{mn}}(1+\lambda m^{-\frac13})\text{ for some $\lambda\in \R$}.
\end{align}
Let us mention that  under~\eqref{hyp: L} and \eqref{hyp: H}, limit distributions for the component sizes of $G(n,m,p)$ have been investigated in Federico \cite{Federico19+}. 

\subsection{The continuum Erd\H{o}s--R\'enyi graph}
\label{sec: ER-intro}

Before stating our results on the convergence of critical random intersection graphs, let us give a brief account on what happens in $G(n, p_{e})$. 
Assume $p_{e}=\frac1n(1+2\lambda n^{-1/3})$ for some $\lambda\in \R$. (The reason for the unconventional constant $2$ will become clearer shortly.) 
Denote by $G^{n}_{k}$  the $k$-th largest connected component of $G(n, p_{e})$. Let $\dgr^{\mathrm{ER}}$ stand for the graph distance of $G(n, p_{e})$ and for $\alpha>0$, let $\alpha\cdot\mu^{\ER}_{k}$ be the finite measure which assigns a mass of $\alpha$ to each vertex of $G^{n}_{k}$.  
Equipped with the rescaled distance $n^{-1/3}\cdot \dgr^{\ER}$ and the finite measure $n^{-2/3}\cdot \mu^{\ER}_{k}$,  $(G^{n}_{k}, \dgr^{\ER}, \mu^{\ER}_{k})$ is an instance of (random) measured metric space. 
It is shown in \cite{ABBrGo12} that there exists a sequence of (random) measured metric spaces
\[
\mathcal G^{\mathrm{ER}}(\lambda)=\big\{\big(\cC^{\lambda, \infty}_{k}, d^{\lambda,\infty}_{k}, \mu^{\lambda,\infty}_{k}\big): k\ge 1\big\}
\]
so that the following convergence in distribution takes place as $n\to\infty$: 
\begin{equation}
\label{cv: ER-graph}
\Big\{\Big(G^{n}_{k}, \, n^{-\frac13}\cdot \dgr^{\mathrm{ER}}, \, n^{-\frac23}\cdot\mu^{\ER}_{k}\Big): k\ge 1\Big\} \;\Longrightarrow\; \mathcal G^{\mathrm{ER}}(\lambda)
\end{equation}
We will come back to the construction of $\mathcal G^{\mathrm{ER}}(\lambda)$ in Section~\ref{sec: ERgraph} and the precise meaning of the previous convergence in Section~\ref{sec: discuss}. 

\section{Main results}
\label{sec: results}

Recall \eqref{hyp: L}, \eqref{hyp: M}, \eqref{hyp: H}, the assumptions of critical window in the respective light, moderate and heavy clustering regimes. 
For $k\ge 1$, we denote by $C^{n, m}_{k}$ the $k$-th largest connected component of $G(n, m, p)$. Denote by $\dgr^{\mathrm{RIG}}$ the graph distance of $G(n, m, p)$. %There are two 
For $\alpha>0$, let $\alpha\cdot\mu^{\RIG}_{k}$ be the finite measure on $C\n_{k}$ which assigns a mass of $\alpha$ to each of its vertices.  

\begin{thm}[Scaling limit in the moderate clustering regime]
\label{thm2}
Under the assumption~\eqref{hyp: M}, 
there exists a sequence of (random) measured metric spaces $\mathcal G^{\mathrm{RIG}}(\lambda, \theta) = \{(\mathcal C^{\lambda,\theta}_{k}, d^{\lambda,\theta}_{k}, \mu^{\lambda,\theta}_{k}): k\ge 1\}$ so that 
 we have the following convergence in distribution as $n\to\infty$:
\begin{equation}
\label{eq: cv-thm2}
\Big\{\Big(C^{n, m}_{k}, \, n^{-\frac13}\cdot \dgr^{\mathrm{RIG}}, \, n^{-\frac23}\cdot\mu^{\RIG}_{k}\Big): k\ge 1\Big\} \;\Longrightarrow\; \mathcal G^{\mathrm{RIG}}(\lambda, \theta)
\end{equation}
with respect to the weak convergence of the product topology induced by the Gromov--Hausdorff--Prokhorov topology. 
\end{thm}

\begin{thm}[Scaling limit in the light clustering regime]
\label{thm1}
Under the assumption~\eqref{hyp: L},
we have the following convergence in distribution as $n\to\infty$:
\begin{equation}
\Big\{\Big(C^{n, m}_{k}, \, n^{-\frac13}\cdot \dgr^{\mathrm{RIG}}, \, n^{-\frac23}\cdot\mu^{\RIG}_{k}\Big): k\ge 1\Big\} \;\Longrightarrow\; \mathcal G^{\mathrm{ER}}(\lambda)
\end{equation}
with respect to the weak convergence of the product topology induced by the Gromov--Hausdorff--Prokhorov topology. 
\end{thm}

\begin{thm}[Scaling limit in the heavy clustering regime]
\label{thm3}
Under the assumption~\eqref{hyp: H}, we have the following convergence in distribution as $n\to\infty$:
\begin{equation}
\Big\{\Big(C^{n, m}_{k}, \, m^{-\frac13}\cdot \dgr^{\RIG}, \, m^{-\frac16}n^{-\frac12}\cdot \mu^{\RIG}_{k}\Big): k\ge 1\Big\} \;\Longrightarrow\; \mathcal G^{\mathrm{ER}}(\lambda)
\end{equation}
with respect to the weak convergence of the product topology induced by the Gromov--Hausdorff--Prokhorov topology. 
\end{thm}

To recap, in both light and heavy clustering regimes, we find the same limit object as in the critical Erd\H{o}s--R\'enyi graph. However, the diameters of $C\n_{k}$ scale differently in these two regimes: as $m=o(n)$ under \eqref{hyp: H}, Theorem \ref{thm3} shows that the diameters are much smaller than those under \eqref{hyp: L} or those of critical $G(n, p_{e})$. In the moderate clustering regime, the scaling limit of $C\n_{k}$ belongs to a two-parameter family  of measured metric spaces. As we will see in Section~\ref{sec: RIGgraph}, $\mathcal G^{\mathrm{RIG}}(\lambda, \theta)$ still bears a strong affinity with the continuum Erd\H{o}s--R\'enyi graph $\mathcal G^{\mathrm{ER}}(\lambda)$.

\medskip
Let us recall that $G(n, m, p)$ is a deterministic function of the bipartite graph $B(n, m, p)$. Furthermore, there is a natural bijection between the vertex set $[n]$ of $G(n, m, p)$ and the vertex set $\cV$ of $B(n, m, p)$ which maps $i\in [n]$ to $w_{i}\in \cV$. Denote this bijection as $\rho_{n}$. 
Let $\dgr^{\bi}$ stand for the graph distance of $B(n,m,p)$. 
Denote by $I_{k}$ the  smallest vertex label in $C\n_{k}$ and let $\hat C\n_{k}$ be the connected component of $B(n, m, p)$ containing $w_{I_{k}}$, $k\ge 1$.  
Let $\mu^{\bi}_{k}$ be the measure of $\hat C\n_{k}$ that assigns a unit mass to each $\cV$-vertex in $\hat C\n_{k}$. 
We write $\dghp$ for the Gromov--Hausdorff--Prokhorov distance between two measured metric spaces; see Section~\ref{sec: ghp} for a formal definition. 
Let us point out the previous results on the scaling limit of $G(n, m, p)$ translates into analogous results for the bipartite graph $B(n, m, p)$, thanks to the following observation. 

\begin{prop}
\label{prop: approx-gr}
For all $n, m\in \N$, $\rho_{n}$ is an isometry from $([n], 2\cdot\dgr^{\RIG})$ to $(\cV, \dgr^{\bi})$. As a result, the following holds almost surely for each $k\ge 1$: 
\[
\dghp\Big(\big(\hat C\n_{k}, \dgr^{\bi}, \mu^{\bi}_{k}\big), \big(C\n_{k}, 2\cdot \dgr^{\RIG}, \mu^{\RIG}_{k}\big)\Big)\le 1.
\]
\end{prop}

\subsection{Convergence of the graph exploration processes}
\label{sec: intro-proc}

We obtain the aforementioned results on the scaling limit of $G(n, m, p)$ by studying certain 
stochastic processes that arise from a depth-first exploration of the bipartite graph $B(n, m, p)$. A related exploration of $G(n, m, p)$ appeared in~\cite{Federico19+}. Noting  Proposition~\ref{prop: approx-gr}, which tells us that $B(n, m, p)$ will offer the same scaling limit as $G(n, m, p)$, we prefer to run the exploration on the bipartite graph as it contains additional information useful for triangle counts. 
We will use the following notation: for a vertex $w\in \cU\cup\cV$ of $B(n, m, p)$, we denote by
\[
\cB(w) = \{w'\in \cU\cup \cV: \dgr^{\bi}(w', w)=1\}, 
\]
the neighbourhood of $w$ in the graph. For a subgraph $G$ of $B(n, m, p)$, we write $\cV(G)$ (resp.~$\cU(G)$) for the set of $\cV$-vertices (resp.~$\cU$-vertices) contained in $G$; we also denote by $\cE(G)$ the edge set of $G$.  
The following algorithm explores the neighbourhood of each vertex in a ``depth-first'' fashion,  and outputs a spanning forest $\cF$ of $B(n, m, p)$ in the end. To do so, the algorithm keeps track of the vertices that have not been explored at step $k$  using two subsets $\cV_{k}$ and $\cU_{k}$. The vertex explored at step $k$ is denoted as $v_{k}$, which is always a $\cV$-vertex. 
There is a further subset of $\cV$-vertices that do not belong to $\cV_{k}$ nor to $\{v_{1}, \dots, v_{k}\}$: these are the neighbours of the explored vertices whose neighbourhoods are yet to be explored; such a vertex is said to be {\it active} at step $k$. 
Information about active vertices are stored in an ordered list (i.e.~sequence) $\cA_{k}$, and the ordering in $\cA_{k}$ determines which vertex to be explored next.

\begin{alg}[Depth-first exploration of $B(n, m, p)$]
\label{alg}
Initially, set $\cF_{0}$ to be the empty graph and $\cA_{0}=\cA_{0}^{\ast}$ to be the empty sequence. 
Set $\cV_{0}^{\ast}=\cV$ and $\cU_{0}=\cU$. 

\noi
{\bf Step $k\ge 1$.}
If $\cA^{\ast}_{k-1}=\varnothing$, let $v_{k}$ be a uniform element from $\cV^{\ast}_{k-1}$; otherwise, let $v_{k}$ be the first element in $\cA^{\ast}_{k-1}$. Set $\cV_{k-1}=\cV_{k-1}^{\ast}\setminus\{v_{k}\}$, and let $\cA_{k-1}$ be the sequence obtained from $\cA^{\ast}_{k-1}$ by removing $v_{k}$. 
Let 
\[
\cM_{k}=\mathcal B(v_{k})\cap \mathcal U_{k-1} \quad\text{ and }\quad X\n_{k}=\#\cM_{k}. 
\]
If $\cM_{k}\ne\varnothing$, sort its elements in an increasing order of their labels and denote by $(u_{k, i})_{1\le i\le X\n_{k}}$ the sorted list. For $1\le i\le X\n_{k}$, define
\[
\cN_{k, i} = \cB(u_{k, i})\cap \cV_{k-1}\setminus \bigcup_{j=1}^{i-1}\cN_{k, j} \quad \text{and} \quad \cN_{k}=\bigcup_{1\le i\le X\n_{k}}\cN_{k, i}. 
\]
Let $\cF_{k}$ be the bipartite graph on the vertex sets $\cV(\cF_{k})=\cV(\cF_{k-1})\cup\{v_{k}\}\cup\cN_{k}$ and $\cU(\cF_{k})=\cU(\cF_{k-1})\cup \cM_{k}$ with its edge set given by 
\[
\cE(\cF_{k})= \cE(\cF_{k-1})\bigcup\big\{\{v_{k}, u_{k, i}\}: 1\le i\le X\n_{k}\big\}\bigcup \bigcup_{1\le i\le X\n_{k}}\big\{\{u_{k, i}, v\}: v\in \cN_{k, i}\big\}.
\]
Set $\cV_{k}^{\ast}=\cV\setminus\cV(\cF_{k})=\cV_{k-1}\setminus \cN_{k}$ and $\cU_{k}=\cU\setminus\cU(\cF_{k})=\cU_{k-1}\setminus \cM_{k}$. If $\cN_{k}$ is nonempty, sort its elements in an increasing order of their labels and add the ordered list to the top of $\cA_{k-1}$; call the new sequence $\cA^{\ast}_{k}$. 

\noi
{\bf Stop when $\cV_{k}^{\ast}=\varnothing$.} Add each member of $\cU_{k}$ to $\cF_{k}$ as an isolated vertex; call the resulting graph $\cF$. 
\end{alg}

We note that the algorithm always terminates at step $n$, when all the $\cV$-vertices are explored. Moreover, for each $1\le k\le n$, the sets $\cN_{k, i}, 1\le i\le X\n_{k}$, are pairwise disjoint. Note also that $\cF_{k}$ is loop-free, and consequently the final graph $\cF$ is a spanning forest of $B(n, m, p)$. In particular, each connected component of $\cF$ is a spanning tree of a connected component of $B(n, m, p)$. If a connected component of $\cF$ contains more than one vertex, then it must contain some $v_{k}$ and we root the tree at the $v_{k}$ with the smallest $k$; otherwise, we root the tree at the sole vertex. 

\begin{figure}[tp]
\centering
\begin{minipage}{.5\textwidth}
  \centering
  \includegraphics[height = 5cm]{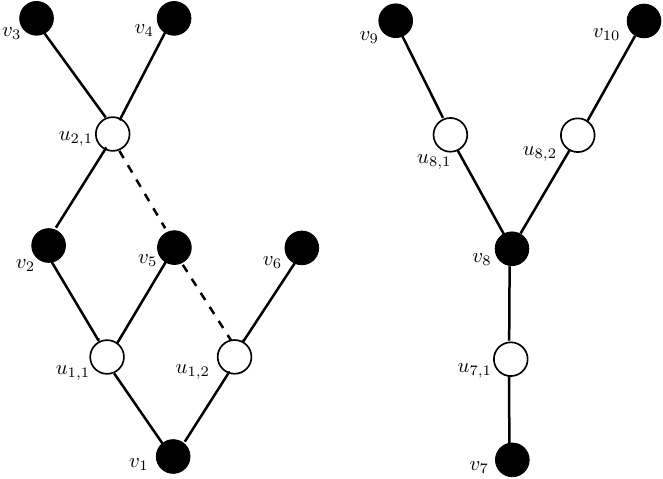}
  
  \vspace{1cm}
  \includegraphics[height = 3.4cm]{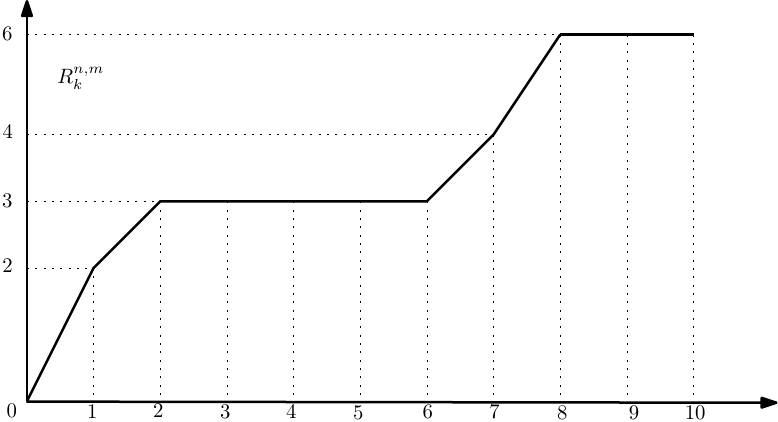}
  %\captionof{figure}{A figure}
  %\label{fig:test1}
\end{minipage}%
\begin{minipage}{.5\textwidth}
  \centering
  \includegraphics[height = 4.3cm]{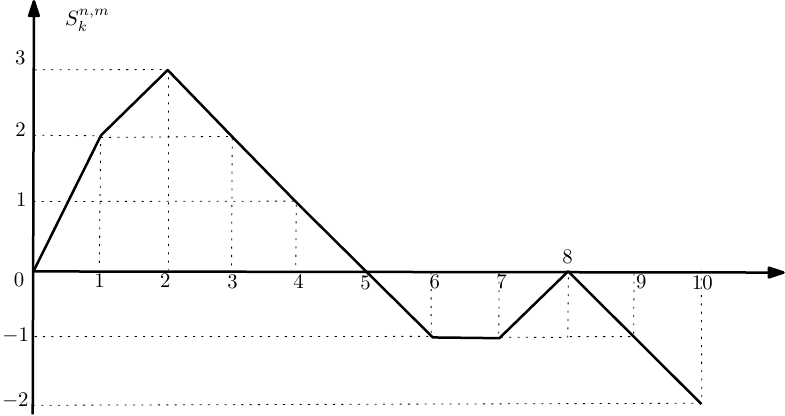}
  
  \vspace{2cm}
  \includegraphics[height = 3cm]{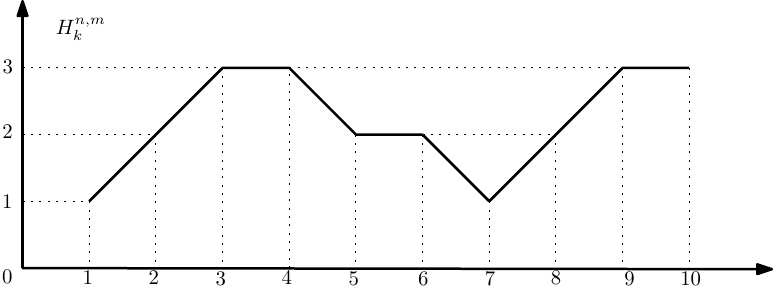}
 \end{minipage}
\caption{\label{fig:2}Upper left: an example of $B(n, m, p)$ and its exploration; dashed edges are those not in the spanning forest. Upper right:  the corresponding $S\n$. Lower left: the corresponding $R\n$. Lower right: the corresponding $H\n$. The successive values of $\#\cN_{k}$ are $3, 2, 0, 0, 0, 0, 1, 2, 0, 0$; the values of $\#\cM_{k}$ are resp.~$2, 1, 0, 0, 0, 0, 1, 2, 0, 0$. The values of $\#\cA_{k}$ are resp.~$2, 3, 2, 1, 0, 0, 0, 1, 0, 0$. }
\end{figure}

The following stochastic processes are the key tools in our study of $B(n, m, p)$ and $G(n, m, p)$. Set $R\n_{0}=S\n_{0}=0$. For $k\ge 1$, let 
\begin{equation}
\label{def: dps}
R\n_{k}-R\n_{k-1} =\#\cM_{k}, \quad  S\n_{k}-S\n_{k-1}=\#\mathcal N_{k}-1.
\end{equation}
Meanwhile, for a sequence of integers $(s_j)_{0\le j\le k}$, we define
\begin{equation}
\label{def: ht-proc}
\mathscr H_{k}\big((s_j)_{0\le j\le k}\big)=\#\Big\{0\le j\le k-1: s_j=\min_{j\le i\le k-1}s_i\Big\}. 
\end{equation}
We then set
\[
H\n_{k} = \mathscr H_{k}((S\n_j)_{0\le j\le k}), \quad k\ge 1.
\]
We say the {\it height} of a vertex $v$ in $\cF$ is the graph distance from $v$ to the root of the tree component of $\cF$ that contains $v$. 
The usefulness of $H\n:=(H\n_{k})_{k\ge 1}$ lies in the following property, whose proof is found towards the end of Section~\ref{sec: pf}. 
\begin{lem}
\label{lem: ht}
For $n\in \N$ and $1\le k\le n$, $2H\n_{k}-2$ is the height of $v_{k}$ in $\cF$. 
\end{lem}
We will refer to the pair $(R\n, S\n):=(R\n_{k}, S\n_{k})_{k\ge 0}$ as the {\it depth-first walk} of $B(n,m, p)$, and $(H\n_{k})_{k\ge 0}$ as the {\it height process}; see Fig.~\ref{fig:2} for an example. 
Unlike usual height processes, our definition of $H\n$ only pertains to vertices at even heights. This is nevertheless sufficient for our purpose as the rest of $\cF$ are at distance 1 from those vertices. 
The triplet $(R\n, S\n, H\n)$ forms the graph exploration process for $B(n, m, p)$ and their scaling limit is our next focus.  

Let $(W_{t})_{t\ge 0}, (W^{\ast}_{t})_{t\ge 0}$ be two independent standard linear Brownian motions. Denote $\R_{+}=[0,\infty)$ and let $d\in\N$. 
For  a closed interval $I\subseteq \R_{+}$ and $d\in \N$, let $\mathbb D^{\ast}(I, \R^{d})$ be the set of all c\`adl\`ag (i.e.~right-continuous with left-hand limits) maps from $I$ to $\R^{d}$, equipped with the {\it uniform} topology on all compact sets of $I$. 

\begin{prop}[Convergence of the exploration processes in the light clustering regime]
\label{thm: cv-light}
Under the assumption \eqref{hyp: L}, we have the following convergence in distribution 
in $\mathbb D^{\ast}(\R_+, \R^{3})$:
\begin{equation}
\label{cv: light}
\Big\{n^{-\frac16}m^{-\frac12}R\n_{\lf n^{2/3}t\rf} , n^{-\frac13}S\n_{\lfloor n^{2/3}t\rfloor}, n^{-\frac13}H\n_{\lfloor n^{2/3}t\rfloor}: t\ge 0\Big\} \Longrightarrow \Big\{ t, \mathcal S^{\lambda, \infty}_{t}, \mathcal H^{\lambda,\infty}_{t}: t\ge 0\Big\},
\end{equation}
where the limit processes $\cS^{\lambda,\infty}=(\cS^{\lambda, \infty}_{t})_{t\ge 0}, \cH^{\lambda,\infty}=(\cH^{\lambda, \infty}_{t})_{t\ge 0}$ are defined as follows: 
\[
\mathcal S^{\lambda, \infty}_{t} = W_{t} + 2\lambda t - \tfrac12 t^{2},  \quad \mathcal H^{\lambda, \infty}_{t} = 2\Big(\mathcal S^{\lambda, \infty}_{t}-\min_{s\le t}\mathcal S^{\lambda, \infty}_{s}\Big), \quad t\ge 0. 
\]
\end{prop}

\begin{prop}[Convergence of the exploration processes in the moderate clustering regime]
\label{thm: cv-moderate}
Assume \eqref{hyp: M}. Let $\tilde R\n_{k} = R\n_{k}-\sqrt\theta k$, $k\ge 1$. 
We have the following convergence in distribution in $\mathbb D^{\ast}(\R_+, \R^{3})$:
\begin{equation}
\label{cv: moderate}
\Big\{n^{-\frac13}\tilde R\n_{\lfloor n^{2/3}t\rfloor}, n^{-\frac13}S\n_{\lf n^{2/3}t\rf}, n^{-\frac13}H\n_{\lfloor n^{2/3}t\rfloor}: t\ge 0\Big\} \Longrightarrow \Big\{  \cR^{\lambda, \theta}_{t}, \cS^{\lambda, \theta}_{t}, \cH^{\lambda, \theta}_{t}: t\ge 0\Big\},
\end{equation}
where the limit processes $(\cR^{\lambda,\theta},\cS^{\lambda,\theta})=(\cR^{\lambda, \theta}_{t}, \cS^{\lambda, \theta}_{t})_{t\ge 0}, \cH^{\lambda,\theta}=(\cH^{\lambda, \theta}_{t})_{t\ge 0}$ are defined as follows: 
\begin{align}\label{def: W-proc}
& \cR^{\lambda, \theta}_{t} = \theta^{\frac14} W^{\ast}_{t} + \theta^{\frac12}\lambda t-\tfrac{1}{2}t^{2}, \quad 
\cS^{\lambda, \theta}_{t}   = W_{t}+\theta^{-\frac14}W^{\ast}_{t}+2\lambda t-\tfrac12(1+\theta^{-\frac12})t^{2}, \ \text{ and } \\ \label{def: ht}
& \cH^{\lambda, \theta}_{t}  = \frac{2}{1+\theta^{-\frac12}}\Big(\cS^{\lambda, \theta}_{t} - \min_{s\le t}\cS^{\lambda, \theta}_{s}\Big). 
\end{align}
\end{prop}
The process $(\cH^{\lambda,\infty}_{t})_{t\ge 0}$ (resp.~$(\cH^{\lambda,\theta}_{t})_{t\ge 0}$) is the {\it height process} of $(\cS^{\lambda, \infty}_{t})_{t\ge 0}$ (resp.~of $(\cS^{\lambda,\theta}_{t})_{t\ge 0}$). We point to Sections~\ref{sec: ERgraph} and~\ref{sec: RIGgraph} for further information on their definitions. 
Let us also note that by letting $\theta\to\infty$ in the expression of $\cS^{\lambda, \theta}_{t}$, we recover $\cS^{\lambda, \infty}_{t}$, which motivates our notation. 

\medskip
From the definition of $B(n, m, p)$, it is clear that exchanging the roles of $n$ and $m$ leaves the distribution of the graph unchanged. By running Algorithm~\ref{alg} on $B(m, n, p)$ instead of $B(n,m,p)$ and letting $(R\m_{k}, S\m_{k})_{k\ge 0}, (H\m_{k})_{k\ge 0}$ be the corresponding depth-first walk and height process, Proposition~\ref{thm: cv-light} has an immediate consequence as follows.

\begin{cor}[Convergence of the exploration processes in the heavy clustering regime]
Under Assumption \eqref{hyp: H}, the following convergence in distribution holds in $\mathbb D^{\ast}(\R_+, \R^{3})$:
\begin{equation}
\label{cv: heavy}
\Big\{m^{-\frac16}n^{-\frac12}R\m_{\lf m^{2/3}t\rf} , m^{-\frac13}S\m_{\lfloor m^{2/3}t\rfloor}, m^{-\frac13}H\m_{\lfloor m^{2/3}t\rfloor}: t\ge 0\Big\} \Longrightarrow \Big\{ t, \mathcal S^{\lambda, \infty}_{t}, \mathcal H^{\lambda, \infty}_{t}: t\ge 0\Big\},
\end{equation}
where $\cS^{\lambda, \infty}_{t}, \cH^{\lambda, \infty}_{t}$ are the same as in Proposition~\ref{thm: cv-light}. 
\end{cor} 

Proof of Propositions~\ref{thm: cv-light} and~\ref{thm: cv-moderate} will occupy a major part of Section~\ref{sec: proof}. 

\subsection{Triangle counts and surplus edges}
\label{sec: sps}

Given a connected graph $G$ and a spanning tree $T$ of $G$, we say an edge $e$ of $G$ is a surplus edge if it is not present in $T$. Although the notion of surplus edge depends on the spanning tree, the number of surplus edges does not and is given by $\# E(G)-\#G +1$, where $\# E(G)$ stands for the number of edges in $G$ and $\#G$ counts the number of vertices in $G$. In a nutshell, results here tell us that large connected components of the critical bipartite graph $B(n, m, p)$ contain very few surplus edges even as their sizes tend to infinity, which is reminiscent of the situation for critical Erd\H{o}s--R\'enyi graphs (\cite{Al97, ABBrGo12}). In contrast, as we will see in Theorems~\ref{thm: triangle-moderate}-\ref{thm: triangle-heavy} below, $G(n,m,p)$ can exhibit a range of  behaviours regarding surplus edges. Nevertheless, the presence of these surplus edges in $G(n, m, p)$ does not change significantly the metric aspect of the graph, as indicated in Proposition~\ref{prop: approx-gr}. 

Let us start with the bipartite graph $B(n, m, p)$. From now on, we will refer to an edge of $B(n, m, p)$ as a (bipartite) surplus edge if it is not contained in the spanning forest $\cF$ from Algorithm~\ref{alg}. We will need information on the distribution of these bipartite surplus edges to identify the scaling limit of $B(n,m,p)$. 
Let us denote by $\cE^{\rS}(n, m)$ the set of surplus edges of $B(n, m, p)$. Recall from Algorithm~\ref{alg} the sets $\cM_{k}, \cA_{k}, \cN_{k, i}$. We have the following description of $\cE^{\rS}(n, m)$. 

\begin{lem}
\label{lem: surplus}
If $e=\{u, w\}\in \cE^{\rS}(n, m)$, then there is a unique $k\in [n]$ so that one of the two following cases occurs:
\begin{enumerate}[(i)]
\item
either $u\in \cM_{k}$ and $w\in \cA_{k-1}$;
\item
or there exists $1\le i<j\le \#\cM_{k}$ satisfying $w\in \cN_{k, i}$ and $u=u_{k, j}$.
\end{enumerate}
\end{lem}

Proof for Lemma~\ref{lem: surplus} is short; we therefore include it here. See also the example illustrated in Fig.~\ref{fig:2}. 

\begin{proof}
Since $\cM_{k}, k\in [n]$, are disjoint and $\cup_{k\in [n]}\cM_{k}$ comprises all non isolated $\cU$-vertices,  there must be a unique $k\in \N$ satisfying $u\in \cM_{k}$. In that case, we must have the neighbourhood $\cB(u)$ of $u$ contained in $\cN_{k}\cup\cA_{k-1}\cup\{v_{k}\}$. However, if $w=v_{k}$, then $e=\{u, w\}$ would be an edge of $\cF$. So $w$ belongs to either $\cA_{k-1}$ or $\cN_{k}$. In the latter case, as the edge $e$ is not present in $\cF$, it can only happen when $w$ is adjacent to at least two elements of $\cM_{k}$. 
\end{proof} 

To study the limit distribution of $\cE^{\rS}(n, m)$, we first turn it into a point measure as follows. For $e=\{u, w\}\in \cE^{\rS}(n, m)$,  let $k(e)$ be the unique $k\in [n]$ identified in Lemma~\ref{lem: surplus}. If $w\in \cA_{k-1}$, recalling that $\cA_{k-1}$ is an ordered list, 
we let $l(e)$ be the rank of the element $w$ in $\cA_{k-1}$. If $w\notin \cA_{k-1}$, we set $l(e)=0$. We introduce the following point measure on $\R^{2}$: 
\begin{equation}
\label{def: Pn}
\tilde P\n(dx, dy) = \sum_{e\in\cE^{\rS}(n,m)}\delta_{(n^{-2/3}k(e), \,n^{-1/3}l(e))}.
\end{equation}
In the example shown in Fig.~\ref{fig:2}, there are two surplus edges: $\{u_{1,2}, v_5\}$, which is of type (ii), and $\{u_{2, 1}, v_5\}$, of type (i). In this case, $\tilde P\n$ is the point measure containing two atoms respectively at $(n^{-2/3}\cdot 1, 0)$ and $(n^{-2/3}\cdot 2, n^{-1/3}\cdot 1)$, since $v_5$ sits at the top of $\cA_1$.   
We note that $\tilde P\n$ is not necessarily simple, i.e.~it may contain multiple atoms at the same location. We let $P\n$ be the simple point measure with the same support as $\tilde P\n$. 
Given $\cF$ and $P\n$, we will construct in Section~\ref{sec: cv-ghp} a graph that has the same scaling limit as $B(n, m, p)$.  To describe the limit of $P\n$, we need the following objects. Recall from Propositions~\ref{thm: cv-light} and~\ref{thm: cv-moderate} the processes $\cS^{\lambda, \infty}$ and $(\cR^{\lambda,\theta},\cS^{\lambda,\theta})$. 
Conditional on $\cS^{\lambda, \infty}$, let $\mathcal P^{\lambda, \infty}$ be a Poisson point measure on $\R^{2}$ with intensity measure $dx dy \mathbf 1_{\{(x, y)\in \cD^{\lambda,\infty}\}}$, where 
\[
\cD^{\lambda,\infty}=\big\{(x, y) : x\in\R_{+}, 0\le y\le \cS^{\lambda, \infty}_{x}-\min_{u\le x}\cS^{\lambda, \infty}_{u} \big\}.
\]
Similarly, conditional on $(\cR^{\lambda,\theta}, \cS^{\lambda,\theta})$, let $\mathcal P^{\lambda,\theta}$ be a Poisson point measure on $\R^{2}$ with intensity measure $dx dy \mathbf 1_{\{(x, y)\in \cD^{\lambda,\theta}\}}$, where 
\[
\cD^{\lambda,\theta}=\big\{(x, y) : x\in\R_{+}, 0\le y\le \cS^{\lambda, \theta}_{x}-\min_{u\le x}\cS^{\lambda, \theta}_{u} \big\}.
\]

\begin{prop}[Convergence of the bipartite surplus edges]
\label{prop: surplus-bi}
Under the assumption \eqref{hyp: L},  jointly with the convergence \eqref{cv: light} in Proposition \ref{thm: cv-light}, we have the weak convergence of the random measure $P\n$ to $\mathcal P^{\lambda, \infty}$ on every compact set of $\R^{2}$. Under the assumption \eqref{hyp: M},  jointly with the convergence  \eqref{cv: moderate} in Proposition \ref{thm: cv-moderate}, we have the weak convergence of the random measure $P\n$ to $\mathcal P^{\lambda,\theta}$ on every compact set of $\R^{2}$. 
\end{prop}

Proof of Proposition~\ref{prop: surplus-bi} is given in Section~\ref{sec: surplus}. 
We now turn our attention to the random intersection graph $G(n, m, p)$. %
Let us denote by $\rL\n_{k}$ the number of triangles in the $k$-th largest connected component of $G(n, m, p)$, noting that each triangle must contain at least one surplus edge. 
Recall that the excursions of $\cS^{\lambda,\infty}$ above its running infimum can be ranked in a decreasing order of their lengths, as shown by Aldous in~\cite{Al97}. This also holds true for $\cS^{\lambda,\theta}$; see Section~\ref{sec: RIGgraph}. Hence, for $\theta\in (0,\infty)\cup\{\infty\}$ and $k\in\N$, we can define
\begin{equation}
\label{def: zeta}
\zeta^{\lambda,\theta}_{k} = \text{length of the $k$-th longest excursion of } \big(\cS^{\lambda,\theta}_{t}-\min_{s\le t}\cS^{\lambda,\theta}_{s}\big)_{t\ge 0}. 
\end{equation}

\begin{thm}[Triangle counts in the moderate clustering regime]
\label{thm: triangle-moderate}
Under the assumption \eqref{hyp: M},  for each $k\ge 1$, we have the following joint convergence in distribution
\[
n^{-\frac23}\cdot \big(\rL\n_{1}, \rL\n_{2}, \rL\n_{3}, \cdots, \rL\n_{k}\big) \ \Longrightarrow\  c_{\theta}\cdot \big(\zeta^{\lambda,\theta}_{1}, \zeta^{\lambda,\theta}_{2}, \zeta^{\lambda,\theta}_{3}, \cdots, \zeta^{\lambda,\theta}_{k}\big). 
\]
The constant $c_{\theta}$ is given by 
\begin{equation}
\label{def: ctheta}
c_{\theta} = \frac{1}{2\sqrt{\theta}}+\frac{1}{6\theta} = \frac{\sqrt{\theta}}{6}\mathbb E\Big[(Y_{\theta}+1)Y_{\theta}(Y_{\theta}-1)\Big], 
\end{equation}
where $Y_{\theta}$ stands for a Poisson random variable of expectation $\theta^{-1/2}$. 
\end{thm}

\begin{thm}[Triangle counts in the light clustering regime]
\label{thm: triangle-light}
Under the assumption \eqref{hyp: L}, the following statements hold true.
\begin{enumerate}[(i)]
\item
If $\frac{m}{n^{7/3}}\to 0$ as $n\to\infty$, then for each $k\ge 1$, we have the following joint convergence in distribution 
\[
m^{\frac12}n^{-\frac76}\cdot \big(\rL\n_{1}, \rL\n_{2}, \rL\n_{3}, \cdots, \rL\n_{k}\big) \ \Longrightarrow\  \big(\tfrac12\zeta^{\lambda,\infty}_{1}, \tfrac12\zeta^{\lambda,\infty}_{2}, \tfrac12\zeta^{\lambda,\infty}_{3}, \cdots, \tfrac12\zeta^{\lambda,\infty}_{k}\big).
\]
\item
If $\frac{m}{n^{7/3}}\to \infty$ as $n\to\infty$, then for each $k\ge 1$, we have
\[
\max_{j\le k}\rL\n_{j}\longrightarrow 0 \quad \text{ in probability}.
\]
\item
If $m=\lf n^{7/3}\rf$, then for each $k\ge 1$, we have the joint convergence in distribution: 
\[
\big(\rL\n_{1}, \rL\n_{2}, \rL\n_{3}, \cdots, \rL\n_{k}\big) \ \Longrightarrow\  \big(\cL_{1},  \cL_{2}, \cL_{3}, \cdots, \cL_{k}\big),
\]
where conditional on $\cS^{\lambda,\infty}$, $(\cL_{i})_{1\le i\le k}$ is a collection of independent Poisson variables with respective expectations $\tfrac12\zeta^{\lambda,\infty}_{i}$, $1\le i\le k$. 
\end{enumerate}
\end{thm}

\begin{thm}[Triangle counts in the heavy clustering regime]
\label{thm: triangle-heavy}
Under the assumption \eqref{hyp: H},  for each $k\ge 1$, we have the following joint convergence in distribution
\[
m^{\frac56}n^{-\frac32}\cdot \big(\rL\n_{1}, \rL\n_{2}, \rL\n_{3}, \cdots, \rL\n_{k}\big) \ \Longrightarrow\   \big(\tfrac16\zeta^{\lambda,\infty}_{1}, \tfrac16\zeta^{\lambda,\infty}_{2}, \tfrac16\zeta^{\lambda,\infty}_{3}, \cdots, \tfrac16\zeta^{\lambda,\infty}_{k}\big).
\]
%where  $\zeta^{\lambda,\infty}_{i}$ is as in Theorem~\ref{thm: triangle-light}. 
\end{thm}

Proof of Theorems~\ref{thm: triangle-moderate} and~\ref{thm: triangle-light} is found in Section~\ref{sec: triangle} and proof of Theorem~\ref{thm: triangle-heavy} in Section~\ref{sec: heavy}. 

\subsection{Construction of the limit graphs}
\label{sec: scale}

After a quick recap on some topological notions, 
we explain here the constructions of $\cG^{\ER}(\lambda)$ and $\cG^{\RIG}(\lambda,\theta)$ from the Brownian motion. 

\subsubsection{Graphs encoded by real-valued functions}
\label{sec: ghp}
We follow the approaches in \cite{ABBrGo12, BrDuWa21}, with some minor differences in detail. 
The account we give here relies on the notions of real trees \cite{Evans}, measured metric spaces \cite{Gromov}, and Gromov--Hausdorff--Prokhorov topology \cite{Mi09}, which are now standard. We therefore point to the previous references for their definitions and further background. 
We first recall the construction of a real tree from a real-valued function. Let $\zeta_{h}\in (0, \infty)$ and suppose that $h: [0, \zeta_{h}] \to \R_{+}$ is a continuous function. The following symmetric function 
\[
d_{h}(s, t) := h(s)+h(t) - 2 b_{h}(s, t), \quad \text{where} \quad b_{h}(s, t) = \min\{h(u): \min(s, t)\le u\le \max(s, t)\}
\]
defines a pseudo-distance on $[0, \zeta_{h}]$. To turn it into a distance, we introduce an equivalence relation by defining $s\sim_{h} t$ if and only if $d_{h}(s, t)=0$. Then $d_{h}$ induces a distance on the quotient space $T_{h}:=[0, \zeta_{h}]/\sim_{h}$, which we still denote as $d_{h}$. 
Moreover, the compact metric space $(T_{h}, d_{h})$ has the property that every pair of points in it is joined by a unique path which is also geodesic. In other words, $(T_{h}, d_{h})$ is a real tree. Let  $p_{h}$ stand for the canonical projection from $[0, \zeta_{h}]$ to $T_{h}$ and denote by $\mu_{h}$ the push-forward of the Lebesgue measure on $[0, \zeta_{h}]$ by $p_{h}$. Set
\begin{equation}
\label{def: real-tree}
\cT_{h}:= \big(T_{h}, d_{h}, \mu_{h}\big). 
\end{equation}
For two points $x, y$ from $T_{h}$, let us write $\llbracket x, y\rrbracket$ for the unique path in $T_{h}$ from $x$ to $y$, whose length is given by $d_{h}(x, y)$. Suppose that $q\in \N$ and $\Pi=\{(s_{i}, t_{i}): 1\le i\le q\}$ is a collection satisfying $0\le s_{i}\le t_{i}\le \zeta_{h}$ for $1\le i\le q$. 
Their images in $T_{h}$ are denoted as
\[
x_{i}=p_{h}(s_{i}), \quad y_{i}= p_{h}(t_{i}), \quad 1\le i\le q. 
\]
For a pair of $(x, y)\in T_{h}^{2}$, we define a set $\Gamma_{\Pi}^{x, y}$ formed by sequences of paths that take the following form: 
\[
\boldsymbol\gamma=\big(\llbracket u_{1}, u_{2}\rrbracket, \llbracket u_{3}, u_{4}\rrbracket, \dots, \llbracket u_{2\ell-1}, u_{2\ell}\rrbracket\big)\,, 
\]
where $\ell\in \N$,  $u_{1}=x, u_{2\ell}=y$, $u_{2j}\in \{x_{i}: 1\le i\le q\}$, and if $u_{2j}=x_{i}$, then $u_{2j+1}=y_{i}$, for $1\le j\le \ell-1$. 
For such a $\boldsymbol\gamma=\big(\llbracket u_{2j-1}, u_{2j} \rrbracket\big)_{1\le j\le \ell}$, we define its {\it $\Pi$-modified length} as follows:
\[
\ell_{\Pi}(\boldsymbol\gamma) = \sum_{j=1}^{\ell} d_{h}(u_{2j-1}, u_{2j}).
\]
The following defines a pseudo-distance on $T_{h}$: for all $x, y\in T_{h}$, let
\[
d_{h, \Pi}(x, y) = \inf\big\{\ell_{\Pi}(\boldsymbol\gamma): \boldsymbol\gamma\in \Gamma_{\Pi}^{x, y}\big\}.
\]
As previously, we can turn it into a true distance by quotienting the points at $d_{h, \Pi}$-distance 0 from each other. Call the resulting metric space $(G_{h, \Pi}, d_{h, \Pi})$ and denote by $p_{h, \Pi}$ the canonical projection from $(T_{h}, d_{h})$ to $(G_{h, \Pi}, d_{h, \Pi})$. Write $\mu_{h, \Pi}$ for the push-forward of $\mu_{h}$ by $p_{h,\Pi}$. Finally, let us denote
\begin{equation}
\label{def: graph}
\cG(h, \Pi) = \big( G_{h, \Pi}, d_{h, \Pi},  \mu_{h, \Pi}\big), 
\end{equation}
which is a measured metric space. 
Recall that two measured metric spaces, say $\cG_{1}=(G_{1}, d_{1}, \mu_{1})$ and $\cG_{2}=(G_{2}, d_{2}, \mu_{2})$, can be compared using their Gromov--Hausdorff--Prokhorov distance: 
\[
\dghp(\cG_{1}, \cG_{2}) := \inf \big\{ \dhaus\big(\phi_{1}(G_{1}), \phi_{2}(G_{2})\big) + \dpr\big(\mu_{1}\circ \phi_{1}^{-1}, \mu_{2}\circ \phi_{2}^{-1}\big)\big\},
\]
where the infimum is over all Polish spaces $(E, d_{E})$ and isometric embeddings $\phi_{i}$ from $G_{i}\to E$; $\dhaus$ stands for the Hausdorff distance of $E$, and $\dpr$ is the Prokhorov distance for the finite Borel measures on $E$, with $\mu_{i}\circ\phi_{i}^{-1}$ standing for the push-forward of $\mu_{i}$ by $\phi_{i}$, $i=1, 2$. 
We further point out that the space of (equivalence classes) of compact measured metric spaces is a Polish space in the topology induced by $\dghp$ (\cite{Mi09}). 
The following statement, taken from \cite{BrDuWa21}, provides a practical way to prove Gromov--Hausdorff--Prokhorov convergence for measured metric spaces that are constructed from real-valued functions as described above. Let $q\in \N$.  For $j\in \{1, 2\}$, suppose that $h_j: [0, \zeta_{j}]\to \R_{+}$ is a continuous function and $\Pi_{j}=\{(s_{j, i}, t_{j, i}): 1\le i\le q\}$ is a collection of points with $0\le s_{j, i}\le t_{j, i}<\zeta_{j}$ for each $j\le q$. Suppose further that there is some $\delta>0$ verifying
\[
\max_{1\le i\le q}|s_{1, i}-s_{2, i}| \le \delta, \quad \max_{1\le i\le q}|t_{1, i}-t_{2, i}| \le \delta. 
\]
Then we have (\cite{BrDuWa21}, Lemma 2.7)
\begin{equation}
\label{eq: ghp}
\dghp\big(\cG(h_{1}, \Pi_{1}), \cG(h_{2}, \Pi_{2})\big) \le 6 (q+1)\big(\|\hat h_{1}-\hat h_{2}\|_{\infty}+\omega_{\delta}(\hat h_{1})\big) + |\zeta_{1}-\zeta_{2}|,
\end{equation}
where $\hat h_{j}$ is the extension of $h_{j}$ to $\R_{+}$ by setting $\hat h_{j}(x)=0$ for all $x\ge \zeta_{j}$, and $\omega_{\delta}(\hat h_{1})=\sup\{|\hat h_{1}(s)-\hat h_{1}(t)|: |s-t|\le \delta\}$ is the $\delta$-modulus of continuity of $\hat h_{1}$.

\subsubsection{The continuum Erd\H{o}s--R\'enyi graph}
\label{sec: ERgraph}

We recall from Proposition~\ref{thm: cv-light} the stochastic process $\cS^{\lambda, \infty}$ and its height process $\cH^{\lambda, \infty}$: 
\[
\cS^{\lambda,\infty}_{t} = W_{t} + 2\lambda t - \tfrac12 t^{2}, \quad \cH^{\lambda,\infty}_{t}=2\Big(\cS^{\lambda, \infty}_{t}-\min_{s\le t}\cS^{\lambda,\infty}_{s}\Big), 
\]
where $(W_{t})_{t\ge 0}$ is a standard linear Brownian motion. As Proposition~\ref{thm: cv-light} indicates, the pair $(\cS^{\lambda, \infty}, \cH^{\lambda, \infty})$ is a continuum analogue for $(S\n_{k}, H\n_{k})_{k\ge 0}$, the discrete depth-first walk and its height process. We also recall the identity \eqref{def: ht-proc}, which expresses $(H\n_{k})_{k\ge 0}$ as a deterministic function of $(S\n_{k})_{k\ge 0}$. As it turns out, there is an analogue of \eqref{def: ht-proc} for spectrally positive L\'evy processes, discovered by Le Gall and Le Jan \cite{LGLJ98}. Roughly speaking, if $X=(X_{t})_{t\ge 0}$ is a spectrally positive L\'evy process, then for $t>0$, the following limit exists in probability:
\[
\cH^{X}_{t} = \lim_{\epsilon\to 0}\frac{1}{\epsilon} \int_{0}^{t} \mathbf 1_{\{X_{s}\le \inf_{s\le u\le t}X_{u}+\epsilon\}}ds,
\]
leading to the definition of a height process for $X$; see \cite{DuLG02, LGLJ98} for more details. However, in the special case where $X=(c\cdot W_{t})_{t\ge 0}$ for some $c>0$, it turns out that 
\begin{equation}
\label{id: ht-procs}
\cH^{X}_{t} = \frac{2}{c^{2}}\Big(X_{t}-\min_{s\le t}X_{s}\Big) \quad \text{a.s.}
\end{equation}
for each $t$.  This was discovered in \cite{MaMo03}; see also Eq.~(1.7) in \cite{DuLG02}. The identity then extends to $\cS^{\lambda,\infty}$ by Girsanov's Theorem, which explains the expression of $\cH^{\lambda, \infty}$ as seen above. 

We can now give an explicit construction of $\cG^{\mathrm{ER}}(\lambda)$ that appears in the scaling limit of critical $G(n, p_{e})$ (\cite{ABBrGo12}) as well as of critical $G(n,m,p)$ in the light and heavy clustering regimes. For $t>0$, let $g(t)$ (resp.~$d(t)$) be the last zero of $\cH^{\lambda, \infty}$ before $t$ (resp.~first zero of $\cH^{\lambda, \infty}$ after $t$). Alternatively, $(g(t), d(t))$ is the excursion interval of $\cS^{\lambda,\infty}$ above its running infimum that contains $t$. 
It is shown in \cite{Al97} that almost surely 
we can rank these excursions in a decreasing order of their lengths. For $k\ge 1$, let $(\mathrm g^{\lambda, \infty}_{k}, \mathrm d^{\lambda,\infty}_{k})$ be the $k$-th longest such excursion interval, which is unique almost surely. Recall $\zeta^{\lambda, \infty}_{k}$ from \eqref{def: zeta}; we have $\zeta^{\lambda,\infty}_{k}=\mathrm d^{\lambda,\infty}_{k}-\mathrm g^{\lambda,\infty}_{k}$. The excursion of $\cH^{\lambda, \infty}$ running on the interval $(\mathrm g^{\lambda,\infty}_{k}, \mathrm d^{\lambda,\infty}_{k})$ is denoted as $\mathrm e^{\lambda,\infty}_{k}$, namely,  
\[
e^{\lambda,\infty}_{k}: [0, \zeta^{\lambda,\infty}_{k}]\to \R_{+} \quad\text{with}\quad e^{\lambda,\infty}_{k}(s) = \cH^{\lambda, \infty}_{s+\mathrm g^{\lambda,\infty}_{k}}.
\]
Recall from Proposition~\ref{prop: surplus-bi} the Poisson point process $\cP^{\lambda, \infty}$ on $\R^{2}$, which has a finite number of atoms on every $(a, b)\times \R$ for $0\le a<b<\infty$. Let 
\[
p^{\lambda,\infty}_{k}= \#\{\cP^{\lambda, \infty}\cap (\mathrm g^{\lambda,\infty}_{k}, \mathrm d^{\lambda, \infty}_{k})\times \R\}\,, 
\]
and write $\{(x_{k, i}, y_{k, i}): 1\le i\le p^{\lambda,\infty}_{k}\}$ for the elements in $\cP^{\lambda, \infty}\cap (\mathrm g^{\lambda,\infty}_{k}, \mathrm d^{\lambda, \infty}_{k})\times \R$. By the definition of $\cP^{\lambda,\infty}$, we have a.s.
\[
0<y_{k, i}\le \cS^{\lambda, \infty}_{x_{k, i}}-\min_{u\le x_{k, i}}\cS^{\lambda,\infty}_{u}, \quad 1\le i\le p^{\lambda,\infty}_{k}. 
\]
For $1\le i\le p^{\lambda,\infty}_k$, we define
\[
s_{k, i}= \sup\Big\{u\le x_{k, i}: \cS^{\lambda, \infty}_{u}\le \,  \cS^{\lambda,\infty}_{x_{k, i}}-y_{k, i}\Big\}-\mathrm g^{\lambda,\infty}_{k}.
\]
It can be checked that $0\le s_{k, i}\le  t_{k, i}:=x_{k, i}-\mathrm g^{\lambda,\infty}_{k}$. We refer to Section~\ref{sec: cv-ghp} for some intuition behind this definition. 
Set $\Pi^{\lambda,\infty}_{k}=\{(s_{k, i}, t_{k, i}): 1\le i\le p^{\lambda,\infty}_{k}\}$. Recalling from \eqref{def: graph} the definition of $\cG(h, \Pi)$, we define the $k$-th largest component of $\cG^{\mathrm{ER}}(\lambda)$ as 
\[
(\cC^{\lambda, \infty}_{k}, d^{\lambda,\infty}_{k}, \mu^{\lambda,\infty}_{k}) = \cG\big(\mathrm e^{\lambda,\infty}_{k}, \Pi^{\lambda,\infty}_{k}\big).
\]
Let us also point out that conditional on  their respective sizes $\zeta^{\lambda,\infty}_{k}$ and the numbers of shortcuts $p^{\lambda,\infty}_{k}$, the measured metric spaces $(\cC^{\lambda, \infty}_{k}, d^{\lambda,\infty}_{k}, \mu^{\lambda,\infty}_{k})$, $k\ge 1$, are rescaled versions of each other. 
To explain this, let $\varepsilon=(\varepsilon_{t})_{0\le t\le 1}$ be the normalised Brownian excursion of length 1. For $\ell\ge 0$, let $(X_{i}, Y_{i})_{1\le i\le \ell}$ be $\ell$ pairs of random points from $D:=\{(x, y): 0\le x\le 1, 0\le y\le \varepsilon_{x}\}$ whose distribution is characterised as follows: for all suitable test functions $F, G$, we have
\[
\mathbb E\Big[F(\varepsilon)G\big((X_{i}, Y_{i})_{1\le i\le \ell}\big)\Big] = \frac{\mathbb E[F(\varepsilon) G\big((x_{i}, y_{i})_{1\le i\le \ell}\big)\int\prod_{1\le i\le \ell}\mathbf 1_{\{(x_{i}, y_{i})\in D\}}dx_{i}dy_{i}]}{\mathbb E[(\int_0^{1}\varepsilon_{s}ds)^{\ell}]}. 
\]
Let $Y'_{i}=\sup\{u\le X_{i}: \varepsilon_{u}\le \varepsilon_{X_{i}}-Y_{i}\}$, $1\le i\le \ell$, and set $\Pi^{(\ell)}=\{(Y'_{i}, X_{i}): 1\le i\le \ell\}$.  Finally, we define
\begin{equation}
\label{def: cond-gr}
\big(\mathrm C^{(\ell)}, \mathrm d^{(\ell)}, \mathrm m^{(\ell)}\big) = \cG(2\varepsilon, \Pi^{(\ell)})
\end{equation}
Then for each $\sigma\in (0, \infty)$, we have
\[
\text{cond.~on }\zeta^{\lambda,\infty}_{k}=\sigma, p^{\lambda,\infty}_{k}=\ell,\quad (\cC^{\lambda, \infty}_{k}, d^{\lambda,\infty}_{k}, \mu^{\lambda,\infty}_{k}) \eqd (\mathrm C^{(\ell)}, \sqrt{\sigma}\cdot\mathrm d^{(\ell)}, \sigma\cdot\mathrm m^{(\ell)}), \quad k\ge 1. 
\]

\subsubsection{The continuum graph $\mathcal G^{\mathrm{RIG}}(\lambda, \theta)$}
\label{sec: RIGgraph}

Recall from Proposition~\ref{thm: cv-moderate} the process $\cS^{\lambda, \theta}$. We set 
\[
\kappa_{\theta} = \big(1+\theta^{-\frac12}\big)^{\frac13} \quad\text{and}\quad \lambda_{\theta}=\lambda \cdot\kappa_{\theta}^{-2}, 
\]
and define $(B_{t})_{t\ge 0}$ by
\[
\kappa_{\theta}^{3/2} B_{t} = W_{t} + \theta^{-\frac14}W^{\ast}_{t}, \quad t\ge 0.
\]
Observe that $(B_{t})_{t\ge 0}$ is distributed as a standard linear Brownian motion. 
Applying the identity~\eqref{id: ht-procs} to $X=(\kappa_{\theta}^{3/2} B_{t})_{t\ge 0}$ yields 
\[
 \cH^{X}_{t} =  \frac{2}{\kappa_{\theta}^{3}}\Big(X_{t}-\min_{s\le t}X_{s}\Big) \quad \text{a.s.}
 \]
This again extends to $\cS^{\lambda, \theta}$ by Girsanov's Theorem, and thus giving us the expression of $\cH^{\lambda, \theta}$ in \eqref{def: ht}. 
On the other hand, using the scaling property of the Brownian motion, we have
\begin{align}\notag
\big(\cS^{\lambda, \theta}_{t}\big)_{t\ge 0} & = \Big(\kappa_{\theta}^{3/2} B_{t} - \tfrac12\kappa_{\theta}^{3}t^{2} + 2\lambda t\Big)_{t\ge 0}
\eqd \Big(\kappa_{\theta}W_{\kappa_{\theta}t}-\tfrac12\kappa_{\theta}^{3} t^{2}+2\lambda t\Big)_{t\ge 0}\\ \label{eq: scale}
 &=  \kappa_{\theta} \Big(W_{\kappa_{\theta}t}-\tfrac12(\kappa_{\theta}t)^{2} + 2\lambda_{\theta} \kappa_{\theta} t\Big)_{t\ge 0} 
 = \big(\kappa_{\theta}\, \cS^{\lambda_{\theta}, \infty}_{\kappa_{\theta}t}\big)_{t\ge 0}.
\end{align}
In particular, combined with the results from \cite{Al97}, this shows that as previously, we can rank the excursions of $\cH^{\lambda,\theta}$ above $0$ in a decreasing order of their lengths. 
For $k\ge 1$, let $(\mathrm g^{\lambda,\theta}_{k}, \mathrm d^{\lambda,\theta}_{k})$ be the $k$-th longest such excursion interval, so that $\zeta^{\lambda, \theta}_{k}=\mathrm d^{\lambda,\theta}_{k}-\mathrm g^{\lambda,\theta}_{k}$. Let $\mathrm e^{\lambda,\theta}_{k}$ denote the excursion of $\cH^{\lambda, \theta}$ running on $(\mathrm g^{\lambda,\theta}_{k}, \mathrm d^{\lambda,\theta}_{k})$. Define $\Pi^{\lambda,\theta}_{k}$ in the same way as in Section~\ref{sec: ERgraph}, replacing $\cP^{\lambda, \infty}$ with $\cP^{\lambda, \theta}$. Let $p^{\lambda,\theta}=\#\Pi^{\lambda,\theta}_{k}$ be the number of shortcuts. The $k$-th largest 
component of $\cG^{\mathrm{RIG}}(\lambda,\theta)$ is then defined as 
\[
\big(\cC^{\lambda, \theta}_{k}, d^{\lambda,\theta}_{k}, \mu^{\lambda,\theta}_{k}\big) = \cG\big(\mathrm e^{\lambda,\theta}_{k}, \Pi^{\lambda,\theta}_{k}\big).
\]
The identity~\eqref{eq: scale} allows us to compare $\cG^{\mathrm{RIG}}(\lambda,\theta)$ with $\cG^{\mathrm{ER}}(\lambda_{\theta})$. Indeed,  combined with \eqref{def: ht}, it implies that $(\cH^{\lambda,\theta}_{t})_{t\ge 0} \eqd (\kappa_{\theta}^{-2}\cH^{\lambda_{\theta},\infty}_{\kappa_{\theta} t})_{t\ge 0}$. Recall that $p^{\lambda,\theta}_k$ is a Poisson variable of mean 
\[
\int_{g^{\lambda, \theta}_k}^{d^{\lambda, \theta}_k} \Big(\cS^{\lambda, \theta}_s-\min_{u\le s}\cS^{\lambda,\theta}_u\Big)ds.
\]
We deduce from~\eqref{eq: scale} that
\[
\Big(\zeta^{\lambda,\theta}_{k}, p^{\lambda,\theta}_{k}\Big)_{k\ge 1}\eqd \Big(\kappa_{\theta}^{-1}\cdot\zeta^{\lambda_{\theta},\infty}_{k}, p^{\lambda_{\theta}, \infty}_{k}\Big)_{k\ge 1}.
\]
Moreover, for each $\sigma\in (0, \infty)$ and $\ell\in \N\cup\{0\}$, we have 
\[
\text{cond.~on } \zeta^{\lambda,\theta}_{k}=\sigma, p^{\lambda,\theta}_{k}=\ell, \quad \big(\cC^{\lambda, \theta}_{k}, d^{\lambda,\theta}_{k}, \mu^{\lambda,\theta}_{k}\big) \eqd \Big(\mathrm C^{(\ell)}, \tfrac{\sqrt{\sigma}}{\kappa_{\theta}^{3/2}}\cdot\mathrm d^{(\ell)}, \sigma\cdot\mathrm m^{(\ell)}\Big), \quad k\ge 1,
\]
where $(\mathrm C^{(\ell)}, \mathrm d^{(\ell)}, \mathrm m^{(\ell)})$ is defined in \eqref{def: cond-gr}. 

\subsection{Discussion}
\label{sec: discuss}

\paragraph{Asymptotic equivalence with the Erd\H{o}s--R\'enyi graph. }
Let $(G_{n})_{n\in\N}, (G^{\ast}_{n})_{n\in \N}$ be two sequences of random graphs. We say the two sequences are {\it asymptotically equivalent} if
\[
 \dTV\big(\sL(G_{n}), \sL(G_{n}^{\ast})\big) \to 0, \quad n\to\infty,
 \] 
where $\dTV$ stands for the total variation distance between two probability measures, and $\sL(X)$ denotes the law of a random element $X$.
To simplify the presentation, here we assume $m = \lf n^{\alpha}\rf$ with $\alpha\in (1, \infty)$ and $p=\frac{1}{\sqrt{mn}}(1+\lambda n^{-1/3})$ for some $\lambda\in \R$. 
Let $p_{e}$ denote the probability of having an edge between two vertices in $G(n, m, p)$, namely,
\[
p_{e} = 1-(1-p^{2})^{m}\sim \frac{1}{n}\big(1+2\lambda n^{-\frac13}\big), \quad n\to\infty. 
\]
In such case, results from \cite{FSS2000,BrBrNa2020} tell us the threshold of asymptotic equivalence for $(G(n, m, p))_{n\ge 1}$ and $(G(n, p_{e}))_{n\ge 1}$ is at $\alpha=3$. More precisely, the total variation distance between $\sL(G(n, m, p))$ and $\sL(G(n, p_{e}))$ tends to 0 if $\alpha>3$, and the same distance tends to 1 at least for $\alpha\in (2, 3)$. Focusing on the case $\alpha>3$ for the moment, we note as an immediate consequence of this asymptotic equivalence, large connected components of $G(n,m,p)$ will have the same scaling limit as found in critical $G(n, p_{e})$ (i.e.~special case of Theorem~\ref{thm1}) and with high probabilities these large components will contain no triangles (i.e.~special case of Theorem~\ref{thm: triangle-light}), as is the case with $G(n, p_{e})$. On the other hand, Theorems~\ref{thm1} and~\ref{thm: triangle-light} also reveal that if we look at individual functions of the graphs, the two models can exhibit the same asymptotic behaviours at a point much earlier than $\alpha=3$ (1 for scaling limit and $\frac73$ for triangle counts in the large connected components). A related question is as follows: what is the smallest $\alpha$ to ensure that the respective largest connected components in $G(n,m,p)$ and $G(n, p_{e})$ are asymptotically equivalent?  Given that their triangle counts match as soon as $\alpha>\frac73$, can the threshold be less than $3$?

\paragraph{Functional Gromov--Hausdorff--Prokhorov convergence.}
The distributional convergence in \eqref{cv: ER-graph} of $G(n, p_e)$, as shown by Addario-Berry, Broutin and Goldschmidt~\cite{ABBrGo12}, takes place in a stronger topology than the product topology appearing in Theorems~\ref{thm2}-\ref{thm3}. Indeed, 
for two sequences of compact measured metric spaces $\mathbf G=\{(G_{k}, d_{k}, \mu_{k}): k\ge 1\}$ and $\mathbf G'=\{(G'_{k}, d'_{k}, \mu'_{k}): k\ge 1\}$ satisfying respectively $\sum_{k\ge 1}\mathrm{diam}(G_{k})^{4}+\mu_k(G_k)^4<\infty$ and $\sum_{k\ge 1}\mathrm{diam}(G'_{k})^{4}+\mu'_k(G'_k)^4<\infty$ with $\mathrm{diam}(\cdot)$ standing for the diameter of the relevant metric space, define 
\[
\delta(\mathbf G, \mathbf G') = \Big(\sum_{k\ge 1}\dghp\big((G_{k}, d_{k}, \mu_{k}), (G'_{k}, d'_{k}, \mu'_{k})\big)^{4}\Big)^{\frac14}.
\]
The convergence in \eqref{cv: ER-graph} in fact holds in the weak topology induced by $\delta$. 
Let us briefly discuss what will be required to obtain a similar $\delta$-weak convergence for $G(n, m, p)$. Following the strategy of~\cite{ABBrGo12}, we can first show that the convergence of the component sizes of $G(n,m,p)$ takes place in the weak $\ell^2$-topology on the space of square summable sequences. In the light and moderate clustering regimes, this is a direct consequence of Aldous' theory on size-biased point processes; see~\eqref{eq: cv-l2} below. In the heavy clustering regime, this has been proved by Federico~\cite{Federico19+} via concentration inequalities. Back to $G(n, p_e)$, it is shown in \cite{ABBrGo12} that following the $\ell^{2}$-convergence of the component sizes, we have
\begin{equation}
\label{eq: diam-tail}
\lim_{N\to\infty}\limsup_{n\to\infty}\mathbb P\Big(\sum_{k\ge N} n^{-\frac43}\mathrm{diam}(G^n_k)^4\ge \epsilon\Big)\to 0.
\end{equation}
We note that the derivation of \eqref{eq: diam-tail} in~\cite{ABBrGo12} leans on the fact that 
each connected component in $G(n, p_e)$, when conditioned on its numbers of vertices and edges, is uniformly distributed. In contrast, when conditioned on their respective sizes, the components of $G(n,m,p)$ have much less tractable distributions. Therefore, it is not clear how to obtain an analogous estimate for random intersection graphs in general. 

\section{Proof of the main results}
\label{sec: proof}

This proof section is organised as follows. Our first goal is to prove Proposition~\ref{thm: cv-moderate}, which identifies the scaling limit of the graph exploration processes in the moderate clustering regime. This will take place in Sections~\ref{sec: pf}-\ref{sec: cv-den}, beginning with an overview of the proof in Section~\ref{sec: pf}. 
In Section~\ref{sec: light}, we prove Proposition~\ref{thm: cv-light}, the analogue of Proposition~\ref{thm: cv-moderate} for the light clustering regime. 
Proposition~\ref{prop: surplus-bi}, which concerns the limit of bipartite surplus edges, is proved in Section~\ref{sec: surplus}. Equipped with Propositions~\ref{thm: cv-moderate}, \ref{thm: cv-light}, and Proposition~\ref{prop: surplus-bi}, we are able to give the proof for our first two main results: Theorem~\ref{thm2} and~\ref{thm1}; this is done in Section~\ref{sec: cv-ghp}. 
We next proceed to the triangle counts and show Theorems~\ref{thm: triangle-moderate} and~\ref{thm: triangle-light} in Section~\ref{sec: triangle}. 
Up to that point, our arguments only concern the moderate and light clustering regimes. 
The heavy clustering regime is dealt with in Section~\ref{sec: heavy}. 

Unless otherwise specified, all random variables in the sequel are defined on some common probability space $(\Omega, \mathcal F, \mathbb P)$. We write Poisson$(c)$ to denote a Poisson distribution with expectation $c$, and Binom$(n, p)$ to denote a Binomial distribution with parameters $n\in \N$ and $p\in [0, 1]$. We use the notation $\Rightarrow$ to indicate convergence in distribution. 

\subsection{A preview of the main proof}
\label{sec: pf}

\subsubsection{Preliminaries on the depth-first exploration}

We discuss here some combinatorial features of Algorithm~\ref{alg}, which are later used in the proof of  Proposition~\ref{thm: cv-moderate}. 
We first observe that the sets $\cV_{k-1}$, $\cA_{k-1}$ and $\{v_{1}, v_{2}, \dots, v_{k}\}$ form a  partition of $\cV$.  This  yields the identity:
\begin{equation}
\label{id: VA}
\forall\, 1\le k\le n: \quad \#\cV_{k-1}+\#\cA_{k-1}+k= n.
\end{equation}
We also recall the spanning bipartite forest $\cF$ output by Algorithm~\ref{alg}. Let $\cT^{(i)}$ be the $i$-th connected component of $\cF$ ranked in the order of appearance when running the algorithm. Denote respectively by $v_{k_{-}(i)}$ and $v_{k_{+}(i)}$ the first and last vertices of $\cT^{(i)}$ explored in the algorithm. We note that the vertices of $\cT^{(i)}$ are precisely those $v_{k}$ with $k_{-}(i)\le k\le k_{+}(i)$. We also have the following observation. 
\begin{lem}
\label{lem: cA}
For all $n\ge 1$ and $1\le k\le n$,  the following statements hold true.
\begin{enumerate}[(i)]
\item
We have
\begin{equation}
\label{id: Ak}
S\n_{k}-\min_{j\le k}S\n_{j} = \#\cA_{k}= (\# \mathcal A^{\ast}_{k}-1)_{+},
\end{equation}
where $a_{+}$ denotes the positive part of a real number $a$. 
\item
$v_{k}$ is the last explored vertex in its connected component  if and only if $S\n_{k}=\min_{0\le j< k}S\n_{j}-1$.
\item
$v_{k}\in \cT^{(i)}$ if and only if $i=-\min_{j\le k-1}S\n_{j}+1$.
\end{enumerate}
\end{lem}

\begin{proof}
We first consider the vertices in $\cT^{(1)}$, that is, $1=k_{-}(1)\le k\le k_{+}(1)$. By definition, we have 
\[
S\n_{k}=\sum_{j\le k}(\#\cN_{j}-1) =\#\Big\{\bigcup_{j\le k}\cN_{j}\Big\}-k=\#\Big(\bigcup_{j\le k}\cN_{j}\setminus\{v_{2}, v_{3}, \dots, v_{k}\}\Big)-1=\#\cA^{\ast}_{k}-1,
\]
where we have used  the facts that $\cN_{j}, j\le k$, are disjoint in the second identity, that $v_{2}, v_{3}, \cdots, v_{k}$ all belong to $\cup_{j\le k}\cN_{j}$ in the third identity, and the definition of $\cA^{\ast}_{k}$ in the last one. We note that $k_{+}(1)$ is the first $k$ for which $\cA^{\ast}_{k}=\varnothing$. The above then implies that for all $k<k_{+}(1)$, $S\n_{k}\ge 0$ and that $S\n_{k_{+}(1)}=-1$. This proves (i)-(iii) for all $k$ up to $k_{+}(1)$. Proceeding to the vertices of $\cT^{(2)}$, we have
\[
S\n_{k}-S\n_{k_{+}(1)} = \sum_{k_{-}(2)\le j\le k}(\#\cN_{j}-1), \quad k_{-}(2)\le k\le k_{+}(2).
\]
By the same argument as before, this is again equal to $\#\cA^{\ast}_{k}-1$. Proceeding as in the previous case, we then extend the statements (i)-(iii) to all $k$ up to $k_{+}(2)$. Iterating this procedure completes the proof. 
\end{proof}

Recall from Algorithm~\ref{alg} the sets $\cU_k$ and $\cV_k$;  let us denote 
\begin{equation}
\label{def: uv}
V_{k}=\#\cV_{k-1}\quad \text{and} \quad U_{k}=\#\,\cU_{k-1}, \quad k\ge 1. 
\end{equation}
Recall from~\eqref{def: dps} that $R\n_{k-1}=\sum_{j\le k-1}\#\mathcal M_j=\sum_{j\le k-1}X\n_j$.
We note that $U_{k}$ and $V_{k}$ are respectively determined by $(R\n_{j})_{j<k}$ and $(S\n_{j})_{j<k}$, as we deduce from~\eqref{id: VA} and \eqref{id: Ak} that
\begin{equation}
\label{id: UkVk}
U_{k}=m-\sum_{j<k}X\n_{j}=m-R\n_{k-1}, \quad V_{k}=n-k-(S\n_{k-1}-\min_{j\le k-1}S\n_{j}). 
\end{equation}

Another by-product of Lemma~\ref{lem: cA} is the following proof of Lemma~\ref{lem: ht}. 
\begin{proof}[Proof of Lemma~\ref{lem: ht}]
This is adapted from the usual combinatorial arguments; see for instance \cite{LeG05}. 
Viewing the rooted tree $\cT^{(i)}$ as a family tree, we note that the height of a vertex in $\cT^{(i)}$ is the number of ancestors of $v$ ($v$ itself excluded). 
Since $\cT^{(i)}$ is a bipartite tree, the latter is equal to twice the number of ancestors of $v$ that belong to $\cV$. 
As the right-hand side of \eqref{def: ht-proc} always counts the term $k-1$, the conclusion will follow once we show that $v_{j+1}$ is an ancestor of $v_{k}$ if and only if $j< k-1$ and $\min_{j\le i\le k-1} S\n_{i}=S\n_{j}$. 
If $v_{j+1}$ is an ancestor of $v_{k}$, then necessarily $j<k-1$; moreover, the depth-first nature of Algorithm~\ref{alg} implies that each $v_{i}$, $j+1<i\le k$, is a descendent of $v_{j+1}$. As a result, we have $\cA_{j}\subseteq \cA_{i}$, $j<i< k$. 
Together with \eqref{id: Ak}, this entails $S\n_{j}-\min_{l\le j}S\n_{l}\le \min_{j\le i\le k-1}(S\n_{i}-\min_{l\le i}S\n_{l})$.
Meanwhile, we deduce from (iii) of Lemma~\ref{lem: cA} that $i\mapsto \min_{l\le i}S\n_l$ is constant on $[j, k-1]$, since the vertices $v_{i}$, $j+1\le i\le k$, all belong to the same tree. 
Combining this with the previous argument, we find that $S\n_{j}=\min_{j\le i\le k-1}S\n_{i}$. 
If $j<k-1$ and $v_{j+1}$ is not an ancestor of $v_{k}$, then there is some $v_{i}$ with $j< i< k$ which is the last descendant of $v_{j+1}$ (set $v_{i}=v_{j+1}$ if the latter has no descendent).  Then $\cA_{j}=\cA_{i}\cup\{v_{i+1}\}$, where $v_{i+1}$ is either a sibling of $v_{j+1}$ or a sibling of an ancestor of $v_{j+1}$. Using \eqref{id: Ak} yields $\min_{j\le l\le k-1}S\n_{l}\le S\n_{i}< S\n_{j}$. 
\end{proof}

\subsubsection{Outline of the proof of Proposition~\ref{thm: cv-moderate}} 

Let us first point out that the convergence of $(R\n_{k}, S\n_{k})_{k\ge 0}$, as a semi-martingale, can be handled using standard tools from stochastic analysis. %; see for instance the arguments in \cite{Al97}. %Therefore, the main focus of the proof is on the convergence of the height process $(H\n_{k})_{k\ge 0}$. 
However, as its definition \eqref{def: ht-proc} indicates, the height process is not a continuous functional of $(S\n_k)_{k\ge 0}$ in the uniform topology. Therefore, its convergence will not follow automatically from the convergence of $(R\n_{k}, S\n_{k})_{k\ge 1}$. Here, we adopt an approach inspired by \cite{CG23} that will yield simultaneously the convergences of $(R\n_{k}, S\n_{k})_{k\ge 0}$ and $(H\n_k)_{k\ge 0}$.   

To describe this approach, let us first consider the increment distribution of the depth-first walk $(R\n_{k}, S\n_{k})_{k\ge 0}$. Recall $U_{k}$ and $V_{k}$ from~\eqref{def: uv}. We note that conditional on $U_{k}$, the increment $\Delta R\n_{k}=X\n_{k}= \#\cM_{k}$ follows a Binom$(U_{k}, p)$ distribution. %Recall that the elements of $\cM_{k}$ are $u_{k, i}, 1\le i\le X\n_{k}$. 
For $1\le i\le X\n_{k}$, let $V_{k, i}=\#(\cV_{k-1}\setminus \cup_{j<i}\cN_{k, j})$. Then conditional on $V_{k, i}$, $\#\cN_{k, i}$ is distributed as a Binom$(V_{k, i}, p)$ variable, independent of $X\n_{k}$. 
If one ignores the difference between $V_{k}$ and $V_{k, i}$, then $\Delta S\n_{k}+1=\sum_{1\le i\le X\n_{k}}\#\cN_{k, i}$ can be approximated by a Binom$(X\n_{k}V_{k}, p)$ variable. Recall the identity \eqref{id: UkVk} that relates $(U_{k}, V_{k})$ to $(R\n_{j}, S\n_{j})_{j<k}$. 

\medskip
\noindent
{\bf Step 1: Poisson approximation. }
The first step in our approach consists in replacing $(R\n, S\n)$ by a similar walk but with Poisson increments. 
Let $\check U_{1}=m, \check V_{1}=n-1, \check R\n_{0}=0$, and $\check S\n_{0}=0$. For $1\le k\le n$, conditional on $\check U_{k}$, let $\check X\n_{k}:=\check R\n_{k}-\check R\n_{k-1}$ be a Poisson$(\check U_{k}p)$ variable, independent of $(\check R\n_{j},  \check S\n_{j})_{j<k}$; conditional on $\check V_{k}$ and $\check X\n_{k}$, let $\check S\n_{k}-\check S\n_{k-1}+1$ be a Poisson$(\check X\n_{k}\cdot \check V_{k}p)$ variable, independent of $(\check R\n_{j},  \check S\n_{j})_{j<k}$; set also 
\begin{equation}
\label{def: checkUk}
\check U_{k+1}=(m-\check R\n_{k})_{+}, \quad \check V_{k+1}=\Big(n-k-1-(\check S\n_{k}-\min_{j\le k}\check S\n_{j})\Big)_{+}. 
\end{equation}
Recall from \eqref{def: ht-proc} the functional $\mathscr H_{k}$. 
The height process $(\check H\n_k)_{k\ge 1}$ of $(\check S\n_k)_{k\ge 0}$ is  defined as 
\[
\check H\n_k = \mathscr H_{k}\Big((\check S\n_j)_{0\le j\le k}\Big), \quad k\ge 1.
\]
Recall the notation $\dTV$ for the total variation distance. Let $\sL(X)$ stand for the law of a random element $X$. 
Proof of the following assertion can be found in Section \ref{sec: tv}.
\begin{prop}
\label{prop: tv}
Assume either \eqref{hyp: L} or \eqref{hyp: M}. For any $t>0$, we have 
\[
\dTV\Big(\sL\big((X\n_{k}, S\n_{k}, H\n_{k})_{k\le tn^{2/3}}\big), \sL\big((\check X\n_{k}, \check S\n_{k}, \check H\n_{k})_{k\le tn^{2/3}}\big)\Big)\to 0, \quad n\to\infty.
\]
\end{prop}

Thanks to Proposition~\ref{prop: tv}, we only need to prove Proposition~\ref{thm: cv-moderate} for $(\check R\n_k, \check S\n_k, \check H\n_k)_{k\ge 1}$. 

\medskip
\noindent
{\bf Step 2:  Absolute continuity w.r.t.~walks with i.i.d.~increments. } 
Let $\lambda_{+}=\max(\lambda, 0)$. We take some  $\lambda_{\ast}\in (\lambda_{+}, \infty)$ and set
\[
\alpha_{n} = mp-\theta^{\frac12} \lambda_{\ast} n^{-\frac13} \quad\text{and}\quad \beta_{n} = np- \theta^{-\frac12}\lambda_{\ast} n^{-\frac13}. 
\]
We note that under \eqref{hyp: M}, both $\alpha_n$ and $\beta_n$ are positive for large $n$. Let $n_0=\min\{n: \alpha_n>0, \beta_n>0\}$. From now on, we only consider those $n\ge n_0$. 
Let $(\hat R\n_{k})_{k\ge 1}$ be a random walk with independent increments $\hat X\n_{k}:=\hat R\n_{k}-\hat R\n_{k-1}$ following a Poisson$(\alpha_{n})$ distribution. 
Given $(\hat X\n_{k})_{k\ge 1}$, let $(\hat S\n_{k})_{k\ge 0}$ be a random walk with i.i.d.~increments such that $\hat S\n_{0}=0$ and $\Delta \hat S\n_{k}+1$ is distributed according to Poisson$(\hat X\n_{k}\cdot \beta_{n})$. Define $\hat H\n_k=\mathscr H_{k}((\hat S\n_j)_{0\le j\le k})$ with the functional $\mathscr H_{k}$ from \eqref{def: ht-proc}. 
A theorem due to Duquesne and Le Gall (see Section \ref{sec: cvrw} for more details) gives us the following convergence of $(\hat R\n_{k}, \hat S\n_{k}, \hat H\n_{k})_{k\ge 1}$. 

\begin{prop}
\label{prop: cvrw}
Under  \eqref{hyp: M}, we have the following convergence in distribution in $\mathbb D^{\ast}(\R_+, \R^3)$: 
\begin{equation}
\label{cv: moderate-tree}
\Big\{n^{-\frac13}\hat R\n_{\lfloor n^{2/3}t\rfloor}-n^{\frac13}\theta^{\frac12} t, n^{-\frac13}\hat S\n_{\lf n^{2/3}t\rf}, n^{-\frac13}\hat H\n_{\lfloor n^{2/3}t\rfloor}: t\ge 0\Big\} \Longrightarrow \Big\{  \hat{\cR}_{t}, \hat{\cS}_{t}, \hat{\cH}_{t}: t\ge 0\Big\},
\end{equation}
where $\hat{\cR}_{t}, \hat{\cS}_{t}, \hat{\cH}_{t}$ are defined as follows: 
\[
 \hat{\cR}_{t} = \theta^{\frac14} W^{\ast}_{t} + \theta^{\frac12}(\lambda-\lambda_{\ast})t, \ 
\hat{\cS}_{t}  = W_{t}+\theta^{-\frac14}W^{\ast}_{t}+2(\lambda-\lambda_{\ast}) t, \  \hat{\cH}_{t} = \frac{2}{1+\theta^{-\frac12}}\Big(\hat{\cS}_{t} - \min_{s\le t}\hat{\cS}_{s}\Big). 
\]
\end{prop}

Our choice of $\lambda_{\ast}$ ensures $\mathbb E[\hat\cS_{t}]<0$ for $t>0$; this is required for the  application of the aforementioned theorem of Duquesne and Le Gall.   
On the other hand, Girsanov's Theorem implies the following absolute continuity relation between $(\cR^{\lambda,\theta}, \cS^{\lambda,\theta})$ and $(\hat\cR, \hat\cS)$: let $t>0$ and $F: \mathbb C([0, t], \R^{2})\to \R_+$ be a measurable function; then 
\begin{equation}
\label{id: abs-cont}
\mathbb E\Big[F\big((\cR^{\lambda,\theta}_{s}, \cS^{\lambda, \theta}_{s})_{0\le s\le t}\big)\Big] = \mathbb E\Big[F\big((\hat\cR_{s}, \hat\cS_{s})_{0\le s\le t}\big)\cdot \cE_{t}\big((\hat\cR_{s}, \hat\cS_{s})_{0\le s\le t}\big)\Big],
\end{equation}
where 
\begin{multline*}
\mathcal E_{t}\big((\hat\cR_{s}, \hat\cS_{s})_{0\le s\le t}\big) =\\
 \exp\Big(-\int_{0}^{t} s \, d\hat\cS_{s}+\lambda_{\ast}\hat\cS_{t}+\lambda_{\ast}(1-\theta^{-\frac12})\hat\cR_{t}+\frac12(1+\theta^{\frac12})(\lambda_{\ast}^{2}-2\lambda\lambda_{\ast})t+\lambda t^{2}-\frac16(1+\theta^{-\frac12})  t^{3}\Big). 
\end{multline*}
We refer to Section \ref{sec: Girsanov} for a proof of \eqref{id: abs-cont}. 
An analogue of \eqref{id: abs-cont} also exists for the discrete walks. 
Fix any $N\in \mathbb N$ and let $F: (\mathbb Z^2)^{N}\to \R_{+}$ be a measurable function; we have
\begin{equation}
\label{id: dens-dis}
\mathbb E\Big[F\big((\check R\n_k, \check S\n_{k})_{0\le k\le N}\big)\Big] = \mathbb E\Big[F\big((\hat R\n_k, \hat S\n_{k})_{0\le k\le N}\big)\cdot E\n_N\big((\hat R\n_{k}, \hat S\n_{k})_{0\le k\le N}\big)\Big],
\end{equation}
where the discrete Radon--Nikodym derivative is defined as
\begin{equation}
\label{def: En}
E\n_N\big((\hat R\n_{k}, \hat S\n_{k})_{0\le k\le N}\big) = \!\!\!\prod_{1\le k\le N}\Big(\frac{\hat U_{k}p}{\alpha_{n}}\Big)^{\Delta \hat R\n_k}\Big(\frac{\hat V_{k}p}{\beta_{n}}\Big)^{\Delta \hat S\n_{k}+1}e^{\alpha_{n}-\hat U_{k}p+\Delta \hat R\n_k(\beta_{n}-\hat V_{k}p)}.
\end{equation}
with $\hat U_{k}, \hat V_{k}$ defined by
\begin{equation}
\label{def: hatUk}
\hat U_{k}=(m-\hat R\n_{k-1})_{+}, \quad \hat V_{k}= \Big(n-k-(\hat S\n_{k-1}-\min_{j\le k-1}\hat S\n_{j})\Big)_{+}.
\end{equation}
Proof of \eqref{id: dens-dis} can be found in Section \ref{sec: cv-den}. 

\medskip
\noindent
{\bf Step 3: Convergence of the Radon--Nikodym derivatives. } The final ingredient of our proof is provided by the following proposition. 

\begin{prop}
\label{prop: cvdens}
Assume \eqref{hyp: M}. For all $t\ge 0$,  the sequence 
\[
\Big\{E\n_{\lf tn^{2/3}\rf}\Big((\hat R\n_{k}, \hat S\n_{k})_{0\le k\le tn^{2/3}}\Big): n\ge 1\Big\}
\]
 is uniformly integrable.  Moreover, jointly with the convergence in  \eqref{cv: moderate-tree}, for any $t\ge 0$, 
\begin{equation}
\label{eq: cvdens}
E\n_{\lf tn^{2/3}\rf}\Big((\hat R\n_{k}, \hat S\n_{k})_{0\le k\le tn^{2/3}}\Big) \ \Longrightarrow\  \mathcal E_{t}\big((\hat\cR_{s}, \hat\cS_{s})_{0\le s\le t}\big).
\end{equation}
\end{prop}

\begin{proof}[Proof of Proposition~\ref{thm: cv-moderate} subject to Propositions~\ref{prop: tv},~\ref{prop: cvrw} and~\ref{prop: cvdens}]
Since the height process is a measurable function of the depth-first walk, both in the discrete and continuum cases, we deduce from \eqref{id: dens-dis}, the uniform integrability of $(E\n_{\lf tn^{2/3}\rf})_{n\ge 1}$, the convergence in \eqref{eq: cvdens}, as well as \eqref{id: abs-cont}, that for any continuous and bounded function $F: \mathbb D^{\ast}([0, t], \R^{3})\to \R_+$, we have
\begin{align*}
&\quad  \mathbb E\Big[F\Big(\Big(n^{-\frac13}\check R\n_{\lf n^{2/3}s\rf}, n^{-\frac13}\check S\n_{\lf n^{2/3}s\rf}, n^{-\frac13}\check H\n_{\lf n^{2/3}s\rf}\Big)_{s\le t}\Big)\Big] \\
&=  \mathbb E\Big[F\Big(\Big(n^{-\frac13}\hat R\n_{\lf n^{2/3}s\rf}, n^{-\frac13}\hat S\n_{\lf n^{2/3}s\rf}, n^{-\frac13}\hat H\n_{\lf n^{2/3}s\rf}\Big)_{s\le t}\Big)\cdot E\n_{\lf tn^{2/3}\rf}\Big((\hat R\n_{k}, \hat S\n_{k})_{k\le tn^{2/3}}\Big)\Big]\\
& \to \mathbb E\Big[F\big((\hat \cR_s, \hat \cS_s, \hat \cH_s)_{0\le s\le t}\big)\cdot \mathcal E_{t}\big((\hat\cR_{s}, \hat\cS_{s})_{s\le t}\big)\Big] = \mathbb E\Big[F\big((\cR^{\lambda,\theta}_s, \cS^{\lambda,\theta}_s, \cH^{\lambda,\theta}_s)_{s\le t}\big)\Big].
\end{align*}
Together with Proposition~\ref{prop: tv}, this yields the convergence in \eqref{cv: moderate}. 
\end{proof}

\medskip
\noi
{\bf Some elementary facts used in the subsequent proofs. }
For all $a_{1}, a_{2}\in \R$, we have $(a_{1}+a_{2})^{2}\le 2a_{1}^{2}+2a_{2}^{2}$; more generally for $k\ge 2$ and $a_{i}\in \R, 1\le i\le k$, we have
\begin{equation}
\label{bd: sum}
\Big(\sum_{i=1}^{k}a_{i}\Big)^{2} \le 2^{k-1}\sum_{i=1}^{k}a_{i}^{2}.
\end{equation}
Let $p\in [0, 1]$ and $1\le n_{1}\le n_{2}$ be two integers. Let $X_{1}$ (resp.~$X_{2}$) be a Binom$(n_{1}, p)$ variable (resp.~Binom$(n_{2}, p)$ variable). Then $X_{1}$ is stochastically bounded by $X_{2}$, namely, $\mathbb P(X_{1}\ge a)\le \mathbb P(X_{2}\ge a)$ for all $a\in\R$. 
Similarly, for $0<c_{1}<c_{2}$ and two Poisson variables $Y_{1}, Y_{2}$ of respective expectations $c_{1}$ and $c_{2}$, $Y_{1}$ is stochastically bounded by $Y_{2}$. 
We also require some estimates on the Binomial and Poisson distributions, which are collected in Appendix~\ref{sec: binom}. 

\medskip
\noi
{\bf Big O and small o notation. }
For a function $g$ with the variables $n, m$ that potentially depends on some additional parameters and a positive function $f$ that does not depend on these parameters, we write 
\[
g(n, m)= \cO(f(n, m))
\]
to indicate that there is some constant $C\in (0, \infty)$ that may well depend on all these parameters so that $|g(n, m)|\le C\cdot f(n, m)$ for all $n, m\ge 1$. Similarly, we write $g(n, m)=o(f(n, m))$ if $g(n, m)/f(n,m)\to 0$ as $n\to\infty$ for each fixed set of parameters. 

\subsection{Poisson approximation: Proof of Proposition~\ref{prop: tv}} 
\label{sec: tv}

Let us recall that for  two probability measures $\mu$ and $\nu$ supported on $\R$, $\dTV(\mu, \nu)<\epsilon$ if and only if we can find a random vector $(X, Y)$ with $X\sim\mu$ and $Y\sim\nu$ satisfying $\mathbb P(X\ne Y)<\epsilon$. We will refer to the vector $(X, Y)$ as an {\it (optimal) coupling} of $\mu$ and $\nu$. 
Here we prove Proposition~\ref{prop: tv} by finding an appropriate coupling between the distributions of random walks.  We start with a well-known total variation bound between Binomial and Poisson distributions (see for instance Theorem 2 in \cite{Cam60}): %Let $N\in \N$ and $p\in [0, 1]$. Denote by $\mu_{p}$ the law Binom$(N, p)$ and by $\nu_{p}$ the law Poisson$(Np)$. We have  
\begin{equation}
\label{bd: binom-tv}
\dTV\big(\text{Binom}(N, p), \text{Poisson}(Np)\big)\le Np^{2} \quad \text{for all } N\ge 1 \text{ and } p\in [0, 1]. 
\end{equation}
Since the sum of two independent Poisson variables also follows a Poisson distribution, we also have 
\begin{equation}
\label{bd: poisson-tv}
\dTV\big(\text{Poisson}(c_{1}), \text{Poisson}(c_{2})\big)\le \mathbb P(\text{Poisson}(c_{2}-c_{1})>0)\le  c_{2}-c_{1}, \quad \text{for all } 0<c_{1}<c_{2}.
\end{equation} 
Since $X\n_{1}$ (resp.~$\check X\n_{1}$) follows a Binom$(m, p)$ (resp.~Poisson$(mp)$) distribution, the previous bound 
implies that we can find two random variables $\rX^{n}_{1}, \check\rX^{n}_{1}$ satisfying $\rX^{n}_{1}\eqd X\n_{1}$, $\check\rX^{n}_{1}\eqd \check X\n_{1}$ and $\mathbb P(\rX\n_{1}\ne \check\rX\n_{1}) \le mp^{2}$.
From now on, instead of introducing a set of new random variables that realises the optimal coupling, we will simply say  that there is a coupling so that $\mathbb P(X\n_{1}\ne \check X\n_{1})\le mp^{2}$. As we are interested in the distributions of $X\n_{k}$ and $\check X\n_{k}$ rather than the random variables themselves, we believe this abuse of notation would not lead to confusion. 

Recall from \eqref{def: checkUk} the definition of $\check U_{k}$ and that $U_{2}=m-X\n_{1}$. Using \eqref{bd: binom-tv} again allows us to find a coupling so that
\[
\mathbb P(X\n_{2}\ne \check X\n_{2}\,|\, X\n_{1}=\check X\n_{1}) \le \mathbb E\big[\mathbb P(X^{n}_{2}\ne \check X\n_{2}\,|\, U_{2}, \check U_{2}=U_{2})\big]\le mp^{2}.
\]
Iterating this procedure, we can find a coupling between $(X\n_{k}, U_{k})_{1\le k\le n}$ and $(\check X\n_{k}, \check U_{k})_{1\le k\le n}$ for each $n\in \N$ so that
\begin{equation}
\label{bd: tv1}
\mathbb P\Big(\exists\, k\le tn^{2/3}: X\n_{k}\ne \check X\n_{k}\Big) \le\!\! \sum_{k\le tn^{2/3}}\!\!\mathbb P(X\n_{k}\ne \check X\n_{k}\,|\,X\n_{j}=\check X\n_{j}, j<k) \le t n^{\frac23}mp^{2}, 
\end{equation}
which tends to 0 under either the assumption~\eqref{hyp: L} or~\eqref{hyp: M}. We now seek to expand this coupling to include $(S\n_{k})_{k\ge 1}$ and $(\check S\n_{k})_{k\ge 1}$. Recall that for $1\le i\le X\n_{k}$, $\#\cN_{k, i}$ is a Binom$(V_{k, i}, p)$ variable, where $V_{k, i}=V_{k}-\sum_{j<i}\#\cN_{k, j}$. In particular, $\#\cN_{1, 1}$ is a Binom$(n-1, p)$ variable. 
By \eqref{bd: binom-tv}, we can find a Poisson$((n-1)p)$-variable ${\mathrm N}_{1, 1}$ satisfying 
\[
\mathbb P(\#\cN_{1, 1}\ne {\mathrm N}_{1, 1}) \le np^{2}.
\]
Let ${\mathrm V}_{1, 2}=n-1-{\mathrm N}_{1, 1}$ and repeat this procedure for $i\ge 2$. We end up with a collection $({\mathrm N}_{1, i})_{i\ge 1}$ of respective Poisson$({\mathrm V}_{1, i}p)$ distributions with ${\mathrm V}_{1, i}=(n-1-\sum_{j<i}{\mathrm N}_{1, j})_{+}$ that satisfy
\begin{equation}
\label{tvbd1}
\mathbb P\Big(\exists\, 1\le i\le X\n_1: \#\cN_{1, i} \ne {\mathrm N}_{1, i}\,\Big|\,X\n_{1}=\check X\n_{1}\Big)\le \mathbb E[X\n_{1}] \cdot np^{2} \le mnp^{3}.  
\end{equation}
On the other hand, recalling~\eqref{bd: poisson-tv}, we can find a coupling between $\mathrm N_{1, i}$ and a Poisson$((n-1)p)$ variable $\check{\mathrm N}_{1, i}$ so that 
\[
\mathbb P(\mathrm N_{1, i}\ne \check{\mathrm N}_{1, i})\le \mathbb E\big[(n-1-\mathrm V_{1, i})p\big]\le \sum_{j<i}p\,\mathbb E\big[\mathrm N_{1, j}]\le inp^2, 
\]
where we have used the fact that each $\mathrm N_{1, j}$ is stochastically bounded by Poisson$(np)$. 
Summing over $i$ yields 
\begin{equation}
\label{bd: tvre}
\mathbb P\Big(\exists\, 1\le i\le X\n_{1}: \mathrm N_{1, i}\ne \check{\mathrm N}_{1, i}\big)\le np^2\,\mathbb E\Big[\sum_{1\le i\le X\n_{1}} i\,\Big]
 \le np^{2}\mathbb E\big[(X\n_{1})^{2}\big]\le m^{2}np^{4}+mnp^{3}. 
\end{equation}
%where we have noted that $\mathrm N_{1, i}$ is stochastically bounded by Poisson$(np)$. 
Set $\check S\n_{1}+1= \sum_{1\le i\le \check X\n_{1}}\check{\mathrm N}_{1, i}$, which then follows a Poisson$(\check X\n_{1}(n-1)p)$ distribution. Together with~\eqref{tvbd1}, the previous bound implies
\begin{equation}
\label{tvbd2}
\mathbb P\big(\check S\n_{1}\ne S\n_{1}\,\big|\, X\n_{1}=\check X\n_{1}\big) \le m^{2}np^{4}+2mnp^{3}. 
\end{equation}
Note from~\eqref{id: UkVk} and~\eqref{def: checkUk} that $V_{k}$ and $\check V_{k}$ are the same functional of $(S\n_{j})_{j<k}$ and $(\check S\n_{j})_{j<k}$ respectively.  Repeating the previous procedure finds us a coupling between $(S\n_{k})_{1\le k\le n}$ and $(\check S\n_{k})_{1\le k\le n}$ so that
\[
\mathbb P\Big(\exists\, k\le tn^{2/3}: S\n_{k} \ne \check S\n_{k}\,\Big|\, X\n_{k}= \check X\n_{k}, \forall\, k\le tn^{2/3}\Big) \le tn^{\frac23}(2mnp^{3}+m^{2}np^{4}) \to 0,
\]
under both~\eqref{hyp: M} and~\eqref{hyp: L}. 
Together with \eqref{bd: tv1}, this proves Proposition~\ref{prop: tv}, as $H\n_{k}$ and $\check H\n_{k}$ are the same functional  of $(S\n_{j})_{1\le j\le k}$ and $(\check S\n_{j})_{1\le j\le k}$, respectively. 

\subsection{Convergence of the i.i.d.~random walks and the associated height process}
\label{sec: cvrw}

We provide here a proof of Proposition~\ref{prop: cvrw}, relying upon prior results from \cite{DuLG02, BrDuWa21}. To that end, we note that $(\hat R\n_{k}, \hat S\n_{k})_{k\ge 0}$ is a bi-dimensional random walk with i.i.d.~increments. Under the assumption \eqref{hyp: M}, we have
\begin{align}
 \label{eq: moments}
&\mathbb E[\hat R\n_{1}]-\theta^{\frac12} = \theta^{\frac12}(\lambda-\lambda_{\ast}) n^{-\frac13}+o(n^{-\frac13}), \quad
\mathbb E[\hat S\n_{1}] = 2(\lambda-\lambda_{\ast}) n^{-\frac13}+o(n^{-\frac13}),\\ \notag
& \Var(\hat R\n_{1})=\theta^{\frac12} + o(1), \quad
\Var(\hat S\n_{1}) = 1+\theta^{-\frac12}+o(1), \quad
\mathrm{Cov}(\hat R\n_{1}, \hat S\n_{1}) = 1+o(1).
\end{align}
The calculations for the estimates in~\eqref{eq: moments} are mostly straightforward. We only give details for the last two. For the variance of $\hat S\n_1$, we note that the Law of total variance implies
\begin{align*}
\Var(\hat S\n_{1}) &= \mathbb E\big[\Var(\hat S\n_1\,|\,\hat X\n_1)\big]+\Var\big(\mathbb E[\hat S\n_1\,|\, \hat X\n_1]\big)\\
&=\beta_n\mathbb E[\hat X\n_1]+\beta_n^2\Var(\hat X\n_1) =\alpha_n\beta_n(1+\beta_n).
\end{align*}
Under~\eqref{hyp: M}, we have $\alpha_n\beta_n=1+o(1)$ and $\beta_n=\theta^{-1/2}+o(1)$; thus $\Var(\hat S\n_{1}) = 1+\theta^{-\frac12}+o(1)$. For the covariance, by first conditioning on $\hat X\n_1$, we find that
\begin{align*}
\mathrm{Cov}(\hat R\n_{1}, \hat S\n_{1}) & = \mathbb E[\hat R\n_1\cdot \hat S\n_1]-\mathbb E[\hat R\n_1]\cdot \mathbb E[\hat S\n_1]\\
&=\beta_n\mathbb E\big[(\hat X\n_1)^2\big]-\beta_n\big(\mathbb E[\hat X\n_1]\big)^2= \beta_n\Var(\hat X\n_1) = \alpha_n\beta_n.
\end{align*}
The conclusion then follows. From~\eqref{eq: moments}, we apply the martingale central limit theorem (recalled in Appendix~\ref{sec: mart}) to derive the following convergence. 

\begin{prop}
Under~\eqref{hyp: M}, we have the following convergence in distribution in $\mathbb D^{\ast}(\R_+, \R^{2})$:
\begin{equation}
\label{eq: cv-firsthalf}
\Big\{n^{-\frac13}\hat R\n_{\lfloor n^{2/3}t\rfloor}-n^{\frac13}\theta^{\frac12}t, \; n^{-\frac13}\hat S\n_{\lf n^{2/3}t\rf}: t\ge 0\Big\} \Longrightarrow \Big\{  \hat{\cR}_{t}, \hat{\cS}_{t}: t\ge 0\Big\}.
\end{equation}
\end{prop}

\begin{proof}
For each $n\ge 1$, let us denote 
\[
E^{n,1}_t = \mathbb E\Big(\hat R\n_{\lfloor t\rfloor}\Big) \quad \text{and} \quad E^{n,2}_t = \mathbb E\Big(\hat S\n_{\lfloor t\rfloor}\Big), \quad t\ge 0.
\]
Let $\mathcal G^n_t$ be the $\sigma$-algebra generated by $(\hat R\n_k, \hat S\n_k), k\le t$. Since $(\hat R\n, \hat S\n)$ has independent increments, it is not difficult to see that 
\[
M^{n, 1}_t=\hat R\n_{\lfloor t\rfloor}-E^{n, 1}_t, \quad M^{n, 2}_t = \hat S\n_{\lfloor t\rfloor}-E^{n, 2}_t, \quad t\ge 0,
\]
is a pair of local martingales with respect to $(\mathcal G^n_t)_{t\ge 0}$. Let us show that 
\begin{equation}
\label{cv: mart-bid}
\Big\{n^{-\frac13}M^{n, 1}_{n^{2/3}t} \;, n^{-\frac13}M^{n, 2}_{n^{2/3}t}: t\ge 0\Big\} \Longrightarrow \Big\{  \theta^{\frac14}W^{\ast}_t, W_t+\theta^{-\frac14}W^{\ast}_t: t\ge 0\Big\}.
\end{equation}
To that end, we use a version of martingale central limit theorem as recalled in Section~\ref{sec: mart} with  $A^n(t)=(A^n_{i,j}(t))_{1\le i, j\le 2}$ chosen as follows:
\begin{align*}
A^n_{1,1}(t) &= n^{-\frac23}\Var\big(M^{n, 1}_{n^{2/3}t}\big)=n^{-\frac23}\Var\big(\hat R\n_{\lf n^{2/3}t\rf }\big), \ A^n_{2,2}(t) = n^{-\frac23}\Var\big(M^{n, 2}_{n^{2/3}t}\big) =n^{-\frac23}\Var\big(S\n_{\lf n^{2/3}t\rf}\big)\\ 
A^n_{1,2}(t) &=A^n_{2, 1}(t)= n^{-\frac23}\mathrm{Cov}\big(M^{n, 1}_{n^{2/3}t}, M^{n,2}_{n^{2/3}t}\big)=n^{-\frac23}\mathrm{Cov}\big(\hat R\n_{\lf n^{2/3}t\rf}, \hat S\n_{\lf n^{2/3}t\rf}\big).
\end{align*}
Thanks again to the independence of the increments, we readily check that $M^{n, i}_{n^{2/3}t}M^{n, j}_{n^{2/3}t}-A^n_{i,j}(t)$ is a local martingale with respect to $(\mathcal G^n_t)_{t\ge 0}$, for $1\le i, j\le 2$. Conditions (i) and (ii) in Section~\ref{sec: mart} also hold. Moreover, the last three estimates in~\eqref{eq: moments} imply that~\eqref{cvcon: var} is true with $C_{1,1}=\theta^{1/2}, C_{2, 2}=1+\theta^{-1/2}$ and $C_{1, 2}=C_{2, 1}=1$. For~\eqref{cvcon: An}, we note that $A^n$ has fixed jump sizes with 
\[
\sup_{s\ge 0} \|A^n(s)-A^n(s-)\|\le \|A^n(n^{-\frac23})\|.
\]
\eqref{cvcon: An} then follows from~\eqref{eq: moments}. To check~\eqref{cvcon: jump}, let us we first establish the so-called Lindeberg's condition: for all $\epsilon>0$, 
\begin{equation}
\label{eq: bd-rw}
\mathbb E\Big[(\hat R\n_{1})^{2}\mathbf 1_{\{|\hat R\n_{1}|\ge \epsilon n^{1/3}\}}\Big] + \mathbb E\Big[(\hat S\n_{1})^{2}\mathbf 1_{\{|\hat S\n_{1}|\ge \epsilon n^{1/3}\}}\Big] \to 0.
\end{equation}
To see why~\eqref{eq: bd-rw} holds true, recall that $\hat S\n_1+1$ is a mixed Poisson variable with mean $\hat X\n_1\cdot \beta_n$ and $\hat X\n_1\sim$ Poisson$(\alpha_n)$. Applying successively~\eqref{id: poisson-moment} and~\eqref{id: poisson-moment'}, we find that 
\begin{align*}
\mathbb E\big[(\hat S\n_1+1)^4\big]&= \sum_{j=1}^4 {4 \brace j}\beta_n^j \cdot \mathbb E\Big[\big(\hat X\n_1\big)^j\Big] = \sum_{j=1}^4 {4 \brace j}\beta_n^j\sum_{i=1}^j {j\brace i} \mathbb E\Big[(\hat X\n_1)_i\Big]\\
&=\sum_{j=1}^4 {4 \brace j}\beta_n^j\sum_{i=1}^j {j\brace i} \alpha_n^i,
\end{align*}
where we have used~\eqref{id: poisson-moment} again in the second equality and \eqref{id: poisson-moment'} for the last line. This shows that $\mathbb E[(\hat S\n_1+1)^4]$ is bounded by a polynomial in $\alpha_{n}$ and $\beta_{n}$ of degree at most $4$. Since both $\alpha_{n},\beta_{n}$ remain bounded under \eqref{hyp: M}, we have
\[
\sup_{n\ge 1}\mathbb E\big[(\hat S\n_1+1)^4\big] <\infty. 
\]
With a similar calculation for $\hat R\n_{1}$, we deduce that 
\[
\sup_{n\ge 1}\mathbb E\Big[(\hat R\n_{1})^{4}\Big]+\sup_{n\ge 1}\mathbb E\Big[(\hat S\n_{1})^{4}\Big]<\infty.
\]
An application of the Markov inequality then yields~\eqref{eq: bd-rw}. From~\eqref{eq: bd-rw}, we deduce that 
\begin{align*}
\mathbb E\Big[\max_{k\le tn^{2/3}} (\Delta \hat S\n_k)^2\Big]
&\le \epsilon n^{\frac23}+ \sum_{k\le tn^{2/3}} \mathbb E\Big[(\Delta\hat S\n_k)^2\mathbf 1_{\{|\hat S\n_k|\ge \epsilon n^{1/3}\}}\Big]\\
&\le \epsilon n^{\frac23}+tn^{\frac23} \mathbb E\Big[(\hat S\n_1)^2\mathbf 1_{\{|\hat S\n_1|\ge \epsilon n^{1/3}\}}\Big]
\end{align*}
By first letting $n\to\infty$ and then $\epsilon\to 0$, we conclude that $\mathbb E\Big[\max_{k\le tn^{2/3}} (\Delta \hat S\n_k)^2\Big]=o(n^{2/3})$. A similar property holds for $\hat R\n$. Since we have
\[
n^{-\frac23}\sup_{s\le t} \|M^n(s)-M^n(s-)\|^2 \le n^{-\frac23}\max_{k\le tn^{2/3}} \Big(( \Delta \hat R\n_k)^2+(\Delta\hat S\n_k)^2+\big(\mathbb E[\Delta \hat R\n_k]\big)^2 + \big(\mathbb E[\Delta \hat S\n_k]\big)^2\Big), 
\]
the condition~\eqref{cvcon: jump} holds as a consequence of the previous arguments and~\eqref{eq: moments}. The martingale central limit theorem now applies to assert that~\eqref{cv: mart-bid} holds.  To conclude, we note that the first two estimates in~\eqref{eq: moments} ensure that 
\[
n^{-\frac13} E^{n, 1}_{tn^{2/3}}-n^{\frac13}\theta^{\frac12}t \to \theta^{\frac12}(\lambda-\lambda_{\ast})t \quad \text{and}\quad
n^{-\frac13} E^{n, 2}_{tn^{2/3}}\to 2(\lambda-\lambda_{\ast})t,
\]
uniformly on the compact sets. 
Recall from Proposition~\ref{prop: cvrw} the definitions of $(\hat\cR, \hat\cS)$. The conclusion follows. 
\end{proof}

Let $\mu_{n}$ denote the law of  $\hat S\n_{1}+1$, and let $g_{n}$ be the probability generating function of $\mu_{n}$, that is, 
\[
g_{n}(u)=\sum_{k\ge 0} \mu_{n}(k)u^{k}, \quad 0\le u\le 1.
\]
We also write $g^{\circ k}_{n}=\overbrace{g_{n}\circ g_{n}\circ\cdots \circ g_{n}}^{k}$ for the $k$-th iterated composition of $g_{n}$. 
Thanks to our choice of $\lambda_{\ast}$, we can find some $N\in \N$ so that $\mathbb E[\hat S\n_{1}]< 0$ for all $n\ge N$. As a result, the sequence $(\mu_{n})_{n\ge 1}$ is eventually subcritical.  
Assuming that there exists some $\delta>0$ satisfying
\begin{equation}
\label{eq: ext-tight}
\limsup_{n\to\infty}n^{1/3}\big(1-g^{\circ\lf \delta n^{1/3}\rf}_{n}(0)\big) <\infty, 
\end{equation}
For $n\ge N$,  let $(Z^n(k))_{k\ge 0}$ be a Bienaym\'e branching process with offspring distribution $\mu_n$ and $Z^n(0)=\lf n^{1/3}\rf$, and denote by $f(n, \delta)=\mathbb P(Z^n(\lf\delta n^{1/3}\rf)=0)$. We note that~\eqref{eq: ext-tight} is equivalent to saying that $\liminf_{n\to\infty}f(n, \delta)>0$. 
Using the monotonicity in $\delta$ and the branching property it can be checked  that if~\eqref{eq: ext-tight} for some $\delta>0$, then it holds for every $\delta>0$.
Theorem 2.3.1 in \cite{DuLG02} (see also Corollary 2.5.1 and Eq.~(1.7) there) then tells us that 
\begin{equation}
\label{eq: cv-secondhalf}
\Big\{n^{-\frac13}\hat S\n_{\lfloor n^{2/3}t\rfloor}, n^{-\frac13}\hat H\n_{\lf n^{2/3}t\rf}: t\ge 0\Big\} \Longrightarrow \Big\{  \hat{\cS}_{t}, \hat{\cH}_{t}: t\ge 0\Big\} \quad \text{in $\mathbb D^{\ast}(\R_+, \R^{2})$.} 
\end{equation}
Let us put aside the verification of \eqref{eq: ext-tight} for a moment and proceed with the rest of the proof. The convergences in \eqref{eq: cv-firsthalf} and \eqref{eq: cv-secondhalf} imply that the sequence of distributions of 
\[
\Big\{n^{-\frac13}\hat R\n_{\lfloor n^{2/3}t\rfloor}-n^{\frac13}\theta^{\frac12}t, n^{-\frac13}\hat S\n_{\lf n^{2/3}t\rf}, n^{-\frac13}\hat H\n_{\lfloor n^{2/3}t\rfloor}: t\ge 0\Big\} 
\]
is tight. To conclude, we only need to show the sequence has a unique limit point. Suppose that $(\cR', \cS', \cH')$ and $(\cR^{\ast}, \cS^{\ast}, \cH^{\ast})$ both appear as weak limits along  respective subsequences.  We must  have  $(\cR', \cS')\eqd (\cR^{\ast}, \cS^{\ast})$ as a result of \eqref{eq: cv-firsthalf}. Moreover, the expression of $\hat{\cH}_{t}$ in Proposition~\ref{prop: cvrw} tells us that $\hat\cH$ is distributed as a deterministic function of $\hat\cS$. 
It follows that $\cH'$ (resp.~$\cH^{\ast}$) is distributed as the same function of $\cS'$ (resp.~$\cS^{\ast}$). The joint law of $(\cR', \cS', \cH')$ is therefore the same as $(\cR^{\ast}, \cS^{\ast}, \cH^{\ast})$. 

It remains to check \eqref{eq: ext-tight}. A sufficient condition for \eqref{eq: ext-tight}, in the context of a specific family of mixed Poisson distributions, first appeared in \cite{BrDuWa21}. However, a closer look at the proof in \cite{BrDuWa21} reveals that it actually works for a general offspring distribution; we therefore re-state it as follows: 

\begin{prop}[Proposition 7.3 in~\cite{BrDuWa21}]
\label{prop: tight-suff}
Suppose that $\mu_{n}$ is a probability measure on $\Z_{+}=\N\cup\{0\}$ satisfying $\sum_{k\ge 1}k\mu_n(k)\le 1$. Suppose that $(a_{n})_{n\ge 1}, (b_{n})_{n\ge 1}$ are two sequences of positive real numbers that increase to $\infty$. Write $g_{\mu_{n}}$ for the probability generation function of $\mu_{n}$ and define
\[
\Psi_{n}(u) = a_{n}b_{n}\Big(g_{\mu_{n}}(1-\tfrac{u}{a_{n}})-1+\tfrac{u}{a_{n}}\Big), \quad 0\le u\le a_{n}.
\]
Assuming further that  
\begin{equation}
\label{eq: ext-tight'}
\lim_{y\to\infty}\limsup_{n\to\infty}\int_{y}^{a_{n}} \frac{du}{\Psi_{n}(u)} = 0,
\end{equation}
then we can find some $\delta>0$ so that 
\[
\limsup_{n\to\infty}a_{n}\big(1- g_{\mu_{n}}^{\circ\lf \delta b_{n}\rf}(0)\big)<\infty. 
\]
\end{prop}

\begin{proof}
The proof works verbatim as the proof of Proposition 7.3 in~\cite{BrDuWa21}. In particular, Lemma 7.4,  the key ingredient of the proof, which compares the lifetime of a discrete-time Bienaym\'e process to that of a continuous-time one, is valid for all (sub)critical offspring distributions. Note that what we have denoted as $b_n$ corresponds to $b_n/a_n$ there. 
\end{proof}

\begin{proof}[Proof of \eqref{eq: ext-tight}]
Our intention is to apply Proposition~\ref{prop: tight-suff} with $a_{n}=b_{n}=n^{1/3}$. To that end, 
let $\hat X$ be a Poisson variable of mean $\alpha_{n}=mp-\theta^{\frac12}\lambda_{\ast}n^{-\frac13}$. Recall that $\beta_{n}=np-\theta^{-\frac12}\lambda_{\ast}n^{-\frac13}$. We have
\[
g_{n}(u)=\mathbb E\Big[e^{-\beta_{n}(1-u)\hat X}\Big] = \exp\big\{-\alpha_{n}\big(1-e^{-\beta_{n}(1-u)}\big)\big\} \ge \exp\big\{-\alpha_{n}\beta_{n}(1-u)\big\},
\]
where we have used the inequality $1-e^{-x}\le x$ for all $x\ge 0$. Next, using the inequality $e^{-x}-1+x\ge \tfrac12 x^{2}e^{-x}$ valid for all $x\ge 0$, we obtain that 
\begin{align*}
\Psi_{n}(u)& = n^{\frac23}\Big(g_{n}(1-n^{-\frac13}u)-1+n^{-\frac13}u\Big) \ge n^{\frac23}\big(\exp(-\alpha_{n}\beta_{n}n^{-\frac13}u)-1+n^{-\frac13}u\big) \\
& \ge \tfrac12 \big(\alpha_{n}\beta_{n}u)^{2} e^{-\alpha_{n}\beta_{n}n^{-\frac13}u} +(1-\alpha_{n}\beta_{n})n^{\frac13}u.
\end{align*}
Since $\alpha_{n}\beta_{n}=1+2(\lambda-\lambda_{\ast})n^{-1/3}+o(n^{-1/3})$, we can find some $c\in (0, \infty)$ so that 
\[
\Psi_{n}(u)\ge c u^{2}, \quad \text{ for all } u\in [0, n^{\frac13}] \text{ and $n$ sufficiently large}.
\]
\eqref{eq: ext-tight'} now follows as a result. 
\end{proof}

To summarise, we have shown that \eqref{eq: ext-tight} is verified with the help of Proposition~\ref{prop: tight-suff}. Together with the previous arguments, this completes the proof of Proposition~\ref{prop: cvrw}.

\subsection{Convergence of the Radon--Nikodym derivatives: Proof of Proposition~\ref{prop: cvdens}} 
\label{sec: cv-den}

We first show the absolute continuity identity \eqref{id: dens-dis} for the depth-first walks with Poisson increments. 

\begin{proof}[Proof of~\eqref{id: dens-dis}]
We prove here a slightly more general version of~\eqref{id: dens-dis} for later applications. Recall from Section~\ref{sec: tv} that for each $k\ge 1$, $\check{\mathrm N}_{k, i}, i\ge 1$ are i.i.d.~Poisson$(\check V_k p)$, so that $\Delta\check S\n_k+1=\sum_{1\le i\le \check X\n_k}\check{\mathrm N}_{k, i}$ follows a Poisson$(\check X\n_k\check V_kp)$ distribution. 
Let $\hat N_{k, i}, i\ge 1, k\ge 1$ be i.i.d.~Poisson$(\beta_n)$ variables. We assume that $\Delta\hat S\n_k +1= \sum_{1\le i\le \hat X\n_k} \hat N_{k, i}$. Recall that $\check X\n_k$ (resp.~$\hat X\n_k$) is a Poisson$(\check U_k p)$ variable (resp.~Poisson$(\alpha_n)$ variable). 
Let $k\in \N$ and $x_{k}, n_{k, i}\in \Z_{+}=\N\cup\{0\}$; we have
\begin{align*}
&\; \mathbb P(\check X\n_{k}=x_{k}, \check{\mathrm N}_{k, i}=n_{k, i}, 1\le i\le x_k  \,|\, \check U_{k},\check V_{k}) = e^{-\check U_{k}p}\frac{(\check U_{k}p)^{x_{k}}}{x_{k}!}\cdot \prod_{1\le i\le x_k}e^{-\check V_{k}p}\frac{(\check V_{k}p)^{n_{k, i}}}{n_{k, i}!} \\
& = e^{-\alpha_{n}}\frac{(\alpha_{n})^{x_{k}}}{x_{k}!}\cdot \prod_{1\le i\le x_k}e^{-\beta_{n}}\frac{\beta_{n}^{n_{k, i}}}{n_{k, i}!} \cdot \Big(\frac{\check U_{k}p}{\alpha_{n}}\Big)^{x_{k}} \Big(\frac{\check V_{k}p}{\beta_{n}}\Big)^{\sum_i n_{k, i}}e^{\alpha_{n}-\check U_{k}p+x_{k}(\beta_{n}-\check V_{k}p)}\\
& = \mathbb P(\hat X\n_{k}=x_{k}, \hat N_{k, i}= n_{k, i}, 1\le i\le x_k)\cdot\Big(\frac{\check U_{k}p}{\alpha_{n}}\Big)^{x_{k}}\Big(\frac{\check V_{k}p}{\beta_{n}}\Big)^{\sum_i n_{k, i}}e^{\alpha_{n}-\check U_{k}p+x_{k}(\beta_{n}-\check V_{k}p)}.
\end{align*}
We recall that conditional on $(\check U_{N}, \check V_{N})$, $(\check X\n_{N}, (\check{\mathrm N}_{N, i})_{i\ge 1})$ are independent of $(\check X\n_{k}, (\check{\mathrm N}_{k, i})_{i\ge 1})_{k<N}$. 
Since $\check U_{N}, \check V_{N}$ are determined by $(\check R\n_{k}, \check S\n_{k})_{k<N}$ via \eqref{def: checkUk}, on the event that $\check X\n_{k}=x_{k}, \check{\mathrm N}_{k, i}=n_{k, i}, 1\le i\le x_k,  k<N$, let us write $\check u_{N}:=(m-\sum_{k<N}x_{k})_{+}$, $y_k:=\sum_{1\le i\le x_k}n_{k, i}$ and $\check v_{N}:=(n-N-(\sum_{k< N}(y_{k}-1)-\min_{k<N}\sum_{j\le k}(y_{j}-1))_{+}$. 
Then we have
\begin{align*}
&\; \mathbb P\big(\check X\n_{N}=x_{N}, \check{\mathrm N}_{N, i}=n_{N, i}, 1\le i\le x_N\,\big|\, \check X\n_{k}=x_{k}, \check{\mathrm N}_{k, i}=n_{k, i},1\le i\le x_k, 1\le k\le N-1 \big)\\ 
= &\; \mathbb P\big(\check X\n_{N}=x_{N}, \check{\mathrm N}_{N, i}=n_{N, i}, 1\le i\le x_N\,|\, \check U_{N}=\check u_{N}, \check V_{N}=\check v_{N}) \\
 = &\; \mathbb P(\hat X\n_{N}=x_{N}, \hat N_{N, i}= n_{N, i}, 1\le i\le x_N)\cdot\Big(\frac{\check u_{N}p}{\alpha_{n}}\Big)^{x_{N}}\Big(\frac{\check v_{N}p}{\beta_{n}}\Big)^{y_{N}}e^{\alpha_{n}-\check u_{N}p+x_{N}(\beta_{n}-\check v_{N}p)}.
\end{align*}
By iterating this procedure and comparing $\check u_{N}, \check v_{N}$ with \eqref{def: hatUk}, we find that for any suitable test function $F$: 
\begin{equation}
\label{id: dens-comp}
\mathbb E\Big[F\Big(\big(\check X\n_k, (\check{\mathrm N}_{k, i})_{i\le \check X\n_k}\big)_{k\le N}\Big)\Big] = \mathbb E\Big[F\Big(\big(\hat X\n_k, (\hat{N}_{k, i})_{ i\le \check X\n_k}\big)_{k\le N}\Big)E\n_N\big((\hat R\n_k, \hat S\n_k)_{k\le N})\Big], 
\end{equation}
where $E\n_N$ is defined in~\eqref{def: En}. 
Since $\Delta\check S_k$ is a function of $(\check{\mathrm N}_{k, i})$ and $\Delta\hat S_k$ is the same function of $(\hat N_{k, i})$, the conclusion follows. 
\end{proof}

We now proceed to the proof of Proposition \ref{prop: cvdens}.  Some of the straightforward computations used in the proof are summarised in the following lemma. %whose proof is left to the reader.  

\begin{lem}
\label{lem: interm}
As $n\to\infty$, 
\begin{equation}
\label{id: sum}
\sum_{k=1}^{\lf tn^{2/3}\rf}\Big(\frac{k}{n}\Big)^{2} = \frac13 t^{3}+o(1) \quad \text{ and } \ \sum_{k=1}^{\lf tn^{2/3}\rf}kn^{-\frac43} =\frac12 t^{2}+o(1).   
\end{equation}
Moreover, 
\[
\sum_{k=1}^{ \lf tn^{2/3}\rf}k\cdot \Delta \hat S\n_{k} =\sum_{k=1}^{ \lf tn^{2/3}\rf} k\cdot (\hat S\n_{k}-\hat S\n_{k-1}) =  \lf tn^{2/3}\rf \hat S\n_{\lf tn^{2/3}\rf} -\sum_{k=1}^{\lf tn^{2/3}\rf-1}\hat S\n_{k}.
\]
It follows that jointly with the convergence in~\eqref{cv: moderate-tree}, we have for any $(c_1, c_2)\in \R^2$:
\begin{equation}
\label{id: ipp}
\Big(c_1 n^{-\frac13}\hat R\n_{\lf tn^{2/3}\rf}-c_{1}n^{\frac13}\theta^{\frac12}t, c_2 n^{-\frac13}\hat S\n_{\lf tn^{2/3}\rf}, \sum_{k=1}^{ \lf tn^{2/3}\rf}\frac{k}{n}\cdot \Delta \hat S\n_{k} \Big)  \ \Longrightarrow\  \Big(c_1 \hat\cR_t, c_2 \hat\cS_t,  \int_{0}^{t} s\, d\hat\cS_{s}\Big). 
\end{equation}
\end{lem}

\begin{proof}
We briefly explain~\eqref{id: ipp}; the rest is rather straightforward. Let $\mathbf x=(x_s)_{0\le s\le t}$ be a continuous function. Since 
\[
\int_0^t x_s ds = \lim_{N\to\infty} \frac{1}{N}\sum_{i=1}^N x_{\frac{i}{N}t},
\]
it is not difficult to check that the mapping $\mathbf x\mapsto \int_0^t x_s ds$ is continuous in the uniform topology on $[0, t]$. We then derive from the convergence in~ \eqref{cv: moderate-tree} that
\begin{align*}
\sum_{k=1}^{ \lf tn^{2/3}\rf}\frac{k}{n}\cdot \Delta \hat S\n_{k}&=\frac{\lf tn^{2/3}\rf}{n} \hat S\n_{\lf tn^{2/3}\rf}-
\frac{1}{n}\sum_{k=1}^{\lf tn^{2/3}\rf-1}\hat S\n_{k} = \frac{\lf tn^{2/3}\rf}{n} \hat S\n_{\lf tn^{2/3}\rf}-\frac1n\int_0^{\lf tn^{2/3}\rf}\hat S\n_{\lf s\rf}ds \\
& = \frac{\lf tn^{2/3}\rf}{n} \hat S\n_{\lf tn^{2/3}\rf}-n^{-\frac13}\int_{0}^{\lfloor tn^{2/3}\rfloor/n^{2/3}} \hat S\n_{\lfloor sn^{2/3}\rfloor }ds 
\Longrightarrow\ t\hat\cS_t-\int_0^t \hat\cS_s ds=\int_0^t s\,d\hat\cS_s,
\end{align*}
by an integration by parts. The joint convergence in~\eqref{id: ipp} then follows, as it can be readily checked that all converging subsequences converge to a unique limit. 
\end{proof}

Recall from~\eqref{def: hatUk} the quantities $\hat U_{k}, \hat V_{k}$. 
The next lemma will provide us with the necessary estimates required in the proof of Proposition \ref{prop: cvdens}. 

\begin{lem}
\label{lem: hatUk}
Assume \eqref{hyp: M}. For all $t>0$, we have
\begin{equation}
\label{bd: hatUk}
\max_{k\le tn^{2/3}}\mathbb E\big[(m-\hat U_{k})^{2}\big]=\cO(n^{\frac43}), \quad \max_{k\le tn^{2/3}}\mathbb E\big[(n-\hat V_{k})^{2}\big]=\cO(n^{\frac43}).
\end{equation}
Let $\epsilon_{k}=m-\hat U_{k}-kmp$ and $\delta_{k}=n-\hat V_{k}-k$, $k\ge 1$. For all $t>0$, we also have 
\begin{equation}
\label{bd: rebduk}
\max_{k\le tn^{2/3}}\mathbb E\big[\epsilon_{k}^{2}\big] = \cO(n^{\frac23}),\quad \max_{k\le tn^{2/3}}\mathbb E\big[\delta_{k}^{2}\big]=\cO(n^{\frac23}).
\end{equation}
\end{lem}

\begin{proof}
To ease the writing, let us denote $N=\lf tn^{2/3}\rf$. 
By definition, we have
\[
0\le \, \delta_{k}=n-k-\big(n-k-(\hat S\n_{k-1}-\min_{j\le k-1}\hat S\n_{j})\big)_{+}\le \hat S\n_{k-1}-\min_{j\le k-1}\hat S\n_{j}.
\]
It follows that 
\[
\delta_k^2\le \Big(\hat S\n_{k-1}-\min_{j\le k-1}\hat S\n_j\Big)^2\le 2\big(\hat S\n_{k-1}\big)^2+2\Big(\max_{j\le k-1}-\hat S\n_j\Big)^2\le 4 \max_{k\le N} \, \big(\hat S\n_k\big)^2.
\]
Thus, the bound for $\mathbb E[\delta_{k}^{2}]$ will follow once we show $\mathbb E[\max_{k\le N}(\hat S\n_{k})^{2}]=\cO(n^{2/3})$. However, $\hat S\n_{k}-\mathbb E[\hat S\n_{k}]$, $k\ge 1$, is a martingale, which implies, via Doob's maximal inequality, that
\[
\mathbb E\Big[\max_{k\le N}\big(\hat S\n_{k}-\mathbb E[\hat S\n_{k}]\big)^{2}\Big]\le 4 N\cdot \Var(\hat S\n_{1}) = \cO(n^{\frac23}),
\]
thanks to~\eqref{eq: moments}. Combined with $\max_{k\le N}|\mathbb E[\hat S\n_{k}]|=\cO(n^{1/3})$, again due to~\eqref{eq: moments}, this gives the desired bound, which entails the uniform bounds for $\mathbb E[\delta_{k}^{2}]$, which in turn yields the uniform bound for $\mathbb E[(n-\hat V_{k})^{2}]=\cO(n^{\frac43})$ in~\eqref{bd: hatUk}. 
Regarding $\hat U_{k}$,  we have
\[
\epsilon_{k}^{2}=\Big(m-(m-\hat R\n_{k-1})_{+}-kmp\Big)^{2}\le 2 \Big(m-(m-\hat R\n_{k-1})_{+}-(k-1)\alpha_{n}\Big)^{2}+2\Big((k-1)\alpha_{n}-kmp\Big)^{2}.
\]
Since $0\le mp-\alpha_{n}=\cO(n^{-1/3})$ and $mp=\cO(1)$, the second term above can be bounded as follows:
\[
\max_{k\le N}\Big((k-1)\alpha_{n}-kmp\Big)^{2} = \cO(n^{\frac23}).
\]
For the other term, writing $\mathcal E_{k}=\{\hat R\n_{k-1}\le m\}$, we have for $n$ sufficiently large so that $\alpha_n>0$, 
\[
\mathbb E\Big[ \Big(m-(m-\hat R\n_{k-1})_{+}-(k-1)\alpha_{n}\Big)^{2}\Big]\le m^{2}\mathbb P(\mathcal E_{k}^{c})+\mathbb E\Big[\Big(\hat R\n_{k-1}-(k-1)\alpha_{n}\Big)^{2}\Big]. 
 \]
Since the increments of $\hat R\n_{k}$ are i.i.d.~Poisson$(\alpha_{n})$ variables, we deduce that
\[
\mathbb E\Big[\Big(\hat R\n_{k-1}-(k-1)\alpha_{n}\Big)^{2}\Big]= \sum_{j\le k-1}\Var(\hat X\n_{j})\le   kmp=\cO(m^{\frac12}n^{\frac16}). 
\]
Meanwhile, standard computations yield that $\max_{k\le N}\mathbb E[\exp(\hat R\n_{k})]\le \exp(c n^{2/3}mp)$ for some $c\in (0, \infty)$. Markov's inequality implies that  
\[
\max_{k\le N}\, m^{2}\mathbb P(\mathcal E_{k}^{c}) \le m^{2}e^{-m}\max_{k\le N}\mathbb E[\exp(\hat R\n_{k})] \le m^{2}\exp(-m+cn^{2/3}mp)=o(1). 
\]
Combined with the previous arguments, this shows $\mathbb E[\epsilon_{k}^{2}]=\cO(n^{2/3})$. Along with  $(Nmp)^{2}=\cO(mn^{\frac13})$, the bound for $\mathbb E[(m-\hat U_{N})^{2}]$ follows. 
\end{proof}

Let us briefly discuss the proof strategy for Proposition~\ref{prop: cvdens}. 
From now on, we fix some $q\in (0, 1)$. Let us define 
\[
\sigma_n=\inf\big\{k\in [n]: m-\hat U_k\ge (1-q)m \text{ or } n-\hat V_k \ge (1-q)n\big\},
\]
with the convention that $\inf\varnothing = \infty$. 
We note that $E\n_{\lf tn^{2/3}\rf}\big((\hat R\n_{k}, \hat S\n_{k})_{1\le k\le \lf tn^{2/3}\rf}\big)$ is a non negative random variable that has unit expectation, as implied by~\eqref{id: dens-dis}. Similarly, we deduce from~\eqref{id: abs-cont} that $\mathbb E\big[\mathcal E_{t}\big(\hat \cR_{s}, \hat\cS_{s})_{s\le t}\big)\big]=1$. By Lemma 4.8 in~\cite{CG23}, it suffices to show that
\begin{equation}
\label{cv: dens''}
 E\n_{\lf tn^{2/3}\rf}\big((\hat R\n_{k}, \hat S\n_{k})_{1\le k\le \lf tn^{2/3}\rf}\big)\mathbf 1_{\{\sigma_n> tn^{2/3}\}} \Longrightarrow \mathcal E_{t}\big(\hat \cR_{s}, \hat\cS_{s})_{s\le t}\big),
 \end{equation}
which will then entail the convergence in~\eqref{eq: cvdens} as well as the uniform integrability. Let us denote $N=\lfloor tn^{2/3}\rfloor$ to ease the writing. We observe that 
\[
\{\sigma_n> N\} =\cap_{1\le k\le N} A_k, \quad \text{where } A_k =\big\{m-\hat U_k<(1-q)m \text{ and } n-\hat V_k < (1-q)n\big\}.
\]
With the understanding that $\exp(-\infty)=0$, we can split the discrete Radon--Nikodym density as follows: 
\[
E\n_{N}\big((\hat R\n_{k}, \hat S\n_{k})_{1\le k\le N}\big) \mathbf 1_{\{\sigma_n>N\}} = \exp\bigg\{\sum_{1\le k\le N}\Big\{(F_{k}+H_{k}+I_{k}+J_{k}+L_{k})\mathbf 1_{A_k}+Q_k\Big\}\bigg\},
\]
where $Q_k$ is a random variable that takes the value $0$ on $A_k$ and $-\infty$ on $A_k^c$. Observe that if there is at least one $k$ for which $A_k^c$ holds, then the right-hand side above is 0. The other terms above are given by:
\begin{align*}
H_{k} & = (\Delta \hat R\n_{k}-\theta^{\frac12})\cdot \Big(\log\frac{\hat U_{k}}{m}+(n-\hat V_{k})p\Big), \\
I_{k}  &=\theta^{\frac12}\log\frac{\hat U_{k}}{m} + (m-\hat U_{k})p+ \log\frac{\hat V_{k}}{n}+ \theta^{\frac12} p(n-\hat V_{k}), \\
J_{k}& = \alpha_{n}-mp+\theta^{\frac12}\Big(\log\frac{mp}{\alpha_{n}} +\beta_{n}-np\Big)+\log\frac{np}{\beta_{n}}\,,\\
L_{k} &= \Delta \hat S\n_{k}  \log\frac{\hat V_{k}}{n}, \quad 
F_{k} = (\Delta\hat R\n_{k}-\theta^{\frac12})\Big(\log\frac{mp}{\alpha_{n}}+\beta_{n}-np\Big)+\Delta\hat S\n_{k}\log\frac{np}{\beta_{n}}\,. 
\end{align*}
We claim that 
\begin{equation}
\label{cv: sumQ}
\sum_{k\le \lfloor tn^{2/3}\rfloor}Q_k\xrightarrow[]{n\to\infty} 0 \text{ in probability}.
\end{equation}
Indeed, we have 
\begin{align}\notag
\mathbb P(\sigma_n\le N) &\le  \mathbb P\Big(\max_{1\le k\le N}(m-\hat U_k) \ge (1-q)m\Big) + \mathbb P\Big(\max_{k\le N}(n-\hat V_k)\le (1-q) n\Big)\\ \notag
& \le \max_{k\le N}\mathbb P\big(m-\hat U_k\ge (1-q)m\big) + \max_{k\le N}\mathbb P\big(n-\hat V_k\ge (1-q)n\big) \\ \label{bd: sigman}
&\le \max_{k\le N} \frac{\mathbb E[(m-\hat U_k)^2]}{(1-q)^2 m^2} + \max_{k\le N} \frac{\mathbb E[(n-\hat V_k)^2]}{(1-q)^2 n^2} = \cO(n^{-\frac23}),
\end{align}
where we have used the monotonicity of $\hat U_k, \hat V_k$ in the second line, Markov's inequality and~\eqref{bd: hatUk} in the third line. Then~\eqref{cv: sumQ} follows, since $\sum_{k\le N}Q_k\ne 0$ if and only if $\sigma_n\le N$. 
Now the proof of~\eqref{cv: dens''} boils down to proving the following string of lemmas.

\begin{lem}
\label{lem: pfcv1}
Under~\eqref{hyp: M}, we have 
\[
\sum_{k\le\lfloor tn^{2/3}\rfloor} I_k\mathbf 1_{A_k} \xrightarrow[]{n\to\infty} \lambda t^{2}-\frac16(1+\theta^{-\frac12})t^{3} \quad \text{in probability}.
\]
\end{lem}

\begin{lem}
\label{lem: pfcv2}
Under~\eqref{hyp: M}, we have 
\[
\sum_{k\le\lfloor tn^{2/3}\rfloor} H_k\mathbf 1_{A_k} \xrightarrow[]{n\to\infty} 0\quad \text{in probability}.
\]
\end{lem}

\begin{lem}
\label{lem: pfcv3}
Under~\eqref{hyp: M}, we have 
\[
\sum_{k\le\lfloor tn^{2/3}\rfloor} J_k\mathbf 1_{A_k} \xrightarrow[]{n\to\infty} \frac12(1+\theta^{\frac12})(\lambda_{\ast}^{2}-2\lambda_{\ast}\lambda)t\quad \text{in probability}.
\]
\end{lem}

\begin{lem}
\label{lem: pfcv4}
Under~\eqref{hyp: M}, jointly with the convergence in \eqref{cv: moderate-tree}, the following holds:
\[
\Big(\sum_{k\le tn^{2/3}} L_{k}\mathbf 1_{A_k}, \sum_{k\le tn^{2/3}} F_{k}\mathbf 1_{A_k}\Big)\  \Longrightarrow\  \Big(-\int_{0}^{t} s\,d\hat\cS_{s}, \;(1-\theta^{-\frac12})\lambda_{\ast}\hat\cR_{t}+\lambda_{\ast}\hat\cS_{t}\Big).
\]
\end{lem}

\begin{proof}[Proof of Lemma~\ref{lem: pfcv1}]
We write $\log(1-x)=-x-x^{2}/2+\eta(x)$, noting that for the fixed number $q\in (0, 1)$, there is a constant $C=C(q)\in (0,\infty)$ satisfying $|\eta(x)|\le Cx^{3}$ for all $x\in [0, 1-q]$. Recall from Lemma \ref{lem: hatUk} that $m-\hat U_{k}= kmp + \epsilon_{k}$. We find that
\begin{align*}
\theta^{\frac12}\log\frac{\hat U_{k}}{m} + (m-\hat U_{k})p
& = \theta^{\frac12}\cdot \Big\{ - \frac{m-\hat U_{k}}{m}-\frac12 \Big( \frac{m-\hat U_{k}}{m}\Big)^{2} + \eta\Big(\frac{m-\hat U_{k}}{m}\Big)\Big\}+(m-\hat U_{k})p\\
& = \Big(p-\frac{\sqrt\theta}{m}\Big)\cdot (kmp+\epsilon_{k}) -\frac{\sqrt\theta}{2}\Big( \frac{kmp +\epsilon_{k}}{m}\Big)^{2}+\sqrt\theta\eta\Big(\frac{m-\hat U_k}{m}\Big)
\end{align*}
Since $p-\sqrt\theta/m=\lambda\theta^{-1/2}n^{-4/3}(1+o(1))$ under \eqref{hyp: M}, together with \eqref{bd: rebduk}, this implies that
\[
\mathbb E\Big[\sum_{k\le N}\Big|p-\frac{\sqrt\theta}{m}\Big|\cdot |\epsilon_{k}|\mathbf 1_{A_k}\Big]\le \Big|p-\frac{\sqrt\theta}{m}\Big|\cdot \sum_{k\le N}\Big(\mathbb E[\epsilon_{k}^{2}]\Big)^{\frac12} = \cO( n^{-\frac13})\to 0.
\]
Similarly, we have 
\[
\mathbb E\Big[\sum_{k\le N}\Big|\Big(\frac{kmp+\epsilon_{k}}{m}\Big)^{2}-(kp)^{2} \Big|\mathbf 1_{A_k}\Big] \le  2m^{-1}Np\sum_{k\le N}\Big(\mathbb E[\epsilon_{k}^{2}]\Big)^{\frac12} +m^{-2}\sum_{k\le N}\mathbb E[\epsilon_{k}^{2}]= \cO( n^{-\frac13}). 
\]
On the other hand, the previous bound of $\eta$ together with the definition of $A_k$ implies that 
\[
\sum_{k\le N}\mathbb E\Big[\Big|\eta\Big(\frac{m-\hat U_k}{m}\Big)\Big|\mathbf 1_{A_k}\Big]\le C\sum_{k\le N}\mathbb E\Big[\Big(\frac{m-\hat U_k}{m}\Big)^3\Big]
\le C\sum_{k\le N}\mathbb E \Big[\Big(\frac{kmp+\epsilon_{k}}{m}\Big)^{3} \Big]= \cO(n^{-\frac13}), 
\]
where we have used the trivial bound $\epsilon_{k}\le m$, so that $\mathbb E[\epsilon_{k}^{3}]=\cO(n^{5/3})$ following \eqref{bd: rebduk}. 
With \eqref{id: sum}, we find that 
\[
\sum_{k\le tn^{2/3}}\Big(p-\frac{\sqrt\theta}{m}\Big) kmp -\frac{\sqrt\theta}{2}(kp)^{2}\ \longrightarrow\  \frac12 \lambda t^{2}-\frac16\theta^{-\frac12}t^{3}.
\]
We note that 
\begin{equation}
\label{bd: sumAc}
\mathbb E\Big[\sum_{k\le N}\mathbf 1_{A_k^c}\Big] \le  \mathbb E\Big[\sum_{k\le N}\mathbf 1_{\{\sigma_n\le k\}}\Big] \le  
\mathbb E[(N+1-\sigma_n)_+] \le (N+1)\cdot \mathbb P(\sigma_n\le N) = \cO(1),
\end{equation}
 by~\eqref{bd: sigman}. If follows that 
\[
\sum_{k\le N}\Big\{\Big(p-\frac{\sqrt\theta}{m}\Big) kmp -\frac{\sqrt\theta}{2}(kp)^{2}\Big\}\mathbf 1_{A_k^c} \le \Big(\Big|p-\frac{\sqrt\theta}{m}\Big|Nmp+\frac{\sqrt\theta }{2}N^2p^2\Big)\sum_{k\le N}\mathbf 1_{A_k^c} = \cO(n^{-\frac23}). 
\]
Combining all the previous bounds, we conclude that 
\begin{equation}
\label{eq: cvIk}
\sum_{k\le tn^{2/3}} \Big\{\theta^{\frac12}\log\frac{\hat U_{k}}{m} + (m-\hat U_{k})p\Big\}\mathbf 1_{A_k}\ \longrightarrow \ \frac12 \lambda t^{2}-\frac16\theta^{-\frac12}t^{3} \ \text{ in probability}. 
\end{equation}
In the same manner, we write
\begin{align*}
\log\frac{\hat V_{k}}{n}+ \theta^{\frac12} p(n-\hat V_{k}) & = -\frac{n-\hat V_{k}}{n}-\frac12\Big(\frac{n-\hat V_{k}}{n}\Big)^{2} + \eta\Big(\frac{n-\hat V_{k}}{n}\Big)+\sqrt\theta p(n-\hat V_{k})\\
& = \Big(\sqrt\theta p-\frac1n\Big) (k+\delta_{k}) - \frac12\Big(\frac{k+\delta_{k}}{n}\Big)^{2} + \eta\Big(\frac{k+\delta_{k}}{n}\Big)
\end{align*}
By a similar reasoning as before, we deduce that
\[
\sum_{k\le tn^{2/3}}\Big\{\log\frac{\hat V_{k}}{n}+ \theta^{\frac12} p(n-\hat V_{k})\Big\}\mathbf 1_{A_k}\quad\longrightarrow\quad  \frac12 \lambda t^{2}-\frac16 t^{3} \ \text{ in probability}. 
\]
Together with \eqref{eq: cvIk}, this proves the desired convergence of $\sum_{k\le N}I_{k}\mathbf 1_{A_k}$. 
\end{proof}

\begin{proof}[Proof of Lemma~\ref{lem: pfcv2}]
We note that $\tilde\eta(x):=  \log(1-x)+x $ is bounded by $C'x^{2}$ for all $x\in [0, 1-q]$ and some suitable constant $C'\in (0, \infty)$. Recall from Lemma \ref{lem: hatUk} that $m-\hat U_{k}= kmp + \epsilon_{k}$ and $n-\hat V_k=k+\delta_k$. We then have 
\begin{align*}
H_{k} & = (\Delta \hat R\n_{k}-\sqrt\theta) \cdot \Big\{-\frac{m-\hat U_{k}}{m}+\tilde\eta\Big(\frac{m-\hat U_{k}}{m}\Big) + (n-\hat V_{k})p\Big\}\\
 & = (\Delta \hat R\n_{k}-\sqrt\theta) \cdot \Big\{-\frac{kmp+\epsilon_{k}}{m} + \tilde\eta\Big(\frac{m-\hat U_{k}}{m}\Big) + (k+\delta_{k})p\Big\}\\
& = (\Delta \hat R\n_{k}-\sqrt\theta)\cdot\Big\{ -\frac{\epsilon_{k}}{m} + \tilde\eta\Big(\frac{m-\hat U_{k}}{m}\Big)+\delta_{k}p\Big\} =: (\Delta \hat R\n_{k}-\sqrt\theta)\cdot \gamma_{k}\,.   
\end{align*}
From the bound on $\tilde\eta(\cdot)$, it follows that we can find some positive constant $C''$ so that
\[
\gamma_{k}^{2}\,\mathbf 1_{A_{k}}\le C''\Big(\Big(\frac{\epsilon_{k}}{m}\Big)^{2} + k^{4}p^{4}+\delta_{k}^{2}p^{2}\Big). 
\]
Together with \eqref{bd: rebduk}, this implies that
\begin{equation}
\label{l2bd}
\sum_{k\le tn^{2/3}}\mathbb E\big[\gamma_{k}^{2}\mathbf 1_{A_k}\big]  = \cO(n^{-\frac23}) \quad\text{and}\quad \sum_{k\le tn^{2/3}}\sqrt{\mathbb E\big[\gamma_{k}^{2}\mathbf 1_{A_k}\big] } = \cO(1 ). 
\end{equation}
Let us denote 
\[
H_{k, 1} = (\Delta \hat R\n_k-\alpha_n) \cdot \gamma_k, \quad H_{k, 2}=(\alpha_n-\sqrt\theta) \gamma_k,
\]
so that $H_k = H_{k, 1}+H_{k, 2}$. Since $\Delta\hat R\n_k$ has expectation $\alpha_n$ and  is independent of $\mathcal G^n_{k-1}=\sigma((\hat R\n_j, \hat S\n_j)_{j\le k-1})$, while both $\gamma_k$ and $A_k$ are measurable in $\mathcal G^n_{k-1}$, we see that 
\[
\sum_{k\le j} H_{k, 1}\mathbf 1_{A_k}, \quad j\ge 1
\]
is a zero-mean martingale. Therefore, 
\[
\mathbb E\Big[\Big(\sum_{k\le N} H_{k, 1}\mathbf 1_{A_k}\Big)^2\Big] = \sum_{k\le N} \mathbb E\big[H_{k, 1}^2\mathbf 1_{A_k}\big] = \Var\big( \hat R\n_1\big) \sum_{k\le N} \mathbb E\big[\gamma_k^2\mathbf 1_{A_k}\big] = \cO(n^{-\frac23}), 
\]
where we have used~\eqref{l2bd} and the fact $\Var(\hat R\n_1)=\alpha_n=\cO(1)$ from~\eqref{eq: moments}. 
Meanwhile, with the Cauchy--Schwarz inequality, we find that 
\[
\mathbb E\Big[\sum_{k\le N}|H_{k, 2}|\mathbf 1_{A_k}\Big] \le |\alpha_n-\sqrt\theta|\cdot \sum_{k\le N} \sqrt{\mathbb E\big[\gamma_k^2\mathbf 1_{A_k}\big]} = \cO(n^{-\frac13}),
\]
where we have relied on~\eqref{l2bd} and the bound $\alpha_n-\sqrt\theta =\cO(n^{-1/3})$ obtained from the definition of $\alpha_n$.
To sum up, we have shown that 
\[
\sum_{k\le N} H_{k}\mathbf 1_{A_k}=\sum_{k\le N}H_{k, 1}\mathbf 1_{A_k}+\sum_{k\le N} H_{k, 2}\mathbf 1_{A_k}\to 0
\]
in probability. 
\end{proof}

\begin{proof}[Proof of Lemma~\ref{lem: pfcv3}]
We first note that by definition,
\[
\frac{mp-\alpha_{n} }{mp}= \lambda_{\ast}n^{-\frac13}-\lambda\lambda_{\ast}n^{-\frac23}+o(n^{-\frac23}), \quad \frac{np-\beta_{n}}{np} = \lambda_{\ast}n^{-\frac13}-\lambda\lambda_{\ast}n^{-\frac23}+o(n^{-\frac23}).
\]
Recall that $\eta(x)=\log(1-x)+x+x^2/2$, which is bounded by $Cx^3$ for $x\in [0, 1-q]$. 
When $n$ is sufficiently large, this leads to:
\begin{align*}
\alpha_{n}-mp+\theta^{\frac12}\log\frac{mp}{\alpha_{n}} & = \alpha_{n}-mp +\theta^{\frac12}\Big\{\frac{mp-\alpha_{n}}{mp}+\frac12\Big(\frac{mp-\alpha_{n}}{mp}\Big)^{2}-\eta\Big(\frac{mp-\alpha_{n}}{mp}\Big)\Big\}\\
& = -\theta^{\frac12}\lambda \lambda_{\ast} n^{-\frac23}+\tfrac{\sqrt\theta}{2}\lambda_{\ast}^{2}n^{-\frac23}+o(n^{-\frac23}),
\end{align*} 
With a similar calculation, we find that
\[
\theta^{\frac12}(\beta_{n}-np)+\log\frac{np}{\beta_{n}} = -\lambda\lambda_{\ast}n^{-\frac23}+\tfrac12\lambda_{\ast}^{2}n^{-\frac23}+o(n^{-\frac23}). 
\]
Combining this with~\eqref{bd: sumAc}, we deduce that 
\[
\Big|\sum_{k\le N}J_k\mathbf 1_{A_k^c}\Big| \le \cO(n^{-\frac23})\cdot  \sum_{k\le N} \mathbf 1_{A_k^c} \to 0 \quad\text{in probability}. 
\]
Therefore, 
\[
\sum_{k\le N}J_k\mathbf 1_{A_k} = \sum_{k\le N} J_k-\sum_{k\le N}J_k\mathbf 1_{A_k^c}= N\cdot \Big((1+\sqrt\theta)(\tfrac12\lambda_{\ast}-\lambda)\lambda_{\ast}n^{-\frac23}+o(n^{-\frac23})\Big)-\sum_{k\le N}J_k\mathbf 1_{A_k^c},
\]
which yields the desired result. 
\end{proof}

\begin{proof}[Proof of Lemma~\ref{lem: pfcv4}]
Recall that $n-\hat V_k=k+\delta_k$ and $\tilde\eta(x)=\log(1-x)+x$. Let us first show 
\[
\sum_{k\le tn^{2/3}} \Delta \hat S\n_{k}\cdot \Big\{-\frac{\delta_{k}}{n}+\tilde\eta\Big(\frac{k+\delta_{k}}{n}\Big)\Big\}\mathbf 1_{A_k} \to 0 \quad\text{in probability}.
\]
This follows from arguments similar to the ones in the proof of Lemma~\ref{lem: pfcv2}. Indeed, we can find some $C''\in (0, \infty)$ so that 
\[
\tilde\gamma_k^2\mathbf 1_{A_k}\le C''\Big(\frac{\delta^2_k}{n^2}+\frac{k^4}{n^4}\Big)  \quad\text{where } \tilde\gamma_k= -\frac{\delta_{k}}{n}+\tilde\eta\Big(\frac{k+\delta_{k}}{n}\Big).
\]
Combined with~\eqref{bd: rebduk}, this yields 
\[
\sum_{k\le N} \mathbb E[\tilde\gamma_k^2\mathbf 1_{A_k}] = \cO(n^{-\frac23}) \quad\text{and}\quad \sum_{k\le N}\sqrt{\mathbb E\big[\tilde\gamma_{k}^{2}\mathbf 1_{A_k}\big] } = \cO(1 ). 
\]
As in the proof of Lemma~\ref{lem: pfcv2}, we further write
\[
\sum_{k\le N} \Delta \hat S\n_{k}\cdot\tilde\gamma_k\mathbf 1_{A_k} = \sum_{k\le N} \big(\Delta\hat S\n_k-\mathbb E[\hat S\n_1]\big) \cdot \tilde\gamma_k\mathbf 1_{A_k} + \mathbb E[\hat S\n_1]\sum_{k\le N} \tilde \gamma_k\mathbf 1_{A_k}. 
\]
Arguing in the same way as in the previous case, we can show that the first sum is a martingale whose $L^2$-norm is $\mathcal O(n^{-\frac23})$, while the $L^1$-norm of the second sum is $\cO(n^{-\frac13})$ since $\mathbb E[\hat S\n_1]=\cO(n^{-\frac13})$ as seen in~\eqref{eq: moments}. This proves the claimed convergence in probability. Meanwhile, using the independence of $\Delta \hat S\n_k$ and $A_k$, we obtain that  
\begin{equation}
\label{bd: ant}
\mathbb E\Big[\sum_{k\le N} \frac{k}{n}\cdot |\Delta \hat S\n_k| \cdot \mathbf 1_{A_k^c}\Big] =\mathbb E\big[|\hat S\n_1| \big] \cdot \sum_{k\le N}\frac{k}{n}\cdot\mathbb P(A_k^c) \le \mathbb E\big[|\hat S\n_1| \big]\cdot \frac{N}{n} \sum_{k\le N}\mathbb P(A_k^c) \to 0,
\end{equation}
as a combined consequence of~\eqref{bd: sumAc} and the bound $\mathbb E[|\hat S\n_1|]\le \sqrt{\mathbb E[(\hat S\n_1)^2]}=\cO(1)$. 
Since 
\[
\sum_{k\le tn^{2/3}} L_{k}\mathbf 1_{A_k} = -\sum_{k\le tn^{2/3}} \frac{k}{n}\cdot \Delta\hat S\n_{k} \mathbf 1_{A_k}+ \sum_{k\le tn^{2/3}} \Delta \hat S\n_{k}\cdot \Big\{-\frac{\delta_{k}}{n}+\tilde\eta\Big(\frac{k+\delta_{k}}{n}\Big)\Big\}\mathbf 1_{A_k},
\]
combined with \eqref{id: ipp}, the previous arguments imply that jointly with the convergence in \eqref{cv: moderate-tree}, we have
\begin{equation}
\label{eq: cvLk}
\sum_{k\le tn^{2/3}} L_{k}\mathbf 1_{A_k}  \quad \Longrightarrow \quad -\int_{0}^{t} s\, d\hat\cS_{s}. 
\end{equation}
We turn to $F_{k}$. On the one hand, we have
\[
\Big(\lambda_{\ast}\big(1-\tfrac{1}{\sqrt\theta}\big)n^{-\frac13}\!\!\!\sum_{k\le tn^{2/3}}(\Delta\hat R\n_{k}-\theta^{\frac12}), \, \lambda_{\ast}n^{-\frac13}\!\!\sum_{k\le tn^{2/3}}\Delta\hat S\n_{k}\Big) \ \Longrightarrow\ \Big(\lambda_{\ast}\big(1-\tfrac{1}{\sqrt\theta}\big)\hat\cR_{t}, \;\lambda_{\ast} \hat\cS_{t}\Big),
\]
Moreover, this can be shown to hold jointly with \eqref{cv: moderate-tree} and \eqref{eq: cvLk}; see \eqref{id: ipp}. Arguing exactly as in~\eqref{bd: ant}, we have 
\[
n^{-\frac13}\mathbb E\Big[\sum_{k\le N}|\Delta \hat R\n_k-\sqrt\theta| \mathbf 1_{A_k^c}\Big] + n^{-\frac13}\mathbb E\Big[\sum_{k\le N}|\Delta\hat S\n_k|\mathbf 1_{A_k^c}\Big]\to 0. 
\]
Therefore, we can improve the previous convergence to 
\[
\Big(\lambda_{\ast}\big(1-\tfrac{1}{\sqrt\theta}\big)n^{-\frac13}\!\!\!\sum_{k\le tn^{2/3}}(\Delta\hat R\n_{k}-\theta^{\frac12})\mathbf 1_{A_k}, \, \lambda_{\ast}n^{-\frac13}\!\!\sum_{k\le tn^{2/3}}\Delta\hat S\n_{k}\mathbf 1_{A_k}\Big) \ \Longrightarrow\ \Big(\lambda_{\ast}\big(1-\tfrac{1}{\sqrt\theta}\big)\hat\cR_{t}, \;\lambda_{\ast} \hat\cS_{t}\Big),
\]
jointly with \eqref{cv: moderate-tree} and \eqref{eq: cvLk}. 
For the remaining terms in $F_{k}$, we observe that
\[
\hat\gamma:=\log\frac{mp}{\alpha_{n}}+\beta_{n}-np-\lambda_{\ast}(1-\theta^{-\frac12})n^{-\frac13}=\cO(n^{-\frac23}). 
\]
As $\Delta\hat R\n_k, k\ge 1$ are i.i.d, it follows from~\eqref{eq: moments} that 
\[
\mathbb E\Big[\Big(\sum_{k\le N}(\Delta\hat R\n_k-\sqrt\theta)\cdot \hat\gamma\cdot \mathbf 1_{A_k}\Big)^2\Big] \le \hat\gamma^2\cdot \Big(N \cdot\mathbb E\big[(\hat R\n_1-\sqrt\theta)^2\big]+N^2\big(\mathbb E[\hat R\n_1-\sqrt\theta]\big)^2\Big)=\cO(n^{-\frac23}).
\]
A similar argument leads to
\[
\mathbb E\Big[\Big(\sum_{k\le N}\Delta\hat S\n_k\Big(\log\frac{np}{\beta_n}-\lambda_{\ast}n^{-\frac13}\Big)\mathbf 1_{A_k}\Big)^2\Big]=\cO(n^{-\frac23}). 
\]
Since
\begin{align*}
F_k \mathbf 1_{A_k}& =  \lambda_{\ast}(1-\theta^{-\frac12})n^{-\frac13}(\Delta\hat R\n_k-\sqrt\theta)\mathbf 1_{A_k}+\lambda_{\ast}n^{-\frac13}\Delta\hat S\n_k\mathbf 1_{A_k}\\
& \quad +(\Delta\hat R\n_k-\sqrt\theta)\cdot\hat\gamma\mathbf 1_{A_k}+\Delta\hat S\n_k\Big(\log\frac{np}{\beta_n}-\lambda_{\ast}n^{-\frac13}\Big)\mathbf 1_{A_k},
\end{align*}
the claimed convergence concerning $F_k$ follows. Together with~\eqref{eq: cvLk}, this completes the proof. 
\end{proof}

\begin{proof}[Proof of Proposition~\ref{prop: cvdens}]
The convergence in~\eqref{cv: sumQ} combined with Lemmas~\ref{lem: pfcv1}-\ref{lem: pfcv4} entails the convergence in~\eqref{cv: dens''}. By Lemma 4.8 in~\cite{CG23}, this allows us to conclude. 
\end{proof}

\subsection{Convergence of the graph exploration processes in the light clustering regime}
\label{sec: light}

Our aim here is to prove Proposition~\ref{thm: cv-light}, whose proof is similar but simpler than that of Proposition~\ref{thm: cv-moderate}, its analogue for the moderate regime. Hence, we only outline the main steps and highlight the differences from the previous proof. The first step is identical to that in the moderate clustering regime. We introduce an approximation of $(R\n_k, S\n_k)$ using Poisson increments as in Section~\ref{sec: pf}: $\check X\n_k:=\check R\n_{k}-\check R\n_{k-1}$ has distribution Poisson$(\check U_kp)$ with $\check U_{k}$ defined in the same way as in~\eqref{def: checkUk}, while $\Delta \check S\n_k+1$ has distribution Poisson $(\check X\n_k\check V_k p)$ with $\check V_{k}$ defined as in~\eqref{def: checkUk}.
We then set  $\check H\n_k=\mathscr H_k((\check S\n_j)_{j\le k})$, $k\ge 1$. 
Note that Proposition~\ref{prop: tv} also holds under the assumption \eqref{hyp: L}, which tells us that we can replace $(R\n_k, S\n_k, H\n_k)$ with $(\check R\n_k, \check S\n_k, \check H\n_k)$ for the proof of Proposition~\ref{thm: cv-light}.

\subsubsection{Convergence of the random walks with i.i.d.~increments}
\label{sec: cv-rw-light}

We fix $\lambda_\ast\in (\lambda_{+}, \infty)$ and set
\begin{equation}
\label{def: betan}
\alpha_n = mp, \quad \beta_n = np-2\lambda_{\ast}m^{-\frac12}n^{\frac16}. 
\end{equation}
We observe that under \eqref{hyp: L}, $\alpha_{n}\to\infty$ and $\beta_{n}\to 0$; moreover, $m^{-\frac12}n^{\frac16}=o(np)$, and therefore we can find some $n_{0}\in \N$ so that 
$\beta_n>0$ for $n\ge n_0$. From now on, we restrict our discussion to $n\ge n_0$.  For those $n$, let $(\hat X\n_k)_{k\ge 1}$ be a sequence of i.i.d.~Poisson$(mp)$ variables, and let $\Delta\hat S\n_k+1, k\ge 1$, be independent with respective distributions Poisson$(\hat X\n_k \cdot \beta_n)$. We define
\[
\hat R\n_{k}=\sum_{1\le j\le k}\hat X\n_{k}, \quad \hat S\n_{k}=\sum_{1\le j\le k} \Delta \hat S\n_{k}, \quad \text{and}\quad  \hat H\n_{k}=\mathscr H_{k}\Big((\hat S\n_j)_{0\le j\le k}\Big), 
\]
with the functional $\mathscr H_{k}$ from \eqref{def: ht-proc}. 

\begin{prop}
\label{prop: cv-light-tree}
Assume~\eqref{hyp: L}. The following convergence in probability takes place  uniformly on compacts: 
\begin{equation}
\label{cv: R-light}
\Big\{n^{-\frac16}m^{-\frac12}\hat R\n_{\lf n^{2/3}t\rf} : t\ge 0\Big\} \xrightarrow[]{n\to\infty} \{t: t\ge 0\}. 
\end{equation}
Moreover, the following weak convergence takes place in $\mathbb D^{\ast}(\R_+, \R^{2})$: 
\begin{equation}
\label{cv: light-tree'}
\Big\{n^{-\frac13}\hat S\n_{\lf n^{2/3}t\rf}, n^{-\frac13}\hat H\n_{\lfloor n^{2/3}t\rfloor}: t\ge 0\Big\} \Longrightarrow \Big\{\hat{\cS}_{t}, \hat{\cH}_{t}: t\ge 0\Big\}, 
\end{equation}
where 
\[
\hat{\cS}_{t}  = W_{t}+2(\lambda-\lambda_\ast) t
\quad \text{and}\quad  \hat{\cH}_{t}  = 2\Big(\hat{\cS}_{t} - \min_{s\le t}\hat{\cS}_{s}\Big). 
\]
\end{prop}

\begin{proof}
We first show the convergence of $(\hat R\n_k)_{k\ge 0}$. This is a random walk with independent and identical increments satisfying
\[
\mathbb E[\hat X\n_k] =\Var(\hat X\n_k) = mp = m^{\frac12}n^{-\frac12}(1+\lambda n^{-\frac13}), \quad k\ge 1.
\]
Doob's maximal inequality for martingales then yields that for any $t>0$,
\[
\mathbb E\Big[\max_{s\le t}\Big(\hat R\n_{\lf sn^
{2/3}\rf}-\lf sn^{\frac23}\rf\cdot mp\Big)^2\Big] \le 4 tn^{\frac23}mp,
\]
from which the convergence in \eqref{cv: R-light} follows, since $n^{2/3}mp\sim n^{1/6}m^{1/2}$. 
Regarding $(\hat S\n_k)_{k\ge 0}$, we have
\begin{equation}
\label{eq: moment-light}
\mathbb E[\Delta\hat S\n_k] =\alpha_n\beta_n-1= 2(\lambda-\lambda_\ast) n^{-\frac13}+o(n^{-\frac13}), \quad \Var(\Delta\hat S\n_k) =\alpha_{n}\beta_n+\alpha_{n}\beta_n^2=1+o(1).
\end{equation}
Following the same calculations as in the moderate clustering regime, we have
\[
\mathbb E\big[(\hat S\n_1+1)^4\big] = \sum_{j=1}^4 \sum_{i=1}^{j}{4 \brace j}{j \brace i}\beta_n^j \alpha_{n}^{i}.
\]
Since $\sup_{n}\beta_{n}^{j}\alpha_{n}^{i}<\infty$ for all $j\ge i\ge 0$ under~\eqref{hyp: L}, we still have 
\begin{equation}
\label{bd: fourth}
\sup_{n}\mathbb E\big[(\hat S\n_1+1)^4\big]<\infty. 
\end{equation}
Therefore, we can argue in a similar way as in the moderate clustering regime and conclude that 
\[
\Big\{n^{-\frac13}\hat S\n_{\lf n^{2/3}t\rf}: t\ge 0\Big\} \Longrightarrow \Big\{  \hat{\cS}_{t}: t\ge 0\Big\} \quad \text{in $\mathbb D^{\ast}(\R_+, \R)$. }
\]
The rest of the proof is identical to those in the moderate clustering regime. We omit the details.
\end{proof}

\subsubsection{Convergence of the Radon--Nikodym derivatives}

The identity~\eqref{id: dens-dis} still holds in the light clustering regime, as its proof arguments still work in this case. 
For its continuum analogue, 
since the limit process of $(n^{-\frac16}m^{-\frac12}\hat R\n_k)_{k\ge 1}$ is deterministic (see \eqref{cv: R-light}), the analogue of \eqref{id: abs-cont} only involves $(\hat\cS_{t})_{t\ge 0}$. More precisely, for $t>0$ and a measurable function $F: \mathbb D^{\ast}([0, t], \R)\to \R_+$, we have 
\begin{equation}
\label{id: abs-cont-light}
\mathbb E\Big[F\big((\cS^{\lambda,\infty}_{s})_{s\le t}\big)\Big] = \mathbb E\Big[F\big((\hat\cS_s)_{s\le t}\big)\cdot \cE_{t}\big((\hat\cS_{s})_{s\le t}\big)\Big],
\end{equation}
where 
\begin{equation}
\label{def: Et'}
\mathcal E_{t}\big((\hat\cS_{s})_{s\le t}\big) = \exp\Big(-\int_{0}^{t} s\, d\hat\cS_{s}+2\lambda_{\ast}\hat\cS_t+(2\lambda_\ast^2-4\lambda\lambda_\ast)t+\lambda t^2-\frac16 t^3\Big). 
\end{equation}
The proof of \eqref{id: abs-cont-light} is a simpler version of the arguments presented in Section~\ref{sec: Girsanov}; we therefore leave the details to the reader. 
We also require a version of Lemma~\ref{lem: hatUk} in the light clustering regime. 
Note that $0\le m-\hat U_k\le \hat R\n_{k-1}$. 
The rest of the proof is a simple adaptation of the proof of Lemma~\ref{lem: hatUk} and therefore omitted. 
\begin{lem}
\label{lem: hatUk'}
Assume \eqref{hyp: L}. Let $\hat U_{k}, \hat V_{k}$ be as in \eqref{def: hatUk}. For all $t>0$, we have
\begin{equation}
\label{bd: hatUk'}
\max_{k\le tn^{2/3}}\mathbb E\big[(m-\hat U_{k})^{2}\big]=\cO(mn^{\frac13}), \quad \max_{k\le tn^{2/3}}\mathbb E\big[(n-\hat V_{k})^{2}\big]=\cO(n^{\frac43}).
\end{equation}
Let $\delta_{k}=n-\hat V_{k}-k$, $k\ge 1$. For all $t>0$, we also have 
\begin{equation}
\label{bd: rebduk'}
 \max_{k\le tn^{2/3}}\mathbb E\big[\delta_{k}^{2}\big]=\cO(n^{\frac23}).
\end{equation}
\end{lem}

\begin{prop}
\label{prop: cv-dens-light}
Assume~\eqref{hyp: L}. Let $E\n_{N}$ be as in \eqref{def: En}. For all $t\ge 0$,  the sequence 
\[
\Big\{E\n_{\lf tn^{2/3}\rf}\Big((\hat R\n_{k}, \hat S\n_{k})_{0\le k\le tn^{2/3}}\Big): n\ge 1\Big\}
\]
 is uniformly integrable.  Moreover, jointly with the convergences in \eqref{cv: R-light} and \eqref{cv: light-tree'}, for any $t\ge 0$, 
\begin{equation}
\label{eq: cvdens-light}
E\n_{\lf tn^{2/3}\rf}\Big((\hat R\n_{k}, \hat S\n_{k})_{0\le k\le tn^{2/3}}\Big) \ \Longrightarrow\  \mathcal E_{t}\big((\hat\cS_{s})_{0\le s\le t}\big),
\end{equation}
where $\mathcal E_{t}$ is defined in~\eqref{def: Et'}. 
\end{prop}

\begin{proof}
Denote $N=\lfloor tn^{2/3}\rfloor$ and fix $q\in (0, 1)$. 
Define $\sigma_n=\inf\{k\ge 1: n-\hat V_k\ge (1-q)n\}$ and $A_k = \{ n-\hat V_k\le (1-q)n\}$, so that $\{\sigma_n>N\}=\cap_{1\le k\le N} A_k$. By the same arguments in the proof of Proposition~\ref{prop: cvdens}, it suffices to show that 
\begin{equation}
\label{eq: cvdens-light'}
E\n_{N}\Big((\hat R\n_{k}, \hat S\n_{k})_{0\le k\le N}\Big)\mathbf 1_{\{\sigma_n>N\}} \ \Longrightarrow\  \mathcal E_{t}\big((\hat\cS_{s})_{0\le s\le t}\big), 
\end{equation}
jointly with the convergences in \eqref{cv: R-light} and \eqref{cv: light-tree'}. To that end, let us write
\[
E\n_{N}\big((\hat R\n_{k}, \hat S\n_{k})_{1\le k\le N}\big)  = \exp\bigg\{\sum_{1\le k\le N}\Big\{(F_k+H_{k}+I_{k}+J_k+L_{k})\mathbf 1_{A_k}+Q_k\Big\}\bigg\},
\]
where $Q_k=-\infty\cdot \mathbf 1_{A_k^c}$ and
\begin{align*}
H_{k} & = \Delta \hat R\n_{k}\log\frac{\hat U_{k}}{m}+(\Delta \hat R\n_{k}-mp)(\beta_n-\hat V_{k}p)+(m-\hat U_{k})p,\quad I_k = \log\frac{\hat V_k}{n}+mp^2(n-\hat V_k)\\
J_{k} & =  \log\frac{np}{\beta_n} + m p(\beta_n-np), \quad L_{k} = \Delta \hat S\n_{k} \cdot \log\frac{\hat V_{k}}{n}, \quad F_k=\Delta \hat S\n_{k} \cdot \log\frac{np}{\beta_n}. 
\end{align*}
Equipped with the estimates of $\hat U_{k}, \hat V_{k}$ from Lemma~\ref{lem: hatUk'}, we will show that the following convergences hold in probability:
\[
\sum_{k\le N}Q_k\to 0, \quad \sum_{k\le N}H_{k}\mathbf 1_{A_k}\to 0, \quad \sum_{k\le N} I_{k}\mathbf 1_{A_k}\to -\frac16 t^{3}+\lambda t^{2}, \quad \sum_{k\le N}J_{k}\mathbf 1_{A_k}\to (2\lambda_{\ast}^2-4\lambda\lambda_{\ast})t. 
\]
Meanwhile, the following weak convergences hold jointly with the convergences in \eqref{cv: R-light} and \eqref{cv: light-tree'}: 
\[
\Big(\sum_{k\le N}L_{k}\mathbf 1_{A_k}, \sum_{k\le N}F_k\mathbf 1_{A_k}\Big) \ \Longrightarrow\  \Big(-\int_{0}^{t}s\, d\hat\cS_{s}, \; 2\lambda_{\ast}\hat\cS_t\Big). 
\]
We start with $Q_{k}$. Following the same arguments as in the moderate clustering regime,  we deduce from Markov's inequality and~\eqref{bd: hatUk'} that 
\begin{equation}
\label{bd: goodeve}
\mathbb P(\sigma_n\le N)\le \max_{k\le N} \frac{\mathbb E[(n-\hat V_k)^2]}{(1-q)^2n^2} = \cO(n^{-\frac23}). 
\end{equation}
Since $\sum_k Q_k\ne 0$ if and only if $\sigma_n\le N$, the convergence of $\sum_k Q_k$ follows.  We turn to $H_k$, which we split into the following two terms: 
\[
H_{k, 1} := (\Delta \hat R\n_k-mp)\cdot \gamma_k, \quad H_{k, 2} = mp\cdot \log \frac{\hat U_k}{m} + (m-\hat U_k)p,
\]
with 
\[
\gamma_k := \log \frac{\hat U_k}{m}+\beta_n-\hat V_k p= \log\Big(1-\frac{m-\hat U_k}{m}\Big)+(n-\hat V_k)p + \beta_n-np 
\]
Recall that we can then find some $C'\in (0, \infty)$ so that  $|\log(1-x)+x| \le C' x^{2}$ for all $x\in [0, 1-q]$. Combined with~\eqref{bd: hatUk'} and~\eqref{def: betan}, this implies that there exists some $C''\in (0, \infty)$ so that 
\[
\gamma_k^2\mathbf 1_{A_k} \le C''\Big(\Big(\frac{m-\hat U_k}{m}\Big)^2+(n-\hat V_k)^2p^2+(\beta_n-np)^2\Big), \quad\text{hence}\;  \mathbb E[\gamma_k^2\mathbf 1_{A_k}] = \cO(m^{-1}n^{\frac13}). 
\]
Arguing in the same way as in the previous case, we can show that $\sum_k H_{k, 1}\mathbf 1_{A_k}$ is a zero-mean martingale satisfying
\[
\mathbb E\Big[\Big(\sum_{k\le N}H_{k, 1}\mathbf 1_{A_k}\Big)^2\Big] = \Var(\hat R\n_1) \sum_{k\le N} \mathbb E[\gamma_k^2\mathbf 1_{A_k}] = mp\cdot \cO(n^{\frac23}m^{-1}n^{\frac13})=\cO(np)\to 0,
\]
under~\eqref{hyp: L}. For $H_{k, 2}$, we find from applying the previous bound on $\log(1-x)$ that 
\[
|H_{k, 2}|\mathbf 1_{A_k} \le C'mp\cdot \Big(\frac{m-\hat U_k}{m}\Big)^2, \quad\text{hence}\; \sum_{k\le N} \mathbb E[|H_{k, 2}|\mathbf 1_{A_k}]=\cO(np),
\]
by employing the bounds in~\eqref{bd: hatUk'} again. Combining the bounds from the  previous two displays allows us to conclude that $\sum_k H_k\mathbf 1_{A_k}\to 0$ in probability. 
Let us show the convergence concerning $I_{k}$. We recall  that $\delta_{k}= n-\hat V_{k}-k$ and $\eta(x)=\log(1-x)+x+x^{2}/2$, with the bound $|\eta(x)|\le Cx^{3}$ for all $x\in [0, 1-q]$ with some suitable $C>0$.
We split $I_{k}$ as follows:
\[
I_{k}=k\Big\{mp^{2}-\frac{1}{n}\Big\} + \delta_{k}\Big\{mp^{2}-\frac1n\Big\} -\frac12 \Big(\frac{k}{n}\Big)^{2} +\Big\{\frac12\Big(\frac{k}{n}\Big)^{2}-\frac12\Big(\frac{k+\delta_{k}}{n}\Big)^{2}+\eta\Big(\frac{n-\hat V_{k}}{n}\Big)\Big\}.
\]
Assumption \eqref{hyp: L} and~\eqref{id: sum} imply that
\begin{equation}
\label{sum: ik}
n^{\frac43}\Big\{mp^{2}-\frac1n\Big\} \to  2\lambda \quad \text{so that} \quad \sum_{k\le tn^{2/3}}k\Big\{mp^{2}-\frac{1}{n}\Big\}\to \lambda t^{2}.
\end{equation}
Meanwhile, Lemma~\ref{lem: hatUk'} tells us that
\[
\mathbb E\Big[\Big|\sum_{k\le N}\delta_{k}\Big(mp^{2}-\frac1n\Big)\mathbf 1_{A_k}\Big|\Big]\le \Big|mp^{2}-\frac1n\Big| \cdot \sum_{k\le N}\sqrt{\mathbb E\big[\delta_{k}^{2}\big] }= \cO(n^{-\frac13}). 
\]
With the previous bound of $\eta$, we deduce that 
\[
\Big|\frac12\Big(\frac{k}{n}\Big)^{2}-\frac12\Big(\frac{k+\delta_{k}}{n}\Big)^{2}+\eta\Big(\frac{n-\hat V_{k}}{n}\Big)\Big|\mathbf 1_{A_k}\le \frac{k|\delta_k|}{n^{2}}| + \frac{\delta_k^2}{n^{2}} +C\,\frac{(k+\delta_{k})^{3}}{n^{3}}
\]
An application of the Cauchy--Schwarz inequality and Lemma~\ref{lem: hatUk'} yields that  
\[
\sum_{k\le N}\mathbb E\Big[\Big|\frac12\Big(\frac{k}{n}\Big)^{2}-\frac12\Big(\frac{k+\delta_{k}}{n}\Big)^{2}+\eta\Big(\frac{n-\hat V_{k}}{n}\Big)\Big|\mathbf 1_{A_k}\Big]=\cO(n^{-\frac13}).
\]
As a final ingredient, we follow the same reasoning as in~\eqref{bd: sumAc} to find that 
\begin{equation}
\label{bd: sumAc'}
\mathbb E\Big[\sum_{k\le N}\mathbf 1_{A_k^c}\Big] \le (N+1)\cdot\mathbb P(\sigma_n\le N)=\cO(1),
\end{equation}
according to~\eqref{bd: goodeve}. Combining this with~\eqref{sum: ik}, we deduce that 
\[
\sum_{k\le N} \mathbb E\bigg[\bigg|k\Big\{mp^{2}-\frac{1}{n}\Big\} -\frac12 \Big(\frac{k}{n}\Big)^{2} \bigg|\mathbf 1_{A_k^c} \bigg]\le \cO(n^{-\frac23})\cdot  \mathbb E\Big[\sum_{k\le N}\mathbf 1_{A_k^c}\Big]\to 0.
\]
Together with $\sum_{k\le tn^{2/3}}(k/n)^{2}\to \frac13 t^{3}$ as seen in~\eqref{id: sum}, the previous arguments show that  $\sum_{k\le tn^{2/3}} I_{k}\mathbf 1_{A_k}\to \lambda t^{2}-\frac16 t^{3}$ in probability. 
Using the fact that 
\begin{equation}
\label{asy: beta}
np-\beta_{n} = 2\lambda_{\ast}m^{-\frac12}n^{\frac16}\quad \text{ and } \quad\frac{np-\beta_{n}}{np} = 2\lambda_{\ast}n^{-\frac13}-2\lambda_{\ast}\lambda n^{-\frac23}+o(n^{-\frac23}),
\end{equation}
we find after an elementary computation that $\sum_{k\le N}J_{k}\to (2\lambda_{\ast}^2-4\lambda\lambda_{\ast})t$. It follows $\sum_{k\le N}J_k\mathbf 1_{A_k}$ converging to the same limit, thanks to~\eqref{bd: sumAc'} and the bound $\max_k |J_k|=\cO(n^{-\frac13})$. 
 
The convergence of $\sum_{k\le N}L_{k}\mathbf 1_{A_k}$ can be shown with identical arguments as in the moderate clustering regime; we omit the details. 
Finally, let us write 
\begin{align*}
F_{k} & =\Delta\hat S\n_{k}\frac{np-\beta_{n}}{np} -\Delta \hat S\n_{k}\cdot \tilde\eta\Big(\frac{np-\beta_{n}}{np}\Big)\\
&=\Delta\hat S\n_{k}\cdot\big(2\lambda_{\ast}n^{-\frac13}-2\lambda_{\ast}\lambda n^{-\frac23}+o(n^{-\frac23})\big)-\Delta \hat S\n_{k}\cdot \tilde\eta\Big(\frac{np-\beta_{n}}{np}\Big),
\end{align*}
where we recall $\tilde \eta(x)=\log(1-x)+x$ is bounded by $C'x^{2}$ on $[0, 1-q]$. 
Noting that $\Delta\hat S\n_{k'}$ is independent of $\Delta\hat S\n_k, A_k, A_{k'}$ for $k<k'$, we deduce that 
\begin{align*}
\mathbb E\Big[\Big(\sum_{k\le N}n^{-\frac23}\Delta\hat S\n_{k}\mathbf 1_{A_k}\Big)^{2}\Big]&=n^{-\frac43}\sum_{k\le N}\mathbb E\big[(\Delta\hat S\n_k)^2\mathbf 1_{A_k}\big]+2n^{-\frac43}\sum_{k<k'}\mathbb E\big[\Delta\hat S\n_k\Delta\hat S\n_{k'}\mathbf 1_{A_k}\mathbf 1_{A_{k'}}\big]\\
&\le n^{-\frac43}N\cdot\mathbb E\big[(\hat S\n_{1})^{2}\big] +2n^{-\frac43}\sum_{k<k'} \mathbb E[\Delta\hat S\n_{k'}]\cdot \mathbb E\big[\Delta\hat S\n_k\mathbf 1_{A_k}\mathbf 1_{A_{k'}}\big]\\
&\le \cO(n^{-\frac23})+ 2n^{-\frac43}N\cdot\big|\mathbb E[\hat S\n_1]\big|\sum_{k\le N}  \sqrt{\mathbb E\big[(\Delta \hat S\n_k)^2\big]} = \cO(n^{-\frac13})
\end{align*}
according to~\eqref{eq: moment-light}. 
The bound for $\tilde\eta$ and~\eqref{asy: beta} implies that $\tilde\eta((np-\beta_n)/np)=\cO(n^{-\frac23})$. 
Then a similar calculation shows that 
\[
\mathbb E\bigg[\Big(\sum_{k\le N}\Delta \hat S\n_{k}\cdot \tilde\eta\Big(\frac{np-\beta_{n}}{np}\Big)\mathbf 1_{A_k}\Big)^2\bigg]=\cO(n^{-\frac13}). 
\]
Using the independence between $\Delta\hat S\n_k$ and $A_k$, we deduce that
\[
\mathbb E\Big[\sum_{k\le N}n^{-\frac13}|\Delta\hat S\n_k|\mathbf 1_{A_k^c}\Big] \le n^{-\frac13}\mathbb E\big[|\hat S\n_1|\big] \sum_{k\le N}\mathbb P(A_k^c) = \cO(n^{-\frac13})
\] 
by~\eqref{bd: sumAc'} and~\eqref{eq: moment-light}. 
Combining the previous arguments shows that $\sum_{k} F_{k}\mathbf 1_{A_k}$ converges in distribution to $2\lambda_{\ast}\hat\cS_{t}$, jointly with the convergence of $\sum_{k}L_{k}\mathbf 1_{A_k}$ and  \eqref{cv: R-light}, \eqref{cv: light-tree'}. 
This completes the proof of \eqref{eq: cvdens-light'}. Lemma 4.8 in~\cite{CG23} asserts that this is enough to conclude. 
\end{proof}

\subsubsection{Convergence of the graph exploration processes}

\begin{proof}[Proof of Proposition~\ref{thm: cv-light}]
Take an arbitrary $t\in (0, \infty)$. Let $F: \mathbb D^{\ast}([0, t], \R^{3})\to \R$ be continuous and bounded. Fix some $K\in (0, \infty)$. Writing $E\n_{N}$ as a shorthand for 
\[
E\n_{\lf tn^{2/3}\rf}\Big((\hat R\n_{k}, \hat S\n_{k})_{k\le tn^{2/3}}\Big), 
\]
we have
\begin{align*}
&\quad  \mathbb E\Big[F\Big(\Big(n^{-\frac16}m^{-\frac12}\check R\n_{\lf n^{2/3}s\rf}, n^{-\frac13}\check S\n_{\lf n^{2/3}s\rf}, n^{-\frac13}\check H\n_{\lf n^{2/3}s\rf}\Big)_{s\le t}\Big)\Big] \\
&=\mathbb E\Big[F\Big(\Big(n^{-\frac16}m^{-\frac12}\hat R\n_{\lf n^{2/3}s\rf}, n^{-\frac13}\hat S\n_{\lf n^{2/3}s\rf}, n^{-\frac13}\hat H\n_{\lf n^{2/3}s\rf}\Big)_{s\le t}\Big)\cdot E\n_N\mathbf 1_{\{E\n_N\le K\}}\Big]\\
& \quad +  \mathbb E\Big[F\Big(\Big(n^{-\frac16}m^{-\frac12}\hat R\n_{\lf n^{2/3}s\rf}, n^{-\frac13}\hat S\n_{\lf n^{2/3}s\rf}, n^{-\frac13}\hat H\n_{\lf n^{2/3}s\rf}\Big)_{s\le t}\Big)\cdot E\n_N\mathbf 1_{\{E\n_N> K\}}\Big].
\end{align*}
By choosing $K$ sufficiently large, the last term above does not exceed $\epsilon$ even as $n\to\infty$, due to the uniform integrability of $(E\n_{N})_{n\ge 1}$. Applying dominated convergence together with Proposition~\ref{prop: cv-light-tree}, Proposition~\ref{prop: cv-dens-light} and \eqref{id: abs-cont-light}, we find after taking $K\to\infty$ that 
\[
\mathbb E\Big[F\Big(\Big(n^{-\frac16}m^{-\frac12}\check R\n_{\lf n^{2/3}s\rf}, n^{-\frac13}\check S\n_{\lf n^{2/3}s\rf}, n^{-\frac13}\check H\n_{\lf n^{2/3}s\rf}\Big)_{s\le t}\Big)\Big] \to \mathbb E\Big[F\big(\big(s, \cS^{\lambda,\infty}_{s}, \cH^{\lambda,\infty}_{s}\big)_{s\le t}\big)\Big]. 
\]
Combined with Proposition~\ref{prop: tv}, this concludes the proof. 
\end{proof}

\subsection{Convergence of the bipartite surplus edges}
\label{sec: surplus}

The convergence of graph exploration processes in Propositions~\ref{thm: cv-moderate} and~\ref{thm: cv-light} ensures the convergence of the bipartite spanning tree $\cF$. To upgrade this into a convergence of the graph $B(n, m, p)$, we need to study the limit distribution of the surplus edges, which we do in this subsection. We first require some estimates. 

\subsubsection{Preliminary estimates}

Recall from Algorithm~\ref{alg} the vertex $v_{k}$ explored at step $k$ and the set $\cM_{k}=\{u_{k, i}: 1\le i\le X\n_{k}\}$. Recall that $\cB(w)$ stands for the neighbourhood of a vertex $w$ in $B(n, m, p)$.  Let
\begin{equation}
\label{def: goodevent}
\cI\n_{k} = \Big\{\cB(u_{k, i})\setminus\{v_{k}\}, 1\le i\le X\n_{k}, \text{ are pairwise disjoint}\Big\}, \quad k\ge 1. 
\end{equation}

\begin{lem}
\label{lem: goodevent}
Assume either~\eqref{hyp: L} or~\eqref{hyp: M}. We have
\begin{equation}
\label{bd: goodevent}
1-\mathbb P(\cI\n_k) = \cO(n^{-1}), \quad 1-\mathbb P\Big(\bigcap_{k\le tn^{2/3}}\cI\n_{k}\Big)=\cO(n^{-\frac13}).  
\end{equation}
\end{lem}

\begin{proof}
On the complementary of $\cI\n_{k}$, there must be some $v\ne v_{k}$ satisfying 
 $\#(\cB(v)\cap \cM_{k})\ge 2$. Since $\#(\cB(v)\cap \cM_{k})$ is stochastically bounded by a Binom$(X\n_{k}, p)$. Noting that $X\n_{k}$ is itself stochastically bounded by a Binom$(m, p)$ variable, we deduce from \eqref{bd: binom-tail} that
 \[
 1-\mathbb P(\cI\n_{k}) \le n\cdot \mathbb P(\text{Binom}(X\n_{k}, p)\ge 2) \le np^{2}\mathbb E[X\n_{k}] +np^{2}\mathbb E[(X\n_{k})^{2}]  \le 2mnp^{3}+m^{2}np^{4},
 \]
 which yields the first bound in \eqref{bd: goodevent}. 
The second bound in \eqref{bd: goodevent} then follows through a union bound. 
\end{proof}

Recall that $U_{k}=\#\cU_{k-1}, V_{k}=\#\cV_{k-1}$, $k\ge 1$. 
\begin{lem}
\label{lem: Uk}
Assume either \eqref{hyp: L} or \eqref{hyp: M}. For all $t>0$, we have
\begin{equation}
\label{bd: uk}
\max_{k\le tn^{2/3}}\mathbb E\big[(m-U_k)^{2}\big] = \cO(n^{\frac13}m ) , \quad \max_{k\le tn^{2/3}}\mathbb E\big[(n-V_{k})^{2}\big]=\cO(n^{\frac43}). 
\end{equation}
Using the shorthand notation $\mathbb E_{\cI_{t}}[\cdot]=\mathbb E[\cdot\,|\,\cap_{k\le tn^{2/3}}\cI\n_{k}]$, we also have
\begin{equation}
\label{bd: cuk}
\max_{k\le tn^{2/3}}\mathbb E_{\cI_{t}}\big[(m-U_k)^{2}\big] = \cO(n^{\frac13}m ) , \quad \max_{k\le tn^{2/3}}\mathbb E_{\cI_{t}}\big[(n-V_{k})^{2}\big]=\cO(n^{\frac43}). 
\end{equation}
\end{lem}

\begin{proof}
Let us denote $N=\lf tn^{2/3}\rf$. 
We first show the bound for $\max_{k\le N}(n-V_{k})=n-V_{N}$. 
By \eqref{id: UkVk}, we have $n-V_{N}\le N+S\n_{N-1}-\min_{k\le N-1}S\n_{k}\le 2N+S\n_{N-1}$, since $\Delta S\n_{k}\ge -1$ for all $k$. 
So it suffices to show that $\mathbb E[(S\n_{N-1})^{2}] = \cO(n^{4/3})$. To that end, let us introduce a random walk $(\mathrm T_{k})_{k\ge 0}$ which starts from $0$ and has i.i.d.~increments distributed as Binom$(Yn, p)$ with $Y$ following a Binom$(m, p)$ law. Clearly, $S\n_{N-1}$ is stochastically dominated by $\mathrm T_{N-1}$. An elementary calculation shows that $\mathbb E[\mathrm T_{1}]=mnp^{2}=\cO(1)$ and $\Var(\mathrm T_{1}) \le mnp^{2}+mn^{2}p^{3}=\cO(1)$ under both~\eqref{hyp: M} and~\eqref{hyp: L}. 
It follows that 
\[
\mathbb E\big[(S\n_{N-1})^{2}\big] \le \mathbb E\big[(\mathrm T_{N})^{2}\big] = N\cdot \mathbb E\big[(\mathrm T_{1})^{2}\big]+N(N-1) \big(\mathbb E[\mathrm T_{1}]\big)^{2} =\cO(n^{\frac43}). 
\]
Combined with the previous arguments, this yields the desired bound for $\mathbb E[(n-V_{N})^{2}]$. As for $m-U_{N}$, we deduce from~\eqref{id: UkVk} that $m-U_{N}=\sum_{k\le N-1}X\n_{k}$, where each $X\n_{k}$ is stochastically bounded by Binom$(m, p)$ for each $k\ge 1$. 
It follows that 
\[
\mathbb E\big[(m-U_{N})^2\big]  \le N(m^2p^2+mp) + N^{2}m^2p^2  = \cO(n^{\frac43}m^2p^2)= \cO(n^{\frac13}m),  
\]
since $mnp^2\to 1$. Finally, the bounds in~\eqref{bd: cuk} follows from~\eqref{bd: uk} and~\eqref{bd: goodevent}, since we have $\mathbb E_{\cI_{t}}[\cdot]\le \mathbb E[\cdot]/\mathbb P(\cap_{k\le tn^{2/3}}\cI\n_{k})$. 
\end{proof}

\subsubsection{Proof of Proposition~\ref{prop: surplus-bi}}

Here, we prove Proposition~\ref{prop: surplus-bi}. Before delving into {the details}, let us make a convenient assumption: by Skorokhod's Representation Theorem, we can assume that the convergences in \eqref{cv: moderate} and \eqref{cv: light} take place {\it almost surely} in the respective asymptotic regimes.  
In Lemma~\ref{lem: surplus} we have identified two potential scenarios where a bipartite surplus edge may appear. 
However, if there is a surplus edge between the members of $\cM_{k}$ and $\cN_{k}$, then the complement of $\cI\n_{k}$ occurs. In view of Lemma~\ref{lem: goodevent}, we can expect that only the surplus edges between $\cM_{k}$ and $\cA_{k-1}$ will make into the limit. We have the following description on their distributions.

\begin{lem}
\label{lem: dist-surplus}
Let $2\le k\le n$. Conditional on $(\cA_{j-1})_{j\le k}$ and $(\cM_{j})_{j\le k}$, the events 
\[
\{u, w\} \text{ is an edge of } B(n, m, p), \text{ with } u\in \cM_{k}, w\in \cA_{k-1}, 
\]
are  independent among themselves, and independent of $(\cN_{j})_{j\ge k}$ and $(\cM_{j})_{j>k}$. Moreover, each of these events occurs with probability $p$. 
\end{lem}

\begin{proof}
By the definition of $B(n, m, p)$, 
\[
\rB_{u, v}:=\mathbf 1_{\{\{u, v\}\text{ is an edge of $B(n, m, p)$}\}}, \quad u\in \cU, v\in \cV,
\]
is a collection of i.i.d.~Bernoulli variables with mean $p$. 
Let $\cF_{1}$ be the $\sigma$-algebra generated by $v_{1}, \cM_{1}, \cN_{1}$, $v_{2}$ and $\cM_{2}$. Let us observe that 
$\cM_{1}\in \sigma(\{v_{1}\})\vee\sigma(\{\rB_{u, v_{1}}: u\in \cU\})$, $\cN_{1}\in \sigma(\{\rB_{u, v}: u\in \cM_{1}, v\in \cV\setminus\{v_{1}\}\})$, and $\cM_{2}\in \sigma(\{v_{2}\})\vee\sigma(\{\rB_{u, v_{2}}: u\in \cU\setminus\cM_{1}\})$. We also have $\cA_{1}\in \cF_{1}$. 
Given $\cF_{1}$, $\{\rB_{u, v}: u\in \cM_{2}, w\in \cA_{1}\}$ is a collection of Bernoulli variables disjoint from the Bernoulli variables that have been used to generate $\cF_{1}$; as a result, it is distributed as a collection of independent Bernoulli variables with mean $p$. Moreover, the values of the variables in this collection do not affect $\cN_{2}$, nor any of the future $\cN_{j}$ and $\cM_{j}$, $j>2$. 
An inductive argument on $k$ then allows us to conclude that for any $k\ge 2$, given $\cF_{k}$ the $\sigma$-algebra generated by $(v_{j})_{j\le k}, (\cM_{j})_{j\le k}, (\cN_{j})_{j<k}$, 
the collection of $\rB_{v, w}$ that are responsible for producing surplus edges at step $k$ is disjoint from those generating $\cF_{k}$, as well as those generating any future $\cM_{j}, \cN_{j}, j>k$.  The claimed statement now follows. 
\end{proof}

Recall that $P\n$ is a simple point measure obtained from the point measure $\tilde P\n$ in \eqref{def: Pn}. Lemma~\ref{lem: dist-surplus}, combined with \eqref{id: Ak}, implies that for each $k\ge 1$ and $1\le l\le S\n_{k-1}-\min_{j\le k-1}S\n_{j}$, $P\n$ has an atom at $(n^{-2/3}k, n^{-1/3}l)$ with probability 
\[
q_{k, l}:=1-(1-p)^{\#\cM_{k}}, 
\]
independently of the other pairs. Let us set  
\[
D_{n}=\Big\{(k, l) :  k\in \N, l\in \N \text{ and } l \le S\n_{k-1}-\min_{j\le k-1}S\n_{j} \Big\} \ \text{and}\ \widetilde D_{n}=\Big\{(n^{-\frac23}k, n^{-\frac13}l): (k, l)\in D_{n}\Big\}.
\] 
Let $\mathrm{Leb}$ denote the Lebesgue measure on $\R^{2}$. 
\begin{lem}
\label{lem: cv-pp}
Fix $t, L\in (0, \infty)$. Let $\theta\in (0, \infty]$. We assume
\begin{itemize}
\item
either that \eqref{hyp: L} holds and \eqref{cv: light} takes place almost surely;
\item
or that \eqref{hyp: M} holds and \eqref{cv: moderate} takes place almost surely.
\end{itemize}
In both cases, we have 
\begin{equation}
\label{cv: P-rate}
-\sum_{k=1}^ {\lf tn^{2/3} \rf}\sum_{l=1}^{\lf Ln^{1/3}\rf} \log(1-q_{k, l})\mathbf 1_{\{(k, l)\in D_{n}\}} \xrightarrow{n\to\infty} \mathrm{Leb}\big(\cD^{\lambda, \theta}\cap [0, t]\times [0, L]\big)\quad\text{in probability},
\end{equation}
and 
\begin{equation}
\label{cv: P-bd}
\limsup_{n\to\infty}\sum_{(k, l)\in\N^{2}}\mathbb E\Big[q_{k, l}\mathbf 1_{\{(n^{-\frac23}k, n^{-\frac13}l)\in \widetilde D_{n}\cap C\}}\Big]\le \mathbb E\big[\mathrm{Leb}\big(\cD^{\lambda, \theta}\cap C\big)\big]
\end{equation}
for all compact sets $C\subseteq \R^{2}$. 
\end{lem}

\begin{proof}
The regularity of the Lebesgue measure implies that 
\[
\frac{1}{n}\#\Big\{(n^{-\frac23}k, n^{-\frac13}l): (k, l)\in \cD^{\lambda,\theta}\cap [0, t]\times [0, L]\Big\}\to \mathrm{Leb}\big(\cD^{\lambda,\theta}\cap [0, t]\times [0, L]\big)\quad \text{a.s.}
\]
as well as 
\[
\limsup_{n\to\infty}\frac{1}{n}\#\Big\{(n^{-\frac23}k, n^{-\frac13}l): (k, l)\in \cD^{\lambda,\theta}\cap C\Big\}\le \mathrm{Leb}\big(\cD^{\lambda,\theta}\cap C\big)\quad \text{a.s.}
\]
for a compact set $C\subseteq \R^{2}$. Together with the convergence in \eqref{cv: moderate} (resp.~in~\eqref{cv: light}), it follows that a.s.
\begin{equation}
\label{cv: volume}
\frac{1}{n}\#\Big(\widetilde D_{n}\cap [0,t]\times [0, L] \Big) \to \mathrm{Leb}\big(\cD^{\lambda, \theta}\cap [0, t]\times [0, L]\big), \;
\limsup_{n\to\infty}\frac{1}{n}\#\big(\widetilde D_{n}\cap C \big) \le \mathrm{Leb}\big(\cD^{\lambda,\theta}\cap C\big).
\end{equation}
Recall that $X\n_{k}=\#\cM_{k}$. We can split the left-hand side of \eqref{cv: P-rate} into the following sum:
\[
mp|\log(1-p)| \#\Big(\widetilde D_{n}\cap [0, t]\times [0, L]\Big) + A_{n}\big(\widetilde D_{n}\cap[0, t]\times [0, L]\big)+B_{n}\big(\widetilde D_{n}\cap [0, t]\times [0, L]\big),
\]
where for $D\subseteq [0, t]\times [0, L]$, we have denoted
\begin{align*}
A_{n}(D) & = p|\log(1-p)|\sum_{(k, l)\in \N^{2}}(U_{k}-m)\mathbf 1_{\{(n^{-\frac23}k, n^{-\frac13}l)\in D\}}, \\
B_{n}(D) & = |\log(1-p)| \sum_{(k, l)\in \N^{2}}\big(X\n_{k}-U_{k}p\big)\mathbf 1_{\{(n^{-\frac23}k, n^{-\frac13}l)\in D\}}.
\end{align*}
Noting that $\log(1-p)+p=o(p)$ as $n\to\infty$, we have 
\begin{align*}
mp|\log(1-p)| \cdot\#\Big(\widetilde D_{n}\cap [0, t]\times [0, L]\Big)  &= \frac{1}{n}(1+o(1))\cdot \#\Big(\widetilde D_{n}\cap [0, t]\times [0, L]\Big)\\
&\to \mathrm{Leb}\big(\cD^{\lambda, \theta}\cap [0, t]\times [0, L]\big)\quad \text{a.s.}
\end{align*}
\eqref{cv: P-rate} will follow once we show both $A_{n}\big(\widetilde D_{n}\cap[0, t]\times [0, L]\big)$ and $B_{n}\big(\widetilde D_{n}\cap [0, t]\times [0, L]\big)$ tend to 0 in probability. Starting with $A_{n}(\cdot)$, we have  
\[
\mathbb E\Big[\big|A_{n}\big(\widetilde D_{n}\cap[0, t]\times [0, L]\big)\big|\Big] \le p|\log(1-p)|\sum_{k\le tn^{2/3}}\sum_{l\le Ln^{1/3}} \mathbb E[|m-U_{k}|]\le \cO(p^{2}\cdot n\cdot m^{\frac12}n^{\frac16}) \to 0, 
\]
where we have used the bound $|\log(1-p)|=\cO(p)$ for $p\le \tfrac12$ and the bound for $\mathbb E[(m-U_{k})^{2}]$ from Lemma~\ref{lem: Uk}. On the other hand, writing $a\wedge b=\min(a, b)$, 
we note that
\[
|\log(1-p)|^{-1}B_{n}\big(\widetilde D_{n}\cap[0,t]\times [0, L]\big) = \sum_{k\le tn^{2/3}}(X\n_{k}-U_{k}p)\cdot (S\n_{k-1}-\min_{j\le k-1}S\n_{j})\wedge\lf Ln^{\frac13}\rf. 
\]
Let us denote by $T_{k}$ the $k$-th summand in the previous display. 
Let $k<k'$. Conditional on $U_{k'}$, $X\n_{k'}$ is independent of $X\n_{k}$ and $(S\n_{j})_{j\le k-1}$. It follows that
\[
\mathbb E\big[T_{k}T_{k'}\,|\, U_{k'}] = \mathbb E\big[X\n_{k'}-U_{k'}p|U_{k'}\big]\cdot \mathbb E\Big[T_{k}\cdot (S\m_{k'}-\min_{j'\le k'-1}S\n_{j'})\wedge \lf Ln^{\frac13}\rf\,\Big|\,U_{k'}\Big] =0.
\]
Taking the expectation yields $\mathbb E[T_{k}T_{k'}]=0$ for all $k<k'$. 
We also have 
\[
\mathbb E[T_{k}^{2}] \le L^{2}n^{\frac23}\,\mathbb E\big[(X\n_{k}-U_{k}p)^{2}\big] = L^{2}n^{\frac23}\,\mathbb E\big[U_{k}p(1-p)\big] =  \cO(n^{\frac23}mp). 
\]
Combining this with the previous argument, we deduce that 
\[
\mathbb E\Big[B_{n}\big(\widetilde D_{n}\cap[0,t]\times [0, L]\big)^{2}\Big] = |\log(1-p)|^{2}\sum_{k\le tn^{2/3}}\mathbb E[T_{k}^{2}] =\cO(n^{\frac43}mp^{3})=\cO(n^{\frac13}p)\to 0. 
\]
This shows $B_{n}\big(\widetilde D_{n}\cap [0, t]\times [0, L]\big)\to 0$ in probability, and therefore completes the proof of~\eqref{cv: P-rate}. Regarding~\eqref{cv: P-bd}, it suffices to consider those $C\subseteq [0, t]\times [0, L]$. 
We write $C_{n}=\widetilde D_{n}\cap C$. 
We note that
$q_{k, l}\le X_{n, k}|\log(1-p)|$, since $1-e^{-x}\le x$ for $x\ge 0$. 
Using the fact that conditional on $U_{k}$, $X\n_{k}$ is independent of $\mathbf 1_{\{(k, l)\in D_{n}\}}$, which is determined by $(S\n_{j})_{j\le k-1}$, we deduce that 
\begin{align*}
\sum_{(k, l)\in\N^{2}}\mathbb E\Big[q_{k, l}\mathbf 1_{\{(n^{-\frac23}k, n^{-\frac13}l)\in C_{n}\}}\Big] &\le |\log(1-p)|\sum_{(k, l)\in \N^{2}}\mathbb E\Big[\mathbb E\big[X\n_{k}|U_{k}\big]\,\mathbb P\big( (n^{-\frac23}k, n^{-\frac13}l)\in C_{n}|U_{k}\big)\Big]\\
& \le mp|\log(1-p)|\cdot \mathbb E\big[\#C_{n}\big]. 
\end{align*}
Together with \eqref{cv: volume} and the dominated convergence theorem, this implies \eqref{cv: P-bd}. 
\end{proof}

\begin{proof}[Proof of Proposition~\ref{prop: surplus-bi}]
Let us focus on the moderate clustering regime. The arguments for the light clustering regime are similar. 
Since $\cP^{\lambda,\theta}$ is a simple point measure, according to Proposition 16.17 in \cite{Kal}, to prove the convergence of $(\cP\n)_{n\ge 1}$ to $\cP^{\lambda, \theta}$, it suffices to show the following two conditions are verified. (i) If $B\subset \R^{2}_{+}$ is a finite union of rectangles, then we have 
\begin{equation}
\label{eq: cv-ppp}
\mathbb P\big(P\n(B)=0\big)\xrightarrow{n\to\infty} \mathbb P\big(\cP^{\lambda,\theta}(B) = 0\big) \quad\text{in probability}. 
\end{equation}
(ii) For any compact set $C\subseteq \R^{2}$, we have 
\begin{equation}
\label{eq: cv-ppp'}
\limsup_{n\to\infty}\mathbb E\big[P\n(C)\big]\le \mathbb E\big[\cP^{\lambda,\theta}(C)\big]<\infty.
\end{equation}
Checking~\eqref{eq: cv-ppp}, we let $B=[0, t]\times [0, L]$ for some $t, L>0$. Taking into account the two types of surplus edges, we have
\begin{align*}
 & \mathbb E\Big[\exp\Big(\sum_{k=1}^ {\lf tn^{2/3} \rf}\sum_{l=1}^{\lf Ln^{1/3}\rf} \log(1-q_{k, l})\mathbf 1_{\{(k, l)\in D_{n}\}}\Big)\Big]-\Big\{1-\mathbb P\Big(\bigcap_{k\le tn^{2/3}}\cI\n_{k}\Big)\Big\}\\
&\qquad\qquad \le \mathbb P\big(P\n(B)=0 \big) \\
&\qquad\qquad\qquad\quad \le  \mathbb E\Big[\exp\Big(\sum_{k=1}^ {\lf tn^{2/3} \rf}\sum_{l=1}^{\lf Ln^{1/3}\rf} \log(1-q_{k, l})\mathbf 1_{\{(k, l)\in D_{n}\}}\Big)\Big]+1-\mathbb P\Big(\bigcap_{k\le tn^{2/3}}\cI\n_{k}\Big).
\end{align*}
Lemma~\ref{lem: goodevent}, Lemma \ref{lem: cv-pp}, and the dominated convergence theorem all combined yield the convergence in \eqref{eq: cv-ppp}. 
This then extends to a general rectangle of the form $[t_{1}, t_{2}]\times [L_{1}, L_{2}]$ with the inclusion-exclusion principle. 
The case of a finite union of rectangles can be similarly argued. Regarding \eqref{eq: cv-ppp'}, we note  that $\mathbb E[\cP^{\lambda, \theta}(C)]=\mathbb E[\mathrm{Leb}(\cD^{\lambda,\theta}\cap C)]$ and $P\n$ is a simple point measure. In particular, if there is a surplus edge of type (ii) as described in Lemma~\ref{lem: surplus} associated to $k\in [n]$, then that yields an atom at $(n^{-2/3}k, 0)$ in $P\n$.  It follows that 
\[
\mathbb E\Big[P\n(C)\Big] \le \sum_{(k, l)\in\N^{2}}\mathbb E\Big[q_{k, l}\mathbf 1_{\{(n^{-\frac23}k, n^{-\frac13}l)\in \widetilde D_{n}\cap C\}}\Big]+\sum_{k\in \N}\Big(1-\mathbb P(\cI\n_{k})\Big)\mathbf 1_{\{(n^{-\frac23}k, 0)\in C\}}. 
\]
The desired bound then follows from~\eqref{cv: P-bd} and \eqref{bd: goodevent}. 
\end{proof}

\subsection{Gromov--Hausdorff--Prokhorov convergence of the graphs}
\label{sec: cv-ghp}

Having obtained the uniform convergence of the graph exploration processes, i.e.~Propositions~\ref{thm: cv-moderate} and \ref{thm: cv-light}, as well as the convergence of the point processes generating bipartite surplus edges in Proposition~\ref{prop: surplus-bi}, we explain in this part how to obtain Theorems~\ref{thm2} and~\ref{thm1}. 
We will only detail the arguments in the moderate clustering regime; the case of light clustering regime can be argued similarly. 
Firstly, we provide a proof of Proposition~\ref{prop: approx-gr}.  Recall that $(\hat C\n_{k},  \dgr^{\bi}, \hat\mu\n_{k})$ is the $k$-th largest connected component of $B(n, m, p)$ in $\cV$-vertex numbers, equipped with its graph distance and the counting measure on $\cV$. Proposition~\ref{prop: approx-gr}, once established, will allow us to replace Theorem~\ref{thm2} with the following equivalent convergence under \eqref{hyp: M}: 
\begin{equation}
\label{eq: cv-thm2'}
\Big\{\Big(\hat C^{n, m}_{k}, \, \tfrac12 n^{-\frac13}\cdot \dgr^{\bi}, \, n^{-\frac23}\hat\mu\n_{k}\Big): k\ge 1\Big\} \;\Longrightarrow\; \mathcal G^{\mathrm{RIG}}(\lambda, \theta)
\end{equation}
in the same sense of convergence as \eqref{eq: cv-thm2}. 

\begin{proof}[Proof of Proposition~\ref{prop: approx-gr}]
Take any pair $(i, j)\in [n]^{2}$. Let us show that for each $k\in \N\cup\{0\}$, 
\begin{equation}
\label{id: isometry}
\dgr^{\RIG}(i, j) = k  \quad \Leftrightarrow \quad \dgr^{\bi}(w_{i}, w_{j}) = 2k.
\end{equation}
This is clear for $k=0$. Suppose that $\dgr^{\RIG}(i, j)=k\in\N$. Then there is a sequence of vertices $i=i_{1},  i_{2}, i_{3}, \dots, i_{k}=j$ where $i_{l}$ is adjacent to $i_{j+1}$ in $G(n, m, p)$, $1\le l\le k-1$. By definition, this means that we can find a sequence of  
$\cU$-vertices $u_{i_{1}}, u_{i_{2}}, \dots, u_{i_{k-1}}$ so that $w_{i_{l}}$ and $w_{i_{l+1}}$ are both neighbours of $u_{i_{l}}$ in $B(n, m, p)$, $1\le l \le k-1$. It follows that $\dgr^{\bi}(w_{i}, w_{j})\le 2k$. Conversely, if $\dgr^{\bi}(w_{i}, w_{j})=2k'\ge 2$, then there is a sequence of vertices alternating between $\cV$- and $\cU$-vertices: $w_{i}=w_{i_{1}}, u_{i_{1}}, w_{i_{2}}, u_{i_{2}}, \dots, u_{i_{k'-1}}, w_{i_{k'}}=w_{j}$ where $w_{i_{l}}$ and $w_{i_{l+1}}$ are adjacent to $u_{i_{l}}$, for $1\le l\le k'-1$. This implies that $\dgr^{\RIG}(i, j)\le k'$.  Combining the two inequalities yields the desired identity \eqref{id: isometry}. 

We note that \eqref{id: isometry} implies that two vertices $i, j$ are connected in $G(n, m, p)$ if and only if $w_{i}, w_{j}$ are found in the same connected component of $B(n, m, p)$. Writing $d_{\mathrm{H}}$ for the Hausdorff distance between subsets of $B(n, m, p)$, we have $d_{\mathrm{H}}(\rho_{n}(C\n_{k}), \hat C\n_{k})\le 1$, since every $\cU$-vertex in $\hat C\n_{k}$ is adjacent to some $\cV$-vertex. On the other hand, $\hat\mu\n_{k}$ is the push-forward of $\mu\n_{k}$ by $\rho_{n}$. The bound on $\dghp$ follows. 
\end{proof}

Recall~\eqref{cv: moderate} concerning the convergence of $(\tilde R^{n,m}, S^{n,m}, H^{n,m})$. Let us first point out that it also implies the convergence of the so-called contour processes for the spanning forest $\cF$. Informally, denote by $\cF'$ the sub-forest of $\cF$ on $\cV$ where $w_i$ is a parent of $w_j$, $1\le i, j\le n$ if and only if the former is a grandparent of the latter in $\cF$. Consider each tree component of $\cF'$ embedded into the plane, with each edge having length one and siblings are ordered with the oldest to the left and youngest to the right.  Imagine a particle visiting each tree in $\cF'$ in turn, each time starting at the root, exploring edges from left to right, backtracking as little as possible and travelling at unit speed; record at time $t$ the distance of the particle to the root of the tree currently visited and call this quantity $D^{n,m}_t$. 
By following the same arguments as in Section 2.4 of~\cite{DuLG02}, we can improve the convergence in~\eqref{cv: moderate} to the following in $\mathbb D^{\ast}(\R_+, \R^4)$:
\begin{equation}
\label{cv: moderate'}
\Big\{n^{-\frac13}\tilde R\n_{\lfloor n^{2/3}t\rfloor}, n^{-\frac13}S\n_{\lf n^{2/3}t\rf}, n^{-\frac13}H\n_{\lfloor n^{2/3}t\rfloor}, n^{-\frac13}D\n_{n^{2/3}t}: t\ge 0\Big\} \Longrightarrow \Big\{  \cR^{\lambda, \theta}_{t}, \cS^{\lambda, \theta}_{t}, \cH^{\lambda, \theta}_{t}, \cH^{\lambda, \theta}_{t/2}: t\ge 0\Big\}. 
\end{equation}
In particular, part of the arguments as in Section 2.4 of~\cite{DuLG02} show that there exists a map $f: \R_+\to [n]$ satisfying $f(t) = k$ if $v_k$ is the vertex visited by the particle at time $t$ and $\sup_{t\le t_0}|f(t)-t/2|\to 0$ in probability for all $t_0\ge 0$.  

From now on, we assume that the previous convergence~\eqref{cv: moderate'} holds {\it almost surely}, which is legitimate thanks to Skorokhod's Representation Theorem. 
The excursions of $\cS^{\lambda,\theta}$ above its running infimum can be ranked in a decreasing order of their sizes, as shown by the arguments in Section~\ref{sec: RIGgraph}. Recall the notation 
$\zeta^{\lambda, \theta}_{k}$ for the length of the $k$-th longest such excursion. We denote by $|\mu\n_{k}|=|\hat\mu\n_{k}|$ the number of $\cV$-vertices in $\hat C\n_{k}$.  Let us point out that Aldous' theory on size-biased point processes (see Proposition 15 and Lemma 7 in~\cite{Al97}), applied to Proposition~\ref{thm: cv-moderate}, entails the following:
\begin{equation}
\label{eq: cv-l2}
\Big\{n^{-\frac23}|\mu\n_{k}|: k\ge 1\Big\}\quad \Longrightarrow\quad \big\{\zeta^{\lambda, \theta}_{k}: k\ge 1\big\},
\end{equation}
with respect to the weak topology of $\ell^{2}$. 

We will also require the convergence of the excursion intervals. 
Recall that $(\mathrm g^{\lambda, \theta}_{k}, \mathrm d^{\lambda, \theta}_{k})$ is the $k$-th longest excursion of $\cS^{\lambda,\theta}$ above its running infimum. 
For $n\ge 1$, let $(2\mathrm g\n_{k}, 2\mathrm d\n_{k})$ be the $k$-th longest excursion of $D\n$ above $0$, which is unique for $n$ sufficiently large. Recall that $(\mathrm g^{\lambda, \theta}_{k}, \mathrm d^{\lambda, \theta}_{k})$ is also the $k$-th longest excursion interval of $\cH^{\lambda, \theta}$ above $0$, with the excursion of $\cH^{\lambda,\theta}$ running on the interval denoted as $\mathrm e^{\lambda, \theta}_{k}$. We introduce its discrete counterpart by setting 
\[
e\n_{k}(s)= D\n_{2\mathrm g\n_{k}+2s }, \quad 0\le s\le \zeta\n_{k}:=\mathrm d\n_{k}-\mathrm g\n_{k}. 
\]
By following almost verbatim the proofs of Lemmas 5.2-5.4 in \cite{BrDuWa21}, we can show that  
\begin{equation}
\label{cv: exc-int}
\Big\{\Big(n^{-\frac23}\mathrm g\n_{k}, n^{-\frac23}\mathrm d\n_{k}\Big): k\ge 1\Big\} \quad\Longrightarrow\quad \big\{\big(\mathrm g^{\lambda, \theta}_{k}, \mathrm d^{\lambda, \theta}_{k}\big): k\ge 1\big\}  
\end{equation}
with respect to the product topology of $(\R^{2})^{\N}$; we omit the details. As $\cH^{\lambda,\theta}$ has continuous sample paths, \eqref{cv: exc-int} then implies that for each $t>0$, 
\begin{equation}
\label{cv: h-exc}
\Big\{\Big(n^{-\frac13}e\n_{k}(n^{\frac23}s\wedge t)\Big)_{0\le s\le n^{-\frac23}\zeta\n_{k}}: k\ge 1\Big\} \quad\Longrightarrow\quad \Big\{\big(e^{\lambda, \theta}_{k}(s\wedge t)\big)_{s\le \zeta^{\lambda, \theta}_{k}}: k\ge 1\Big\}
\end{equation}
with respect to the product topology of $(\mathbb D^{\ast}([0, t],\R))^{\N}$. Recall from Section~\ref{sec: sps} the point measures $P\n$ and $\cP^{\lambda, \theta}$, and from Section~\ref{sec: RIGgraph} the finite collection $\Pi^{\lambda, \theta}_{k}$. Replacing $\mathrm e^{\lambda, \theta}_{k},\cP^{\lambda,\theta}$ with $\mathrm e\n_{k}, P\n$ in the construction there yields its discrete counterpart $\Pi\n_{k}$. More precisely, let $(\tilde x_{k, i}, \tilde y_{k, i}),  1\le i\le p\n_{k}$, be the elements of $P\n\cap \{(n^{-2/3}\mathrm g\n_{k}, n^{-2/3}\mathrm d\n_{k})\times \R\}$. For $1\le i\le p\n_{k}$, set 
\[
\check s_{k, i}= \max\Big\{u\le n^{\frac23}\tilde x_{k, i}-1: S\n_{u}\le S\n_{n^{\frac23}\tilde x_{k, i}-1}-n^{\frac13}\tilde y_{k, i}\Big\}+1,\quad  \check t_{k, i}=n^{\frac23}\tilde x_{k, i}. 
\]
{Let us explain the intuition behind the previous definition. Recall from Lemma~\ref{lem: surplus} that there are two types of surplus edges. If $(\tilde x_{k, i}, \tilde y_{k, i})\in P\n$ corresponds to a surplus edge $e=\{u, w\}$ of type (i), then 
we have $n^{2/3}\tilde x_{k, i}=k(e)$ and $n^{1/3}\tilde y_{k, i}=l(e)$, where 
$v_{k(e)}$ is adjacent to $u$ and $l(e)$ is the rank of $w$ in $\cA_{k-1}$. In particular, there is some $j(e)\le k(e)-1$ satisfying $w\in \cN_{j(e)}$.  Since the vertices discovered after $w$ are added ahead of $w$ in the active list, we have $\#\cA_i\ge \#\cA_{k-1}-l(e)+1$ for each $i\in[j(e), k(e)-1]$, and $\#\cA_{j(e)-1}\le \#\cA_{k-1}-l(e)$. 
It follows that $j(e)-1=\max\{i\le k(e)-1: \#\cA_i\le \#\cA_{k(e)-1}-l(e)\}$. 
Combining this with the identity from Lemma~\ref{lem: cA} that $S\n_{i}-\min_{u\le i}S\n_{u}=\#\cA_{i}$, 
we find that $\check s_{k, i}=j(e)$. To recap, for a surplus edge of type (i), the corresponding $(\check s_{k, i}, \check t_{k, i})$ identifies a pair of black vertices, $v_{k(e)}$ and $v_{j(e)}$ more precisely, that are no more than $2$ edges away from the endpoints of $e$. This also works for type (ii) surplus edges. Indeed, let $(\tilde x_{k, i}, \tilde y_{k, i})\in P\n$ be associated to such an edge, say $e=\{u, w\}$; then in this case we have $v_{k(e)}$ adjacent to $u$ and $w$ at $2$ edges away from $v_{k(e)}$. On the other hand, we must have $\tilde y_{k, i}=l(e)=0$ in this case, and thus $\check s_{k, i}=\check t_{k, i}$, so that both $\check s_{k, i}$ and $\check t_{k, i}$ point to $v_{k(e)}$ in this case. 
Denote by $\rho: k\in [n]\mapsto v_k$ and by $\mathrm d$ the graph distance of the bipartite forest $\cF$. Then the previous arguments show that for each surplus edge $e=\{u, w\}$, for the corresponding atom $(\tilde x_{k, i}, \tilde y_{k, i})\in P\n$, we have $\mathrm d(\rho(\check t_{k, i}), u)\le 2$ and $\mathrm d(\rho(\check s_{k, i}), w)\le 2$. }

Recall the map $f$ from below~\eqref{cv: moderate'} and let $f^{-1}$ be its right-continuous inverse. Let $n^{2/3}\tilde s_{k, i}=f^{-1}(\check s_{k, i})-2\mathrm g\n_k$ and $n^{2/3}\tilde t_{k, i}=f^{-1}(\check t_{k, i})-2\mathrm g\n_k$. 
Define $\Pi\n_{k}=\{(\tilde s_{k, i}, \tilde t_{k, i}): 1\le  i\le p\n_{k}\}$. Recalling from \eqref{def: graph} the measured metric space $\cG(h, \Pi)$, we set $(C^{\ast}_{n, k}, d^{\ast}_{n, k}, \mu^{\ast}_{n, k})  =\cG(2\mathrm e\n_{k}, \Pi\n_{k})$. We claim that 
\begin{equation}
\label{bd: ghp}
 \dghp\big((C^{\ast}_{n, k},  d^{\ast}_{n, k}, \mu^{\ast}_{n, k}), (\hat C\n_{k}, \dgr^{\bi}, \hat\mu\n_{k})\big) \le 2(p\n_{k}+1).
\end{equation}
Putting aside the verification of \eqref{bd: ghp} for the moment, let us explain how it completes the proof of Theorem~\ref{thm2}. 

\begin{proof}[Proof of Theorem~\ref{thm2}]
Identifying $\Pi\n_{k}$ as a point measure of $\R^{2}$ with a unit mass located at each of its member, we deduce from Proposition~\ref{prop: surplus-bi}, \eqref{cv: exc-int} and \eqref{cv: h-exc} that 
\[
\big\{\Pi\n_{k}: k\ge 1\big\} \quad \Longrightarrow\quad \big\{\Pi^{\lambda,\theta}_{k}: k\ge 1\big\}, 
\]
with respect to the product topology of weak topology for finite measures. 
This, along with  \eqref{cv: h-exc}, entails the convergence of $\{(C^{\ast}_{n, k}, \tfrac12 n^{-\frac13}d^{\ast}_{n, k}, n^{-\frac23}\mu^{\ast}_{n, k}): k\ge 1\}$ to $\cG^{\RIG}(\lambda, \theta)$, thanks to \eqref{eq: ghp}. The convergence in \eqref{eq: cv-thm2'} then follows 
from \eqref{bd: ghp}. 
\end{proof}

All it remains is to prove~\eqref{bd: ghp}. To that end, let $(T^{\ast}_{k}, d^{\ast}_{k}, m^{\ast}_{k})$ be the measured real tree encoded by $2 \mathrm{e}\n_{k}$ in the sense of \eqref{def: real-tree}. 
Recall the bipartite forest $\cF$ output by Algorithm~\ref{alg}. 
Let $\mathrm T_{k}$ be the $k$-th largest tree component of $\cF$ in $\cV$-vertex numbers. Then $\mathrm T_{k}$ is a spanning tree of $C\n_{k}$, at least when $n$ is sufficiently large.  
We turn it into a measured metric space by equipping it with the graph distance $\mathrm d_{k}$ of $\cF$ restricted to $\mathrm T_{k}$, and with the counting measure $\mathrm m_{k}$ of $\cV$. 
Since $2 H\n_{k}$ corresponds to the height of $v_{k}$ in $\mathrm T_{k}$, we can find an isometric embedding of $(T^{\ast}_{k}, d^{\ast}_{k})$ into $(\mathrm T_{k}, \mathrm d_{k})$. As in the proof of  Proposition~\ref{prop: approx-gr}, this implies that 
\[
\dghp((T^{\ast}_{k}, d^{\ast}_{k}, m^{\ast}_{k}), (\mathrm T_{k}, \mathrm d_{k}, \mathrm m_{k})\big)\le 1. 
\]
Next, we note that $(C^{\ast}_{k}, d^{\ast}_{k})$ is obtained from $(T^{\ast}_{k}, d^{\ast}_{k})$ by introducing shortcuts that are encoded by $\Pi\n_{k}$. 
The definitions of $\tilde s_{k, i}$ and $\tilde t_{k, i}$ imply that for each surplus edge $\{u, w\}$, we can find some $(\tilde s_{k, i}, \tilde t_{k, i})\in \Pi\n_{k}$ satisfying
$\mathrm d_{k}\big(\rho(n^{2/3}\tilde t_{k, i}+\mathrm g\n_k), u\big)\le 2$ and $\mathrm d_{k}\big(\rho(n^{2/3}\tilde s_{k, i}+\mathrm g\n_k), w\big)\le 2$. 
Note also that $p\n_{k}=\#\Pi\n_{k}$ is bounded by the number of surplus edges in $\hat C\n_{k}$. 
Since a geodesic in $C^{\ast}_{k}$ contains at most $p\n_{k}$ shortcuts, the bound in \eqref{bd: ghp} follows. We refer to Lemma 21 in \cite{ABBrGo12} and Appendix C of \cite{BrDuWa21} for similar arguments. 

\begin{proof}[Proof of Theorem~\ref{thm1}]
This is similar to the proof of Theorem~\ref{thm2}. We therefore omit the details. 
\end{proof}

\subsection{Triangle counts in the random intersection graphs}
\label{sec: triangle}

This subsection contains the proof of Theorems~\ref{thm: triangle-moderate} and \ref{thm: triangle-light}, which pertains to triangle counts in the moderate and light clustering regimes. 
We begin with the following observation: a triangle in $G(n, m, p)$  can be classified into two types. Either the  three vertices of the triangle belong to a common group--we call such a triangle of {\it type I}--or there is no such common group; in that case, we say the triangle is of {\it type II}.  As we will see, in most of the cases that we are interested in,  the numbers of type I triangles will dominate those of type II. So the first step of the proof is to establish limit theorems for type I triangles. 
To that end, let us denote for $N\in \Z_{+}$ and $k\in\N$,
\[
(N)_{k} = N(N-1)(N-2)\cdots (N-k+1)
\]
for the falling factorial. 
We observe that inside a complete graph with $N$ vertices, there are precisely $\frac16(N)_{3}$ triangles. This motivates the introduction of the following discrete-time process $(T\n_{k})_{k\ge 0}$: 
let $T\n_{0}=0$, and more generally for $k\ge 1$:
\begin{equation}
\label{def: deltaT}
\Delta T\n_{k}:= T\n_{k}-T\n_{k-1}:= \sum_{i=1}^{X\n_{k}} \tfrac16\big(1+\#\cN_{k, i}\big)_{3}\,,  
\end{equation}
where we recall $X\n_{k}$ and $\cN_{k, i}$ from Algorithm~\ref{alg}. We note that $1+\#\cN_{k, i}$ corresponds to the degree of the vertex $u_{k, i}$ in the bipartite spanning forest $\cF$. 
As $v_{k}$ is always a neighbour of every $u\in \cM_{k}$, this accounts for the term 1 there. 
The scaling limit of the process $(T\n_{k})_{k\ge 0}$, stated in the two following propositions, is the key ingredient in determining the asymptotic numbers of type I triangles.

\begin{prop}
\label{prop: moderate-tri-ct}
Under the assumption~\eqref{hyp: M}, the following convergence takes place in probability: 
\[
\Big\{n^{-\frac23}T\n_{\lf tn^{2/3}\rf}: t\ge 0\Big\} \xrightarrow{n\to\infty} \{c_{\theta}\cdot t: t\ge 0\} \quad\text{in }\mathbb D^{\ast}(\R_{+}, \R), 
\]
where $c_{\theta}$ is the constant defined in \eqref{def: ctheta}. 
\end{prop}

\begin{prop}
\label{prop: light-tri-ct}
Under the assumption~\eqref{hyp: L}, the following statements hold true.
\begin{enumerate}[(i)]
\item
If $\frac{m}{n^{7/3}}\to 0$ as $n\to\infty$, then the following convergence takes place in probability: 
\begin{equation}
\label{eq: cv-T1}
\Big\{m^{\frac12}n^{-\frac76}T\n_{\lf tn^{2/3}\rf}: t\ge 0\Big\} \xrightarrow{n\to\infty} \{\tfrac12 t: t\ge 0\} \quad\text{in }\mathbb D^{\ast}(\R_{+}, \R). 
\end{equation}
\item
If $\frac{m}{n^{7/3}}\to\infty$ as $n\to\infty$, then for any $t\ge 0$, we have 
\begin{equation}
\label{eq: cv-T2}
\max_{s\le tn^{2/3}}T\n_{s} \to 0 \quad \text{in probability}.
\end{equation}
\item
If  $m=\lf n^{7/3}\rf$, then jointly with the convergence in \eqref{cv: light}, the following convergence in distribution takes place with respect to the Skorokhod topology for c\`adl\`ag paths on $\R_{+}$: 
\begin{equation}
\label{eq: cv-T3}
\Big\{T\n_{\lf tn^{2/3}\rf}: t\ge 0\Big\} \ \Longrightarrow\  \{\mathscr T_{t}: t\ge 0\}, 
\end{equation}
where $\{\mathscr T_{t}: t\ge 0\}$ is a Poisson process of rate $\frac12$ on $\R_{+}$, independent of the limit processes in \eqref{cv: light}. 
\end{enumerate}
\end{prop}

\begin{proof}[Proof of Proposition~\ref{prop: moderate-tri-ct}]
We use the method of moments: the conclusion will follow once we show that for all $t>0$, 
\begin{equation}
\label{eq: cvtrE}
 \sup_{s\le t}\Big|n^{-\frac23}\mathbb E\Big[T\n_{\lf sn^{2/3}\rf}\Big]- c_{\theta}s \Big|\to 0 , \ 
 n^{-\frac43}\max_{k\le tn^{2/3}}\Var\big(T\n_k\big) \to 0.
\end{equation}
We recall $V_{k, i}=V_{k}-\sum_{j<i}\#\cN_{k, j}$. Conditional on $V_{k, i}$, $Y_{k, i}:=\#\cN_{k, i}$ has the Binom$(V_{k, i}, p)$ distribution, independent of $Y_{k, j}, j<i$. 
According to~\eqref{id: binom-moment}, we have
\[
\mathbb E\big[(1+Y_{k, i})_{3}\big] = 3\mathbb E[(Y_{k, i})_{2}]+\mathbb E[(Y_{k, i})_{3}]=3p^{2}\mathbb E[(V_{k, i})_{2}]+p^{3}\mathbb E[(V_{k, i})_{3}].
\]
We note that $V_{k+1}\le V_{k, i}\le n$. Combined with~\eqref{hyp: M}, this implies that 
\[
3p^{2}\mathbb E\big[\big(V_{k+1}\big)_{2}\big] + p^{3}\mathbb E\big[\big(V_{k+1})_{3}\big]\le \mathbb E\big[(1+Y_{k, i})_{3}\big]\le 3n^2p^2 + n^3p^3 = \mathbb E[(Y_{\theta}+1)_3]+\cO(n^{-\frac13}), 
\]
where $Y_{\theta}$ is a Poisson$(1/\sqrt\theta)$ variable. 
To deal with the lower bound above, we note that Lemma~\ref{lem: Uk} tells us that $\mathbb E[n-V_{\lf tn^{2/3}\rf+1}]=\cO(n^{2/3})$, from which it follows that
\begin{align*}
\max_{k\le tn^{2/3}}\Big(n^{3}p^{3}-p^{3}\mathbb E\big[(V_{k+1})_{3}\big]\Big)  &\le n^{3} p^{3}-p^{3}\mathbb E\Big[\big(V_{\lf tn^{2/3}\rf+1}-2\big)_{+}^{3}\Big]\\
& \le4n^{2}p^{3}\mathbb E\Big[n+2-V_{\lf tn^{2/3}\rf+1}\Big]  =\cO(n^{-\frac13}),  
\end{align*}
where we have used the elementary inequality $0\le a^{3}-b^{3}\le (a-b)(a+b)^{2}\le (2a)^{2}(a-b)$ for all $0\le b\le a<\infty$. Similarly, we have $\mathbb E[(V_{k})_{2}]p^{2}= \theta^{-1}+\cO(n^{-\frac13})$, where the error term is uniform for all $k\le tn^{2/3}$. We have shown that $\mathbb E[(Y_{k,i}+1)_{3}]= \mathbb E[(Y_{\theta}+1)_{3}]+\cO(n^{-\frac13})$ uniformly for all $i\ge 1$ and $k\le tn^{2/3}$. 
Next, we note that 
\[
\mathbb E\big[\Delta T\n_k\,\big|\, X\n_k\big] =\sum_{i=1}^{X\n_k}\tfrac16 \mathbb E\big[(1+Y_{k, i})_{3}\big]. 
\]
Combining this with the fact that $X\n_{k}$ has the Binom$(U_{k}, p)$ distribution, the assumption~\eqref{hyp: M}, and the bound $\max_{k\le tn^{2/3}}\mathbb E[(m-U_{k})p]=\cO(n^{-\frac13})$ from Lemma~\ref{lem: Uk}, we conclude that 
\begin{equation}
\label{cv: Delta-T}
\max_{k\le tn^{2/3}}\big| \mathbb E[\Delta T\n_{k}] - c_{\theta}\big|\to 0, \ \text{ hence } \  \sup_{s\le t}\Big|n^{-\frac23}\mathbb E\Big[T\n_{\lf sn^{2/3}\rf}\Big]- c_{\theta}s \Big|\to 0,
\end{equation}
where we recall $c_{\theta}$ from~\eqref{def: ctheta}. 
Let us further note that the upper bounds above hold in a stronger sense: indeed, we have 
\begin{equation}
\label{bd: DeltaA}
\mathbb E[(Y_{k, i}+1)_3\,|\, U_k, V_k]\le \mathbb E[(Y_{\theta}+1)_3]+o(1)\  \text{and then}\  \Delta A_k:=\mathbb E[\Delta T\n_k\,|\, U_k, V_k]\le c_{\theta}+o(1).
\end{equation}
We shall need these bounds later on.  
We now turn to the second moments. First, we note that for $1\le j<i$, by the conditional independence between $Y_{k, i}$ and $Y_{k, j}$, we have
\begin{align*}
\mathbb E\big[(1+Y_{k, i})_3(1+Y_{k, j})_3\,\big|\, V_{k, i}, U_k, V_k\big]
&=\mathbb E\big[ (1+Y_{k, j})_3\,\big|\,  V_{k, i}, U_k, V_k]\cdot (3p^2 (V_{k, i})_2+p^3(V_{k, i})_3\big)\\
& \le \mathbb E\big[ (1+Y_{k, j})_3\,\big|\,  V_{k, i}, U_k, V_k]\cdot(3(n)_2p^2+(n)_3p^3)
\end{align*}
Taking expectation and using the fact that conditional on $V_k$, $(1+Y_{k, i})_{3}$ is still stochastically bounded by $(1+Y)_{3}$, where $Y$ is a Binom$(n, p)$ variable, we find that 
\[
\mathbb E\big[(1+Y_{k, i})_3(1+Y_{k, j})_3\,\big|\,U_k, V_k\big]\le (3(n)_2p^2+(n)_3p^3)^2=\big(\mathbb E[(1+Y)_3]\big)^2.
\]
Thus, for each $\ell\in \N$, 
\begin{align*}
&\qquad 36\mathbb E\big[(\Delta T\n_k)^2\mathbf 1_{\{X\n_k=\ell\}}\,\big|\,U_k, V_k\big] \\&
=\sum_{1\le i\le \ell}\mathbb E\Big[\big((1+Y_{k, i})_3\big)^2\,\Big|\,U_k, V_k\Big] + 2\sum_{1\le j<i\le \ell} \mathbb E\Big[(1+Y_{k, i})_3(1+Y_{k, j})_3\,\Big|\,U_k, V_k\Big]\\
&\le \ell \cdot \mathbb E\Big[\big((1+ Y)_{3}\big)^2\Big] + \ell(\ell-1)\cdot \Big(\mathbb E \big[(1+Y)_{3}\big]\Big)^2 
\le \ell^2\, \mathbb E\big[(1+Y)^{6}\big]. 
\end{align*}
As $(np)_{n\ge 1}$ remains bounded under~\eqref{hyp: M}, computations on the binomial distributions (see~\eqref{id: binom-moment} and~\eqref{id: poisson-moment}) show that 
$\mathbb E\big[(1+Y)^{6}\big]=\cO(1)$. On the other hand, $X\n_k$, when conditioned on $U_k, V_k$, is still dominated by Binom$(m, p)$, so that $\mathbb E[(X\n_k)^2\,|\,U_k, V_k]\le (mp)^2+mp=\cO(1)$ under~\eqref{hyp: M}. Summing over $\ell$ in the previous display yields
\[
 \max_{k\le tn^{2/3}} \mathbb E\big[(\Delta T\n_k)^2\,\big|\,U_k, V_k\big]\le \max_{k\le tn^{2/3}} \mathbb E\big[(X\n_k)^2\,\big|\,U_k, V_k\big] \cdot \mathbb E\big[(1+Y)^{6}\big]=\cO(1).
\]
{Recall that $\Delta A_j=\mathbb E[\Delta T\n_j\,|\,U_j, V_j]$. Note that we also have $\mathbb E[\Delta T\n_j\,|\, T\n_{j-1}, U_j, V_j]=\Delta A_j$ and $\mathbb E[(\Delta T\n_j)^2\,|\, T\n_{j-1}, U_j, V_j]=\mathbb E[(\Delta T\n_j)^2\,|\, U_j, V_j]$. It follows that
\begin{align*}
\mathbb E\big[(T\n_k)^2\big] 
& = \mathbb E\Big[\sum_{j\le k}\mathbb E\big[(T\n_j)^2-(T\n_{j-1})^2\,\big|\,T\n_{j-1}, U_j, V_j\big]\Big] \\
&=\sum_{j\le k} 2\,\mathbb E[ \Delta A_j\cdot T\n_{j-1}] +\mathbb E\Big[\sum_{j\le k} \mathbb E[(\Delta T\n_j)^2\,|\, U_j, V_j]\Big]\\
& \le 2\big(c_{\theta}+o(1)\big)\sum_{j\le k}\mathbb E[T\n_{j-1}] + k\cdot \cO(1)\\
&\le 2\big(c_{\theta}+o(1)\big)^2\sum_{j\le k}j + k\cdot \cO(1)\le \big(c_{\theta}+o(1)\big)^2k^2 +k\cdot \cO(1), 
\end{align*}
where we have used in the third line~\eqref{bd: DeltaA} and the previous bound for $\mathbb E[(\Delta T\n_k)^2\,|\,U_k, V_k]$, and in the last line the bound $\mathbb E[T_k]\le k(c_{\theta}+o(1))$ from~\eqref{cv: Delta-T}. 
We deduce again from~\eqref{cv: Delta-T} that $\mathbb E[T_k]\ge k(c_{\theta}+o(1))$. 
It then follows that 
\[
\max_{k\le tn^{2/3}}\Var(T\n_k) = \max_{k\le tn^{2/3}}\Big(\mathbb E[(T\n_k)^2]-(\mathbb E[T\n_k])^2\Big)\le o(n^{\frac43}),
\]
so that \eqref{eq: cvtrE} holds.}
\end{proof}

\begin{proof}[Proof of Proposition~\ref{prop: light-tri-ct}: cases (i) and (ii)]
As in the previous proof, we employ the method of moments. As previously, let us denote $Y_{k, i}=\#\cN_{k, i}$ which has distribution Binom$(V_{k, i}, p)$ with $V_{k, i}=V_k-\sum_{j<i}\#\cN_{k, j}$. We then have
\[
3p^{2}\mathbb E\big[\big(V_{k+1}\big)_{2}\big] + p^{3}\mathbb E\big[\big(V_{k+1})_{3}\big]\le \mathbb E\big[(1+Y_{k, i})_{3}\big]\le 3n^2p^2 + n^3p^3
\]
Since $np\to 0$ under~\eqref{hyp: L}, together with the bound $\mathbb E[n-V_{k+1}]=\cO(n^{2/3})$ from Lemma~\ref{lem: Uk}, we deduce that
\[
\mathbb E\big[(1+Y_{k, i})_{3}\big] = 3(np)^2(1+o(1))\quad\text{and}\quad {\mathbb E\big[(1+Y_{k, i})_{3}\,\big|\,U_k, V_k\big] \le 3(np)^2(1+o(1))},
\]
where the error term is uniform for all $i\ge 1, k\le tn^{2/3}$.  Combined with~\eqref{bd: uk}, which implies that $\mathbb E[X\n_k] = \mathbb E[U_k]p\sim mp$ under~\eqref{hyp: L},  this entails that 
\[
\max_{k\le tn^{2/3}}\Big|\mathbb E[\Delta T\n_{k}] -\frac12 \sqrt{\frac{n}{m}}\Big|\to 0, \quad \text{since}\quad \mathbb E[\Delta T\n_k]=\mathbb E\Big[\sum_{1\le i\le X\n_k}\tfrac16(1+Y_{k, i})_3\Big].
\]
If $n^{7/3}=o(m)$, then $n^{\frac23}\sqrt{n/m}\to 0$. In that case,  
$\mathbb E[T\n_{\lf tn^{2/3}\rf}]\to 0$, from which \eqref{eq: cv-T2} follows, as $(T\n_{k})_{k\ge 0}$ is non decreasing and positive. 
If, instead,  $m=o(n^{7/3})$, then we deduce that 
\[ 
\sup_{s\le t}\Big|m^{\frac12}n^{-\frac76}\mathbb E\Big[T\n_{\lf sn^{2/3}\rf}\Big]- \frac{s}{2} \Big|\to 0. 
\]
Note that in this case we also have $\max_{k\le tn^{2/3}}\mathbb E[\Delta T\n_k\,|\,U_k, V_k]\le \frac12\sqrt{n/m}(1+o(1))$. 
We next need a bound on the variance. To that end, we note from~\eqref{bd: binom2} that 
\begin{equation}
\label{bd: hytt}
\mathbb E\big[\big((1+Y_{k, i})_3\big)^2\,\big|\,U_k, V_k\big] \le 18 n^2p^2 (1+o(1)),
\end{equation}
where the error term above uniform for $i$ and $k$.  
For $1\le j<i$, using the conditional independence between $Y_{k, i}$ and $Y_{k, j}$, we can show that 
\[
\mathbb E\big[(1+Y_{k, i})_3(1+Y_{k, j})_3\,\big|\,U_k, V_k\big] \le  \big(3n^2p^2(1+o(1))\big)^2 = 9n^4p^4(1+o(1)),
\]
where $o(1)$ is uniform for $i, j\ge 1$ and $k\le tn^{2/3}$. 
Together with~\eqref{bd: hytt}, this implies that
\begin{align*}
&\qquad \mathbb E\big[(\Delta T\n_k)^2\,\big|\,U_k, V_k\big]\\
&=\frac{1}{36}\cdot\mathbb E\Big[\sum_{1\le i\le X\n_k} \big((1+Y_{k, i})_3\big)^2\,\Big|\,U_k, V_k\Big] + \frac{1}{18}\cdot\mathbb E\Big[\sum_{1\le j<i\le X\n_k} (1+Y_{k, i})_3(1+Y_{k, j})_3\,\Big|\, U_k, V_k\Big]\\
&\le\frac12 \mathbb E\big[X\n_k\,\big|\,U_k, V_k\big] \cdot n^2p^2\big(1+o(1)\big) + \frac14 \mathbb E\big[(X\n_k)_2\,\big|\,U_k, V_k\big]\cdot n^4p^4\big(1+o(1)\big)\\
& = \frac12 \sqrt{\frac{n}{m}}(1+o(1))+\cO\Big(\frac{n}{m}\Big) = \frac12 \sqrt{\frac{n}{m}}(1+o(1)), 
\end{align*}
where we have used ${\mathbb E[X\n_k\,|\,U_k, V_k] \le} mp$ and ${\mathbb E[(X\n_k)_2\,|\,U_k, V_k]} \le m^2p^2$, and $n/m\to 0$ under~\eqref{hyp: L}. The error term above is uniform for all $k\le tn^{2/3}$.
{Following the same arguments as in the previous proof, we then find that
\begin{align*}
\mathbb E\big[(T\n_k)^2\big]&=2\sum_{j\le k}\mathbb E\big[\Delta T\n_j\cdot T\n_{j-1}\big]+\sum_{j\le k}\mathbb E\big[(\Delta T\n_j)^2\big] \\
& \le \Big(\frac12\sqrt{\frac{n}{m}}\Big)^2k^2\big(1+o(1)\big) + k\cdot \frac12 \sqrt{\frac{n}{m}}(1+o(1)).
\end{align*}
Together with $\mathbb E[T\n_k]=\frac{k}{2}\sqrt{n/m}(1+o(1))$, we deduce that 
\[
mn^{-\frac73} \max_{k\le tn^{2/3}}\Var(T\n_k)\le o(1)+\cO(m^{\frac12}n^{-\frac76})\to 0. 
\]}
The convergence in~\eqref{eq: cv-T1} follows. 
\end{proof}

Let us now explain the proof strategy for the case (iii) in Proposition~\ref{prop: light-tri-ct}. We will rely upon the approximation and change of measure techniques adopted in the proof of Proposition~\ref{thm: cv-light}. 
For each $k\ge 1$, recall that $\check X_k$ is a Poisson$(\check U_k p)$ variable and 
$\check{\mathrm N}_{k, i}, i\ge 1$ are i.i.d.~Poisson$(\check V_k p)$ variables, with $\check U_k, \check V_k$ defined in~\eqref{def: checkUk}. 
From the coupling arguments in Section~\ref{sec: tv}, which are valid under~\eqref{hyp: L}, we deduce that there exists a coupling between $\{\check X\n_k, \check{\mathrm N}_{k, i}, 1\le i\le\check X\n_k, k\ge 1\}$ and $\{X\n_k, \#\cN_{k, i}, 1\le i\le X\n_k, k\ge 1\}$ so that 
\begin{equation}
\label{dtvT}
\mathbb P\Big(\exists\,k\le tn^{2/3}: \check X\n_k\ne X\n_k \text{ or } \check{\mathrm N}_{k, i}\ne \#\cN_{k, i} \text{ for some } i\le X\n_k\Big)\to 0.
\end{equation}
See~\eqref{bd: tv1},~\eqref{tvbd1} and~\eqref{bd: tvre} there. 
We introduce the check version of $T\n_k$: let 
\[
\Delta\check T\n_k := \check T\n_k-\check T\n_{k-1}:=\sum_{1\le i\le \check X\n_k}\tfrac16\big(1+\check{\mathrm N}_{k, i}\big)_3.
\]
Since both $\Delta S\n_k$ and $\Delta T\n_k$ are functions of $\#\cN_{k, i}, 1\le i\le X\n_k$, it suffices to prove a version of the case (iii) with $(T\n_{\lf tn^{2/3}\rf})_{t\ge 0}$ replaced by $(\check T\n_{\lf tn^{2/3}\rf})_{t\ge 0}$. We next move on to the hat processes. Recall $\alpha_n, \beta_n$ from~\eqref{def: betan}. 
Recall that $\hat X_k, k\ge 1$ are i.i.d~Poisson$(\alpha_n)$, $\hat N_{k, i}, i\ge 1, k\ge 1$ are i.i.d~Poisson$(\beta_n)$ variables, and $\Delta\hat S\n_k+1=\sum_{1\le i\le \hat X\n_k} \hat N_{k, i}$. Define 
\[
\Delta\hat T\n_k := \hat T\n_k-\hat T\n_{k-1}:=\sum_{1\le i\le \hat X\n_k}\tfrac16\big(1+\hat N_{k, i}\big)_3.
\]
The discrete Radon--Nikodym identity~\eqref{id: dens-comp} then implies that for any continuous bounded functions $F, G: \mathbb D^{\ast}([0, t], \R)\to \R$, 
\begin{align}\notag
& \mathbb E\Big[F\Big((\check S\n_{\lf sn^{2/3}\rf})_{s\le t}\Big) G\Big((\check T\n_{\lf sn^{2/3}\rf})_{s\le t}\Big)\Big] \\ \label{id: Tdens}
&\qquad\qquad = \mathbb E\Big[F\Big((\hat S\n_{\lf sn^{2/3}\rf})_{s\le t}\Big) G\Big((\hat T\n_{\lf sn^{2/3}\rf})_{s\le t}\Big)E\n_{\lf tn^{2/3}\rf}\Big(\big(\hat R\n_k, \hat S\n_k)_{k\le tn^{2/3}}\big) \Big] 
\end{align}
In the first step of the proof, we will establish the convergence of $(\hat T\n_k)_{k\ge 1}$. Namely, under~\eqref{hyp: L} with $m=\lf n^{\frac73}\rf$, we have
\begin{equation}
\label{cv: hat-T}
\Big\{\hat T\n_{\lf tn^{2/3}\rf}: t\ge 0\Big\} \ \Longrightarrow\  \{\mathscr T_{t}: t\ge 0\}, 
\end{equation}
where $\{\mathscr T_{t}: t\ge 0\}$ is the Poisson process introduced in Proposition~\ref{prop: light-tri-ct}. We next define the following modified version of $\hat S\n_k$: let $\tilde S\n_j=\sum_{k\le j}\Delta \tilde S\n_k$ with 
\[
\Delta \tilde S\n_k = 
\left\{\begin{array}{ll}
\Delta \hat S\n_k & \text{ if } \Delta \hat T\n_k=0\\
\Delta_k & \text{ if } \Delta \hat T\n_k>0
\end{array}\right.
\]
where $\Delta_k, k\ge 1$ are i.i.d.~random variables with the same distribution as $\hat S\n_1$ conditioned on $\hat T\n_1=0$. 
Recall that $\Delta \hat S\n_k$ are i.i.d. It then follows from the previous construction that $\Delta \tilde S\n_k, k\ge 1$ are also i.i.d.~and further independent of $\Delta\hat T\n_k, k\ge 1$. We then proceed to show that 
\begin{equation}
\label{cv: tilde-S}
\Big\{n^{-\frac13}\tilde S\n_{\lf tn^{2/3}\rf}: t\ge 0\Big\} \ \Longrightarrow\  \{\hat\cS_t: t\ge 0\} \quad \text{ in $\mathbb D^{\ast}(\R_+, \R)$}, 
\end{equation}
where $\hat\cS_t=W_t+2(\lambda-\lambda_{\ast})t$. As we have argued that $(\tilde S\n_k)$ is independent of $(\hat T\n_k)$, the convergences in~\eqref{cv: hat-T} and~\eqref{cv: tilde-S} hold jointly with the pair of the limit processes independent of each other.  We will further show that for any $\epsilon>0$, 
\begin{equation}
\label{eq: diff-hat-tilde}
\mathbb P\Big(\sup_{k\le tn^{2/3}}\big|\hat S\n_k-\tilde S\n_k\big|\ge \epsilon n^{\frac13}\Big)\to 0. 
\end{equation} 
It follows that 
\[
\Big\{n^{-\frac13}\hat S\n_{\lf tn^{2/3}\rf}, \hat T\n_{\lf tn^{2/3}\rf}: t\ge 0\Big\} \ \Longrightarrow\  \big\{\hat\cS_t, \mathscr T_t: t\ge 0\big\} \quad \text{ in $\mathbb D^{\ast}(\R_+, \R^{2})$}. 
\]
Combined with~\eqref{id: Tdens} and the convergence of $E\n_{\lf tn^{2/3}\rf}$ seen in Proposition~\ref{prop: cv-dens-light}, this yields
\begin{align*}
\mathbb E\Big[F\Big((n^{-\frac13}\check S\n_{\lf sn^{2/3}\rf})_{s\le t}\Big) G\Big((\check T\n_{\lf sn^{2/3}\rf})_{s\le t}\Big)\Big] \longrightarrow &\,\mathbb E\Big[F\Big((\hat \cS_s)_{s\le t}\Big)G\Big((\mathscr T_s)_{s\le t}\Big)\cE_t\Big((\hat\cS_s)_{s\le t}\Big)\Big]\\
= &\,\mathbb E\Big[F\Big((\hat \cS_s)_{s\le t}\Big)\cE_t\Big((\hat\cS_s)_{s\le t}\Big)\Big]\cdot \mathbb E\Big[G\Big((\mathscr T_s)_{s\le t}\Big)\Big]\\
= &\, \mathbb E\Big[F\Big((\cS^{\lambda,\infty}_s)_{s\le t}\Big)\Big]\cdot \mathbb E\Big[G\Big((\mathscr T_s)_{s\le t}\Big)\Big],
\end{align*}
where we have used the independence in the second line and~\eqref{id: abs-cont-light} in the third. Together with~\eqref{dtvT}, this will allow us to complete the proof of case (iii) in Proposition~\ref{prop: light-tri-ct}. So it remains to prove~\eqref{cv: hat-T},~\eqref{cv: tilde-S} and~\eqref{eq: diff-hat-tilde}. The following lemma collects some preliminary estimates that we need for the proof.

\begin{lem}
\label{lem: prep1}
Under~\eqref{hyp: L} with $m=\lf n^{\frac73}\rf$, the following statements are true. \begin{enumerate}[(a)]
\item
For any $t>0$, we have
\begin{equation}
\label{eq: nodoublejump}
\mathbb P\big(\exists\, k\le tn^{2/3}: \Delta\hat T\n_{k}\ge 2\big) \to 0, 
\end{equation}
and for each $k\ge 1$, 
\begin{equation}
\label{eq: jumprate}
 \mathbb P(\Delta\hat T\n_k\ge 1)\sim \tfrac12 mn^2p^3, \quad n\to\infty. 
\end{equation}
\item
We have 
\begin{equation}
\label{mm: codS}
\mathbb E\big[1+\hat S\n_1\,\big|\,\hat T\n_1\ge 1] = \cO(1), \quad n\to\infty.
\end{equation}
\item
We have 
\begin{equation}
\label{mm: tildeS}
\mathbb E\big[\tilde S\n_1\big]= \mathbb E\big[\hat S\n_1\big]+o(n^{-\frac13})\quad\text{and}\quad \Var\big(\tilde S\n_1\big) = \Var\big(\hat S\n_1\big)+o(1), \quad n\to\infty, 
\end{equation}
as well as $\sup_{n\ge 1}\mathbb E[(\tilde S\n_1+1)^4]<\infty$. 
\end{enumerate}
\end{lem}

\begin{proof}[Proof of Lemma~\ref{lem: prep1}]
We start with the following observation: 
\begin{align*}
\mathbb P\big(\Delta \hat T\n_{k}\ge 2\big) &\le \mathbb P\big(\exists\, i\le \hat X\n_{k} : \hat N_{k,  i}\ge 3\big)
 +\mathbb P\big(\exists\,  i<j\le \hat X\n_{k}: \hat N_{k, i}\ge 2, \hat N_{k, j}\ge 2\big)\\
 & \le  \mathbb E[\hat X\n_{k}]\cdot \mathbb P(\text{Poisson}(\beta_n)\ge 3)+ \mathbb E\big[(\hat X\n_{k})_2\big] \cdot \mathbb P(\text{Poisson}(\beta_n)\ge 2)^{2},
\end{align*}
Recall from~\eqref{def: betan} that  $\alpha_n=mp$ and $\beta_n=np-2\lambda_{\ast}m^{-1/2}n^{1/6}\le np$. 
Together with \eqref{bd: pois-tail} and~\eqref{id: poisson-moment'}, we find
\[
\mathbb P\big(\Delta\hat T\n_{k}\ge 2\big) =\cO(mn^3p^4)+\cO(m^2n^4p^6)= \cO(m^{-1}n) = \cO(n^{-\frac43}),
\]
The convergence in~\eqref{eq: nodoublejump} then follows. Next, we have 
\begin{equation}
\label{eq: mhht}
1-\mathbb P(\Delta\hat T\n_{k}\ge 1) = \mathbb E\Big[\mathbb P\big(\text{Poisson}(\beta_n)< 2\big)^{\hat X\n_k}\Big].
\end{equation}
Let us identify the asymptotic behaviour of $\mathbb P\big(\text{Poisson}(\beta_n)< 2\big)$. We note from~\eqref{bd: pois-tail} that 
\[
\tfrac12 e^{-\beta_n}\beta_n^2\le \mathbb P\big(\text{Poisson}(\beta_n)= 2\big)\le \mathbb P\big(\text{Poisson}(\beta_n)\ge 2\big)\le \tfrac12n^{2}p^{2}.
\]
Since $\beta_n=np(1+o(1))$ and $np\to 0$ under~\eqref{hyp: L}, this implies that 
\[
\mathbb P\big(\text{Poisson}(\beta_n)\ge 2\big)=\frac12 n^2p^2\big(1+o(1)\big).
\]
Recall that $\hat X\n_k\sim$ Poisson$(mp)$ and note that $\mathbb E[s^{\hat X\n_k}]=\exp(-(1-s)mp)$. It follows that
\[
 \mathbb E\Big[\mathbb P\big(\text{Poisson}(\beta_n)< 2\big)^{\hat X\n_k}\Big]=
\mathbb E\Big[\big(1-\tfrac12n^2p^2(1+o(1))\big)^{\hat X\n_k}\Big] = \exp\Big(-\frac12mn^2p^3(1+o(1))\Big).
\]
As $mn^2p^3\to 0$, we deduce~\eqref{eq: jumprate} from~\eqref{eq: mhht}. This completes the proof of (a).

Turning to (b), we note that on the event $\hat T\n_1>0$, we can find some $i\le X\n_1$ such that $\hat N_{1, i}\ge 2$. For $q\in \N$, denoting $H^q_i = \{\hat X\n_1=q, \hat N_{1, i}\ge 2, \hat N_{1, j}<2,  \forall\,j<i\}$, we have the following decomposition into a disjoint union: 
\[
\big\{\hat T\n_1>0, \hat X\n_1=q\big\}=\bigcup_{1\le i\le q}H^q_i.
\]
Hence, 
\begin{equation}
\label{deq: sum}
\mathbb E\Big[\sum_{j=1}^q\hat N_{1, j}\mathbf 1_{\{\hat T\n_1>0, \hat X\n_1=q\}}\Big] = \sum_{i=1}^q \sum_{j=1}^{i-1} \mathbb E[\hat N_{1, j}\mathbf 1_{H^q_i}] + \sum_{i=1}^q \sum_{j=i+1}^q \mathbb E[\hat N_{1, j}\mathbf 1_{H^q_i}] + \sum_{i=1}^q \mathbb E[\hat N_{1, i}\mathbf 1_{H^q_i}]. 
\end{equation}
Let $Y$ be a Poisson$(\beta_n)$ variable. Recall that $\hat N_{1, i}, i\ge 1$ are independent with the same distribution as $Y$. From~\eqref{bd: binom7''} we get  $\mathbb E[Y\,|\, Y<2]\le \beta_n$. If $j\le i-1$, then the previous bound combined with the independence between $(\hat N_{1, j})$ and with $\hat X\n_1$ implies that 
\[
\mathbb E[\hat N_{1, j}\mathbf 1_{H^q_i}] =\mathbb E[\hat N_{1, j}\,|\, H^q_i]\cdot\mathbb P(H^q_i)= \mathbb E[Y| Y<2]\cdot \mathbb P(H^q_i) \le\beta_n\cdot  \mathbb P(H^q_i).
\]
If $j\ge i+1$, then the independence implies that $\mathbb E[\hat N_{1, j}\mathbf 1_{H^q_i}] = \beta_n\cdot \mathbb P(H^q_i)$. For $i=j$, we apply~\eqref{bd: binom7'} to find that 
\[
\mathbb E [\hat N_{1, i}\mathbf 1_{H^q_i}] =\mathbb E[ Y\,|\, Y\ge 2] \cdot \mathbb P(H^q_i) \le 2e^{\beta_n} \mathbb P(H^q_i).
\]
Feeding all these into~\eqref{deq: sum}, we obtain that 
\[
\mathbb E\Big[\sum_{i=1}^q\hat N_{1, i}\mathbf 1_{\{\hat T\n_1>0, \hat X\n_1=q\}}\Big] \le (q\beta_n+2e^{\beta_n}) \sum_{j=1}^q \mathbb P(H^q_j)  = (q\beta_n+2e^{\beta_n})\mathbb P(\hat T\n_1>0, \hat X\n_1=q).  
\]
Summing over $q$ leads to 
\[
\mathbb E\big[\big(\hat S\n_1+1\big)\mathbf 1_{ \{\hat T\n_1>0\}}\big]\le \big(\mathbb E[\hat X\n_1]\beta_n+2e^{\beta_n}\big) \mathbb P(\hat T\n_1>0).
\]
As $\mathbb E[\hat X\n_1]=mp$,  $mp\beta_n \le mnp^2\to 1$ and $\beta_n\to 0$, this yields the bound in~\eqref{mm: codS}. 

For (c), we recall that $\tilde S\n_1$ is distributed as $\hat S\n_1$ conditioned on $\hat T\n_1=0$. As $\hat S\n_1\ge -1$, we have $|\hat S\n_1|\le \hat S\n_1+2$. 
We deduce from~\eqref{mm: codS} and~\eqref{eq: jumprate} that 
\[
\mathbb E[|\hat S\n_1|\,; \hat T\n_1>0] =\cO(1)\cdot \cO(mn^2p^3)=o(n^{-\frac13}). 
\]
It follows that 
\[
\mathbb E[\hat S\n_1; \hat T\n_1=0]=\mathbb E[\hat S\n_1]-\mathbb E[\hat S\n_1; \hat T\n_1>0] = \mathbb E[\hat S\n_1]+o(n^{-\frac13}).
\]
Together with~\eqref{eq: jumprate}, which implies that $1-\mathbb P(\hat T\n_1=0) = o(n^{-1/3})$, this yields 
\[
\mathbb E[\tilde S\n_1] = \frac{\mathbb E[\hat S\n_1; \hat T\n_1=0]}{\mathbb P(\hat T\n_1=0)} = \frac{\mathbb E[\hat S\n_1]+o(n^{-\frac13})}{1+o(n^{-1/3})}=\mathbb E[\hat S\n_1]+o(n^{-\frac13}) . 
\]
For the variance, we use the Cauchy--Schwarz inequality to find that 
\[
\mathbb E[(\hat S\n_1)^2; \hat T\n_1>0]\le \sqrt{ \mathbb E[(\hat S\n_1)^4]\cdot \mathbb P(\hat T\n_1>0)} = o(1),
\]
where we have used $\mathbb E[(\hat S\n_1)^4]=\cO(1)$, a consequence of~\eqref{bd: fourth}. 
Therefore, 
\[
\mathbb E\big[(\hat S\n_1)^2\,\big|\,\hat T\n_1=0\big]  = \frac{\mathbb E[(\hat S\n_1)^2]-\mathbb E[(\hat S\n_1)^2; \hat T\n_1>0]}{\mathbb P(\hat T\n_1=0)} = \frac{\mathbb E[(\hat S\n_1)^2]+o(1)}{1+o(1)}. 
\]
We have previously shown that $\mathbb E[\tilde S\n_1]=\mathbb E[\hat S\n_1]+o(1)$. Note from~\eqref{eq: moment-light} that $\mathbb E[\hat S\n_1]=\cO(1)$, so that $(\mathbb E[\tilde S\n_1])^2=(\mathbb E[\hat S\n_1])^2+o(1)$.  Combining this with the previous bound, we find that
\[
\Var(\tilde S\n_1)=\mathbb E[(\hat S\n_1)^2\,|\,\hat T\n_1=0]-\big(\mathbb E[\tilde S\n_1]\big)^2=\mathbb E[(\hat S\n_1)^2]-(\mathbb E[S\n_1])^2+o(1), 
\]
which yields the statement concerning $\Var(\tilde S\n_1)$ in~\eqref{mm: tildeS}.  
Finally, we have 
\[
\mathbb E\big[(\tilde S\n_1+1)^4\big] \le \big(\mathbb P(\hat T\n_1=0)\big)^{-1}\mathbb E\big[(\hat S\n_1+1)^4\big], 
\]
and is therefore uniformly bounded, since we have seen $\sup_n\mathbb E[(\hat S\n_1+1)^4]<\infty$ in~\eqref{bd: fourth}. This completes the proof. 
\end{proof}

\begin{proof}[Proof of Proposition~\ref{prop: light-tri-ct}: case (iii)]
Define $\tau\n_1=\min\{k\ge 1: \Delta\hat T\n_k\ge 1\}$, and for $q\ge 1$, $\tau\n_{q+1}=\min\{k> \tau\n_q: \Delta\hat T\n_k\ge 1\}$. In words, $\tau\n_q, q\ge 1$,  are the successive jump times of $(\hat T\n_k)_{k\ge 1}$. 
Take an arbitrary $s>0$ and write $N=\lf sn^{\frac23}\rf$.  We deduce from~\eqref{eq: jumprate} that
\[
\mathbb P\big(\tau\n_1 >N\big) =\mathbb P\big(\Delta\hat T\n_k=0, 1\le k\le N\big) = \big(1-\tfrac12 mn^2p^3(1+o(1))\big)^N  \to \exp(-\tfrac12 s).
\]
This shows $\tau\n_1$ converging in distribution to an exponential variable of mean $2$. Since $\Delta\hat T\n_k, k\ge 1$ are i.i.d,  
we can then extend this argument to $\tau\n_{q+1}-\tau\n_q$ to find that it is independent of $\tau\n_q$ and has the same limit distribution. Together with \eqref{eq: nodoublejump}, this shows that $(\hat T\n_{\lf tn^{2/3}\rf})_{t\ge 0}$ converges in distribution to a simple point process that jumps according to exponentially distributed waiting times, which proves~\eqref{cv: hat-T}. For the convergence in~\eqref{cv: tilde-S}, we 
note that in the course of proving Proposition~\ref{prop: cv-light-tree}, the key estimates that have allowed us to apply the martingale central limit theorem and deduce the convergence of $\hat S\n$ are the ones in~\eqref{eq: moment-light} and~\eqref{bd: fourth}. Lemma~\ref{lem: prep1} (c) show that the same estimates hold for $\tilde S\n_1$. The convergence in~\eqref{cv: tilde-S} then follows from the same arguments as seen in the proof of Proposition~\ref{prop: cv-light-tree}. 
Finally, for~\eqref{eq: diff-hat-tilde}, we find from the definition of $(\tilde S\n_k)$ that
\[
\sup_{k\le tn^{2/3}}\big|\tilde S\n_k-\hat S\n_k\big| \le \sum_{k\le tn^{2/3}}\big|\Delta\hat S\n_k-\Delta_k\big|\mathbf 1_{\{\Delta\hat T\n_k>0\}}. 
\]
Noting that $\Delta_k$ is independent of $\Delta T\n_k$, it follows that 
\begin{align*}
\mathbb E\Big[\sup_{k\le tn^{2/3}}\big|\tilde S\n_k-\hat S\n_k\big|\Big]&\le \sum_{k\le tn^{2/3}}\Big(\mathbb E\big[|\Delta\hat S\n_k|\,\big|\,\Delta T\n_k>0\big]+\mathbb E\big[|\Delta_k|\big]\Big)\cdot \mathbb P(\Delta\hat T\n_k>0)\\
& \le \Big(\mathbb E\big[|\hat S\n_1|\,\big|\,\hat T\n_1>0\big]+\mathbb E\big[|\Delta_1|\big]\Big)\cdot \sum_{k\le tn^{2/3}}\mathbb P(\Delta\hat  T\n_k>0). 
\end{align*}
Both $\mathbb E\big[|\hat S\n_1|\,|\,\hat T\n_1>0\big]$ and $\mathbb E\big[|\Delta_1|\big])$ are $\cO(1)$, respective consequences of~\eqref{mm: codS} and~\eqref{mm: tildeS}. Meanwhile,~\eqref{eq: jumprate} implies that the sum above is also $\cO(1)$. An application of the Markov inequality then yields: 
\[
\mathbb P\Big(\sup_{k\le tn^{2/3}}\big|\tilde S\n_k-\hat S\n_k\big| \ge \epsilon n^{\frac13}\Big) =\cO( n^{-\frac13}), 
\]
which proves~\eqref{eq: diff-hat-tilde}. 
The arguments prior to Lemma~\ref{lem: prep1} allow us to conclude. 
\end{proof}

We recall the bijection between the vertex labels: $\rho_{n}: [n]\to \cV$, which maps vertex $i$ of $G(n, m, p)$ to the $i$-th $\cV$-vertex $w_{i}$ of $B(n, m, p)$. Recall from Algorithm~\ref{alg} the sequence $\{v_i: 1\le i\le n\}$, the $\cV$-vertices in the order of exploration.
Let us denote by $\hat T\n_{k}$ the number of triangles in $G(n, m, p)$ that are formed by the vertices $\{\rho_{n}^{-1}(v_{i}): i\le k\}$. Recall from~\eqref{def: goodevent} the event $\cI\n_{k}$. 

\begin{lem}
\label{lem: diff-t}
Assume either \eqref{hyp: L} or \eqref{hyp: M}. 
For all $m, n\ge 1$, $p\in [0, 1]$ and  $t>0$, writing $\cI_{t}=\cap_{k\le tn^{2/3}}\cI\n_{k}$, we have 
\[
\mathbb E\Big[\max_{k\le tn^{2/3}}|\hat T\n_{k}-T\n_{k}|\,; \cI_{t}\Big]\le 2 t mn^{\frac53}p^{3}\cdot \sup_{k\le tn^{2/3}}\mathbb E[\#\cA_{k-1}] + 8t^{3}m^{3}n^{2}p^{6}.
\]
\end{lem}

\begin{proof}
Recall that the three vertices of a type I triangle belong to a common group, and non type I triangles are of type II. 
Let us introduce the following quantities:
\begin{align*}
\eta_1(n, m; k) &= \#\big\{(v_{j}, v_{j'}, u): u\in \cM_{k}, v_{j}\in \cB(u)\cap\cN_{k}, v_{j'}\in \cB(u)\cap \cA_{k-1}\big\},\\
\eta_2(n,m;k) &= \#\big\{(v_{j}, v_{j'}, u): 1\le j<j', u\in \cM_k, v_{j}\in\cB(u)\cap \cA_{k-1}, v_{j'}\in\cB(u)\cap \cA_{k-1}\big\},
\end{align*}
as well as 
\[
\gamma(n, m; k) = \#\big\{\text{triangles of type II with vertices from } \rho_{n}^{-1}(\{v_{1}, v_{2}, \cdots, v_{k}\})\big\}.
\]

We claim that
\begin{equation}
\label{eq: bd1-tri}
0\le \big(\hat T\n_{k}-T\n_{k}\big)\mathbf 1_{\cap_{j\le k}\cI\n_{j}} \le \gamma(n, m; k)+\sum_{j=1}^{k}\big(\eta_{1}(n, m; j)+\eta_{2}(n, m; j)\big). 
\end{equation}
To see why this is true, we first note that $\Delta T\n_{k}$ counts the number of type-I triangles with all three vertices belong to some $\{v_{k}\}\cup\cN_{k, i}$; thus $\hat T\n_{k}\ge T\n_{k}$. Since $\gamma(n, m; k)$ counts the type II triangles in $\hat T\n_{k}$, it remains to see that on the event $\cap_{j\le k}\cI\n_{j}$,  the number of type I triangles missing from $T\n_{k}$ is indeed bounded by the summation term in \eqref{eq: bd1-tri}.  
Let us consider such a triangle. 
Without loss of generality, we can assume that the three vertices of the triangles are $\rho_{n}^{-1}(v_{j_1}), \rho_{n}^{-1}(v_{j_2}), \rho_{n}^{-1}(v_{j_3})$ with $j_1<j_2<j_3\le k$. Since the triangle is of type I, $v_{j_1}, v_{j_2}, v_{j_3}$ has a common neighbour; let $u\in \cM_j$ be their common neighbour with the smallest index $j$ (if there are multiple $u$ satisfying the condition, choose an arbitrary one). We note that necessarily $j\le j_1$ as $u$ is adjacent to $v_{j_1}$. Suppose $j=j_{1}$; then there are four possibilities: (i) we can have both $v_{j_2}$ and $v_{j_3}$ belong to $\cN_{j}$; or (ii) $v_{j_2}\in \cN_{j}$ and $v_{j_3}\in \cA_{j-1}$; or (iii) $v_{j_3}\in \cN_j$ and $v_{j_{2}}\in \cA_{j-1}$; or (iv) both $v_{j_2}$ and $v_{j_3}$ belong to $\cA_{j-1}$. However, in case (i), the event $\cI\n_{j}$ ensures both $v_{j_{2}}$ and $v_{j_{3}}$ are contained in the same $\cN_{j, i^{\ast}}$ for some $i^{\ast}\le X\n_{j}$, and the triangle spanning $v_{k}, v_{j_{2}}, v_{j_{3}}$ would be counted in $\Delta T\n_j$. Case (iii) is also impossible, as $v_{j_3}$ would be ahead of $v_{j_2}$ in $\cA_j$ and be explored first as a result. The other two cases are counted respectively in $\eta_{1}(n, m; j)$ and $\eta_2(n, m; j)$. 
Suppose $j<j_{1}$. Then we can exclude the possibility that all three vertices belong to $\cN_j$ on account of $\cI\n_{j}$. Depending on whether one or at least two of them belong to $\cA_{j-1}$, this scenario is also accounted for by either $\eta_1(n, m; j)$ or $\eta_2(n, m; j)$. 
To sum up, we have shown that on $\cap_{j\le k}\cI\n_{j}$, the difference between $\hat T\n_{k}$ and $T\n_{k}$ is indeed bounded by the right-hand side of \eqref{eq: bd1-tri}. 

Next, let us show the following bounds:
\begin{equation}
\label{eq: eta-bd}
\sup_{k\le tn^{2/3}}\big(\mathbb E[\eta_{1}(n, m;k)]+\mathbb E[\eta_{2}(n, m;k)]\big)\le 2mnp^{3}\cdot\sup_{k\le tn^{2/3}}\mathbb E[\#\cA_{k-1}]. 
\end{equation}
and
\begin{equation}
\label{eq: gamma-bd}
\mathbb E[\gamma(n, m; k)] \le 8k^{3}m^{3}p^{6}.
\end{equation}
For the proof of \eqref{eq: eta-bd}, we note that conditional on $u\in \cM_k$ and on $\cA_{k-1}$, $\#(\cB(u)\cap \cA^c_{k-1})$ and $\#(\cB(u)\cap \cA_{k-1})$ are independent. Moreover, we have $\#(\cB(u)\cap \cN_k)\le \#(\cB(u)\cap \cA^c_{k-1})$, which is further stochastically bounded by a Binom$(n, p)$, and $\#(\cB(u)\cap \cA_{k-1})$ distributed as Binom$(\#\cA_{k-1}, p)$. It follows that
\[
\mathbb E[\eta_{1}(n, m; k)]\le \mathbb E\Big[\sum_{u\in \cM_{k}} \#(\cB(u)\cap \cN_{k}) \cdot \#(\cB(u)\cap \cA_{k-1})\Big]
\le mnp^{3}\,\mathbb E[\#\cA_{k-1}].
\] 
Meanwhile, writing $Y_{u}=\#(\cB(u)\cap \cA_{k-1})$, which has  Binom$(\#\cA_{k-1}, p)$ distribution, we have 
\begin{align*}
\mathbb E[\eta_{2}(n, m; k)]&\le \mathbb E\Big[\sum_{u\in \cM_{k}} Y_{u}(Y_{u}-1)\Big]
=\mathbb E\Big[\mathbb E\Big[\sum_{u\in \cM_{k}} Y_{u}(Y_{u}-1)\,\Big|\, \cM_k\Big]\Big]\\
&=p^2\mathbb E\Big[\# \cM_k\cdot (\#\cA_{k-1})_2\Big]\le p^2\mathbb E\Big[\# \cM_k\cdot (\#\cA_{k-1})^2\Big]\\
&\le mp^{3}\cdot \mathbb E[(\#\cA_{k-1})^{2}]\le mnp^{3}\,\mathbb E[\#\cA_{k-1}], 
\end{align*}
where we have used~\eqref{id: binom-moment} in the second line, $\#\cM_k$ dominated by Binom$(m, p)$ in the third line as well as $\#\cA_{k-1}\le n$.  
Together they give the desired bound in \eqref{eq: eta-bd}. Regarding \eqref{eq: gamma-bd}, we first observe that for a type II triangle, any two of its three vertices belong to a common group which does not contain the third one. As previously, suppose that the three vertices of such a triangle are $\rho_{n}^{-1}(v_{j}), \rho_{n}^{-1}(v_{j'}), \rho_{n}^{-1}(v_{j''})$ with $j<j'<j''\le k$. 
Then we have
\begin{align*}
&\cN:=\cB(v_{j})\cap \cB(v_{j'})\ne\varnothing, \quad \cN' := \big(\cB(v_{j})\!\setminus\! \cN\big)\cap \cB(v_{j''})\ne\varnothing,\\
& \cN'' := \big(\cB(v_{j'})\!\setminus\! \cN\big)\cap \big(\cB(v_{j''})\!\setminus\! \cN'\big)\ne\varnothing.
\end{align*}
Let $X$ be a $\Z_+$-valued random variable and  $Y\sim$ Binom$(X, p)$. Using 
\eqref{bd: binom-tail}, we have $\mathbb P(Y\ge 1\,|\,X)\le p^{2} X+p X$. Since $\cB(v_{j''})$ is independent of $\cB(v_{j'}), \cN, \cN'$, by first conditioning on $X:=\#(\cB(v_{j'})\setminus(\cN\cup\cN'))$ and $\cN, \cN'$, we deduce that
\begin{align*}
\mathbb P(\cN''\ne\varnothing\,|\, \cN, \cN') &\le \mathbb P(\text{Binom}(X, p)\ge 1\,|\, \cN, \cN')\\
&\le p^{2}\mathbb E[X\,|\,\cN,\cN']+p\mathbb E[X\,|\,\cN,\cN']\le mp^3+mp^2\le 2mp^2,
\end{align*}
where we have used the fact that $X$ is dominated by Binom$(m, p)$, even if conditional on $\cN, \cN'$. 
Similarly, both $\mathbb P(\cN\ne\varnothing)$ and $ \mathbb P(\cN'\ne \varnothing \,|\, \cN)$ are bounded by $2mp^2$. 
It follows that 
\[
\mathbb E(\gamma(n,m; k)) \le \sum_{j<j'<j''\le k} \mathbb P\big(\cN\ne\varnothing, \cN'\ne\varnothing, \cN''\ne\varnothing\big)\le k^{3}\cdot \big(\mathbb P(Y\ge 1)\big)^{3}\le 8 k^{3}(mp^{2})^{3},
\]
which proves \eqref{eq: gamma-bd}. Now \eqref{eq: bd1-tri}, \eqref{eq: eta-bd} and \eqref{eq: gamma-bd} put together complete the proof. 
\end{proof}

In view of the upper bound in Lemma~\ref{lem: diff-t}, we require the following estimate on $\mathbb E[\#\cA_{k}]$.
\begin{lem}
\label{lem: bdA}
Assume either \eqref{hyp: L} or \eqref{hyp: M}. 
For all $t>0$, we have
\[
\max_{k\le tn^{2/3}} \mathbb E\big[(\#\cA_{k})^{2}\big] =\cO(n^{\frac23}). 
\]
\end{lem} 

\begin{proof}%[Proof of Lemma \ref{lem: bdA}] 
Let us denote $N=\lf tn^{2/3}\rf$. Thanks to \eqref{id: Ak} and the fact that $-\min_{j\le k}S\n_j\le \max_{j\le k}|S\n_j|$, we only need to show that $\mathbb E\big[\max_{k\le N}(S\n_{k})^{2}\big] =\cO(n^{2/3})$. Let 
\[
\Delta E_k = \mathbb E\big[\Delta S\n_k\,\big|\, U_k, V_k\big], \quad E_k = \sum_{j\le k} \Delta E_j \quad\text{and}\quad M_k = S\n_k-E_k.  
\]
Note that $M_k, k\ge 0$ is a zero-mean martingale. The conclusion will follow once we show that 
$\mathbb E[\max_{k\le N} E_k^2]+\mathbb E[\max_{k\le N} M_k^2] = \cO(n^{2/3})$. Starting with $E_k$, we recall that $\Delta S\n_k+1=\sum_{1\le i\le X\n_k} \#\cN_{k, i}$, where $X\n_k\sim$ Binom$(U_k, p)$ and $\#\cN_{k, i}\sim$ Binom$(V_{k, i}, p)$ with $V_{k+1}\le V_{k, i}\le V_k$. It follows that  for all $k\ge 1$, 
\[
mnp^2-1-(mn-U_kV_{k+1})p^2=U_k V_{k+1}p^2-1\le \Delta E_k \le U_kV_kp^2-1\le mnp^2-1
\]
Summing over $k$ yields 
\[
0\le  k(mnp^2-1) - E_k \le \sum_{j\le k} (mn-U_jV_{j+1})p^2\le k(mn-U_kV_{k+1})p^2, 
\]
from which it follows that
\[
\max_{k\le N} \big|k(mnp^2-1)-E_k\big|\le N (mn-U_NV_{N+1})p^2. 
\]
Under either~\eqref{hyp: L} or~\eqref{hyp: M}, we have $mnp^2-1=\cO(n^{-\frac13})$. Moreover, Lemma~\ref{lem: Uk} tells us that 
\[
0\le \mathbb E[(mn-U_NV_{N+1})^2] \le 2\mathbb E[n^2(m-U_N)^2] + 2\mathbb E[U_N^2(n-V_{N+1})^2] = \cO(mn^{\frac73})+\cO(m^2n^{\frac43}). 
\]
Therefore,  $\mathbb E[(mn-U_NV_{N+1})^2p^4]=\cO(n^{-\frac23})$. Together with the previous bound, this implies that $\mathbb E[\max_{k\le N}E_k^2]=\cO(n^{\frac23})$. We next prove the desired bound for $(M_k)$. Note that
\[
\Var(\Delta M_k\,|\,U_k, V_k)\le \mathbb E\big[(\Delta S\n_k+1)^2\,|\, U_k, V_k\big].
\]
Note that each $\#\cN_{k, i}$ is stochastically bounded by a Binom$(V_k, p)$. 
We then develop the first term as follows: 
\begin{align*}
\mathbb E\big[(\Delta S\n_k+1)^2\,|\, U_k, V_k\big] &= \mathbb E\Big[\sum_{1\le i\le X\n_k}\#\cN_{k, i}^2\,\Big|\, U_k, V_k\Big] + 2\mathbb E\Big[\sum_{1\le i<j\le X\n_k} \#\cN_{k, i}\#\cN_{k, j}\,\Big|\, U_k, V_k\Big]\\
&\le U_kp (V_k^2p^2+V_kp) + U_k^2V_k^2p^4 \le mp(n^2p^2+np)+m^2n^2p^4 =\cO(1),
\end{align*}
thanks to the assumptions~\eqref{hyp: L} and~\eqref{hyp: M}. We have shown that $\max_{k\le N} \Var(\Delta M_k\,|\,U_k, V_k)=\cO(1)$. 
It follows that $\Var(M_N)= \cO(N)$. Doob's maximal inequality then allows us to conclude that
\[
\mathbb E\Big[\max_{k\le N}M_k^2\Big]\le 4\cdot \mathbb E[M_N^2]=\cO(n^{\frac23}).
\]
By the previous arguments, this completes the proof. 
\end{proof}

\begin{proof}[Proof of Theorem~\ref{thm: triangle-moderate}]
Combining Lemma~\ref{lem: diff-t} with Lemmas~\ref{lem: bdA} and~\ref{lem: goodevent} yields
\[
n^{-\frac23}\max_{k\le tn^{2/3}}\Big|\hat T\n_{k}-T\n_{k}\Big|\Big] \xrightarrow{n\to\infty}0, \quad  \text{in probability}. 
\]
Comparing this with Propositions \ref{prop: moderate-tri-ct}, we obtain the following convergence in probability:
\begin{equation}
\label{cv: tri-ct}
\Big\{n^{-\frac23}\hat T\n_{\lf tn^{2/3}\rf}: t\ge 0\Big\} \xrightarrow{n\to\infty} \{c_{\theta}\cdot t: t\ge 0\} \quad\text{in }\mathbb D^{\ast}(\R_{+}, \R), 
\end{equation}
where $c_{\theta}$ is the constant defined in \eqref{def: ctheta}. 
Let $(\hat{\mathrm g}\n_k, \hat{\mathrm d}\n_k)$ be the $k$-th longest excursion of $S\n$ above its running infimum. Then similar to~\eqref{cv: exc-int}, we also have
\begin{equation}
\label{cv: exc-int'}
\Big\{\Big(n^{-\frac23}\hat{\mathrm g}\n_k, n^{-\frac12}\hat{\mathrm d}\n_k\Big): k\ge 1\Big\} \Longrightarrow \big\{(\mathrm g^{\lambda, \theta}_k, \mathrm d^{\lambda, \theta}_k): k\ge 1\big\}
\end{equation}
with respect to the product topology of $(\R^2)^{\N}$. 
It follows that~\eqref{cv: tri-ct} also holds jointly with~\eqref{cv: exc-int'}. By Skorokhod's representation theorem and a diagonal argument, we can assume that both the convergences \eqref{cv: tri-ct} and \eqref{cv: exc-int'} take place a.s. We then deduce that for each $k\ge 1$, 
\[
n^{-\frac23}\Big(\hat T\n_{\hat{\mathrm d}\n_{k}}-\hat T\n_{\hat{\mathrm g}\n_{k}}\Big) \xrightarrow[]{n\to\infty} c_{\theta}\cdot \zeta^{\lambda,\theta}_{k}, \quad\text{a.s.}
\]
Note that the left-hand side counts the number of triangles in the $k$-th largest component of $G(n,m,p)$. The conclusion now follows. 
\end{proof}

\begin{proof}[Proof of Theorem~\ref{thm: triangle-light}]
We first observe that under \eqref{hyp: L}, Lemma~\ref{lem: diff-t} with Lemmas~\ref{lem: bdA} and~\ref{lem: goodevent} combined yield
\[
\max_{k\le tn^{2/3}}\Big|\hat T\n_{k}-T\n_{k}\Big| \to 0, \quad \text{in probability}. 
\]
For (i) and (ii), rest of the proof is similar to that of Theorem~\ref{thm: triangle-moderate}, with Proposition~\ref{prop: moderate-tri-ct} replaced by Proposition~\ref{prop: light-tri-ct}; we omit the details. For (iii), independence between $\mathscr T:=(\mathscr T_{t})_{t\ge 0}$ and $\cS^{\lambda, \infty}$ allows us to assert the joint convergence of $(\hat T\n_{\lf tn^{2/3}\rf})_{t\ge 0}\Rightarrow \mathscr T$ and $(n^{-\frac23}\hat{\mathrm g}\n_{k}, n^{-\frac23}\hat{\mathrm d}\n_{k})_{k\ge 1}\Rightarrow (\mathrm g^{\lambda, \infty}_{k}, \mathrm d^{\lambda, \infty}_{k})_{k\ge 1}$. Moreover, since $\mathscr T$ has no fixed discontinuity points and is independent of $(\mathrm g^{\lambda, \infty}_{k}, \mathrm d^{\lambda, \infty}_{k})_{k\ge 1}$, we can take compositions and deduce that  
\[
\Big\{\hat T\n_{\hat{\mathrm d}\n_{k}}-\hat T\n_{\hat{\mathrm g}\n_{k}}: k\ge 1\Big\} \ \Longrightarrow\ \big\{\mathscr T_{\mathrm d^{\lambda, \infty}_{k}}-\mathscr T_{\mathrm g^{\lambda, \infty}_{k}}: k\ge 1\big\}.
\]
Using the aforementioned independence again, we conclude that the right-hand side is a collection of Poisson$(\tfrac12\zeta^{\lambda, \infty}_{k})$ variables.   
\end{proof}

\subsection{Proof in the heavy clustering regime}
\label{sec: heavy}

Here, we give the proof of Theorems~\ref{thm3} and~\ref{thm: triangle-heavy}. 
As remarked  previously, we study the graph $B(m, n, p)$ instead, which allows us to apply the conclusions from the light clustering regime. In particular, let $\hat C\m_{k}$ be the connected component in $B(m, n, p)$ that has the $k$-th largest number of $\cV$-vertices (breaking ties arbitrarily), and denote by $\hat\mu\m_{k}$ the counting measure for $\cV$-vertices on $\hat C\m_{k}$. Then Theorem~\ref{thm1} tells us that 
\begin{equation}
\label{eq: cv-heavy'}
\Big\{\Big(\hat C\m_{k}, \, \tfrac12 m^{-\frac13}\cdot \dgr^{\bi}, \, m^{-\frac23}\hat\mu\m_{k}\Big): k\ge 1\Big\} \;\Longrightarrow\; \mathcal G^{\mathrm{ER}}(\lambda)
\end{equation}
in the same sense of convergence as \eqref{eq: cv-thm2}. On the other hand, there is a dual procedure to the one described in Section~\ref{sec: intro} that builds $G(n, m, p)$ from $B(m, n, p)$. To do so, we should regard $\cU$-vertices as individuals and $\cV$-vertices as communities. 
Let $\tilde\mu\m_{k}$ be the counting measure for $\cU$-vertices on $\hat C\m_{k}$. 
Recall that for the current discussion we have $\#\cU=n$ and $\#\cV=m$. 
The convergence in~\eqref{eq: cv-heavy'} is close to the one claimed in Theorem~\ref{thm3} apart from two aspects: first, ranking of the components in Theorem~\ref{thm3} corresponds to ranking in $\cU$-vertices rather than $\cV$-vertices; secondly, in the theorem, we use counting measures on the $\cU$-vertices with a different scaling $m^{1/6}n^{1/2}$. This scaling can be understood this way: in a first order approximation, the numbers of $\cU$-neighbours for each $\cV$-vertex is of the order $np\sim \sqrt{n/m}$; thus for a connected component that contains $m^{2/3}$ $\cV$-vertices, there are approximately $m^{1/6}n^{1/2}$ $\cU$-vertices. More specifically, 
Theorem~\ref{thm3} will follow from~\eqref{eq: cv-heavy'} once we show two things: first, with high probability, the ranking of the connected components  in $\cU$-vertices yields the same ranking as in $\cV$-vertices; secondly, in \eqref{eq: cv-heavy'} above, we can replace $m^{-2/3}\hat\mu\m_{k}$ by $m^{-\frac16}n^{-\frac12} \tilde\mu\m_{k}$. Both are proved in the following proposition. 

\begin{prop}
Assume \eqref{hyp: H}.  
The following holds for each $k\ge 1$: 
\[
\mathbb P\Big(\tilde\mu\m_1(\hat C\m_{1})>\tilde\mu\m_2(\hat C\m_{2})>\cdots>\tilde\mu\m_k(\hat C\m_{k})> \sup_{j>k} \tilde\mu\m_j(\hat C\m_{j}) \Big) \to 1, \quad m\to\infty. 
\]
Denoting $\delta^{\mathrm{Pr}}_{m, k}$ for the Prokhorov distance for finite measures on $(\hat C\m_{k}, m^{-\frac13}\dgr^{\bi}$). 
For each $k\ge 1$, we also have 
\begin{equation}
\label{cv: u-mass}
\delta^{\mathrm{Pr}}_{m, k}\big(m^{-\frac16}n^{-\frac12}\tilde\mu\m_k, m^{-\frac23}\hat\mu\m_{k}\big)\to 0 \quad \text{in probability}.
\end{equation}
\end{prop} 

\begin{proof}
We first note that  $\Delta R\m_{k} = X\m_{k}$ follows a Binom$(U_{k}, p)$ law, which is stochastically dominated by Binom$(n, p)$.~\eqref{bd: cc-binom} implies that 
for any $\epsilon>0$, 
\begin{align}\notag
\mathbb P\Big(\exists\,k_{1}<k_{2}\le m: R\m_{k_{2}} - R\m_{k_{1}} > (1+\epsilon) (k_{2}-k_{1}) np\Big) &\le \mathbb P(\exists\, k\le m: \Delta R\m_{k} > (1+\epsilon) np ) \\ \label{eqq1}
&\le m\exp\Big(-\frac{\epsilon^{2}np}{2+2\epsilon/3}\Big)\to 0,
\end{align}
as $np\to\infty$ under~\eqref{hyp: H}. Let $(\hat{\mathrm g}\m_{i}, \hat{\mathrm d}\m_{i})$ denote the $i$-th longest excursion in $(S\m_{k})_{k\ge 1}$. 
Since $R\m_{\hat{\mathrm d}\m_{i}} - R\m_{\hat{\mathrm g}\m_{i}} = \tilde\mu\m_{i}(\hat C\m_{i})$ and $\hat{\mathrm d}\m_{i}-\hat{\mathrm g}\m_{i}=\hat\mu\m_{i}(\hat C\m_{i})$, we have in particular 
\begin{equation}
\label{eq: tail1}
\mathbb P\Big(\exists\, i\ge 1: \tilde\mu\m_{i}(\hat C\m_{i})\ge 2 \hat\mu\m_{i}(\hat C\m_{i})np\Big) \to 0.
\end{equation}
On the other hand, let $\delta\in (0, \epsilon)$; we deduce again from~\eqref{bd: cc-binom} that for each $t\ge 0$, 
\begin{align*}
&\mathbb P\Big(\exists\,k_{1}<k_{2}\le tm^{\frac23}: R\m_{k_{2}} - R\m_{k_{1}} < (1-\epsilon) (k_{2}-k_{1}) np\,; U_{\lf tm^{2/3}\rf} \ge (1-\delta)n\Big)\\ 
&\qquad\qquad\qquad\le  tm^{\frac23}\mathbb P(\text{Binom}((1-\delta)n, p)<(1-\epsilon)np)\le tm^{\frac23}\exp\Big(-\frac{(\epsilon-\delta)^2np}{2(1-\delta)}\Big)\to 0.
\end{align*}
Switching the roles of $n$ and $m$ in Lemma~\ref{lem: Uk} yields  the bound relevant to the heavy clustering regime: 
\begin{equation}
\label{bd: Uk-heavy}
\max_{k\le t m^{2/3}}\mathbb E[(n-U_{k})^{2}] = \cO(m^{\frac13}n). 
\end{equation}
In particular, this yields $\mathbb P(U_{\lf tm^{2/3}\rf}<(1-\delta)n)\to 0$. Combining with the previous argument, we find that
\begin{equation}
\label{eqq2}
\mathbb P\Big(\exists\,k_{1}<k_{2}\le tm^{\frac23}: R\m_{k_{2}} - R\m_{k_{1}} < (1-\epsilon) (k_{2}-k_{1}) np\Big)\to 0
\end{equation}
 Note that the arguments in Section~\ref{sec: cv-ghp} also imply that 
\[
\big\{ m^{-\frac23}(\hat{\mathrm g}\m_{i}, \hat{\mathrm d}\m_{i}): i\ge 1\big\} \Longrightarrow \{(\mathrm g^{\lambda, \infty}_{i}, \mathrm d^{\lambda, \infty}_{i}): i\ge 1\}
\]
with respect to the product topology. Combined with the convergence of $R\m$ in~\eqref{cv: heavy}, this implies that for each $k\ge 1$, 
\[
\Big\{m^{-\frac16}n^{-\frac12}(R\m_{\hat{\mathrm d}\m_i}-R\m_{\hat{\mathrm g}\m_i}): 1\le i\le k\Big\}  \Longrightarrow \{\zeta^{\lambda, \infty}_i: 1\le i\le k\}.
\]
Since the left-hand side corresponds to $m^{-1/6}n^{-1/2} \tilde\mu\m_{i}(\hat C\m_{i})$ and $\zeta^{\lambda,\infty}_{1}>\zeta^{\lambda,\infty}_{2}>\zeta^{\lambda,\infty}_{3}>\cdots$ a.s, this implies in particular that for each $k\ge 1$, 
\begin{equation}
\label{ord1}
\mathbb P\Big(\tilde\mu\m_1(\hat C\m_{1})>\tilde\mu\m_2(\hat C\m_{2})>\cdots>\tilde\mu\m_k(\hat C\m_{k})\Big)\to 0
\end{equation}
The convergence \eqref{eq: cv-l2} applied to $(S\m_{k})_{k\ge 1}$ implies that we can find $N=N(k)\in \N$ so that
\[
\limsup_{m\to\infty}\mathbb P\big(\hat\mu\m_{N}(\hat C\m_{N})> \tfrac14 \hat\mu\m_{k}(\hat C\m_{k})\big) \to 0.  
\]
Combined with~\eqref{eq: tail1}, this yields 
\[
\limsup_{m\to\infty}\mathbb P\big(\tilde\mu\m_{N}(\hat C\m_{N})> \tilde\mu\m_{k}(\hat C\m_{k})\big) \to 0.  
\]
We note that~\eqref{ord1} also holds true for $N=N(k)>k$. Combined with the previous argument, this proves the first statement. 
For \eqref{cv: u-mass}, let us first take $A_{\ell, \ell'}=\{v_{\ell}, v_{\ell+1},\dots, v_{\ell'}\}\subseteq\hat C\m_{k}$. 
Denote by $\cB(A_{\ell, \ell'}) =\cup_{w\in A}\cB(w)$ its 1-neighbourhood in $B(m,n,p)$. We note that its 1-neighbourhood in the spanning forest $\cF$, with the exception of its parent, is contained in $\cup_{\ell\le j\le \ell'} \cM_{j}$. Let $p\m_{k}$ be the number of surplus edges in $\hat C\m_{k}$ and recall that $\Delta R\m_{k}=\#\cM_{k}$. Then we have
\begin{equation}
\label{bd: mea}
R\m_{k_{2}}-R\m_{k_{1}}\le \tilde\mu\m_{k}\big(\cB(A_{s, s'})\big)\le R\m_{k_{2}}-R\m_{k_{1}}+1+p\m_{k}
\end{equation}
We note that for a fixed $k$, $(p\m_k)_{k\ge 1}$ is tight. 
Applying  \eqref{eqq1} and noting that $k_2-k_1+1=\hat\mu\m_{k}(A_{\ell, \ell'})$, we find: 
\[
\mathbb P\Big(m^{-\frac16}n^{-\frac12}\big(R\m_{\ell'}-R\m_{\ell}+1+p\m_k\big)\ge m^{-\frac23}(1+\epsilon) \hat\mu\m_k(A_{\ell, \ell'})\Big)\to 0. 
\]
Similarly, we have $\limsup_m \mathbb P(\hat{\mathrm d}\m_k>tm^{2/3})\le \epsilon$ for some $t>0$, Applying~\eqref{eqq2} to $\ell\le\ell'\le \hat{\mathrm d}\m_k$, we then find that
\[
\mathbb P\Big(m^{-\frac16}n^{-\frac12}\big(R\m_{\ell'}-R\m_{\ell}\big)\le  m^{-\frac23}(1-\epsilon) \hat\mu\m_k(A_{\ell, \ell'})\Big)\to 0. 
\]
Together with~\eqref{bd: mea}, we deduce that 
\[
\sup  \big| m^{-\frac16}n^{-\frac12}\tilde\mu\m_{k}\big(\cB(A_{\ell, \ell'})\big) - m^{-\frac23}\hat\mu\m_{k}(A_{\ell, \ell'}) \big|\to 0 \quad \text{in probability}, 
\]
where the supremum is over all such $A_{\ell, \ell'}\subseteq \hat C\m_k$. A general subset $A$  can be treated similarly after decomposition. This then allows us to derive the desired Prokhorov bound.
\end{proof}

For the triangle counts, since the $\cV$-vertices serve as communities in this context, the analogue for $(T\n_{k})_{k\ge 1}$  is now defined as
\[
T\m_{k} = \sum_{i\le k}\tfrac16 (X\m_{k})_{3},\quad k\ge 0,
\]
where we recall $X\m_{k}=\#\cB(v_{k})\cap \cU_{k}$ is the degree of $v_{k}$ in $\cF$. 

\begin{prop}
Under the assumption~\eqref{hyp: H},
the following convergence takes place in probability: 
\[
\Big\{m^{\frac56}n^{-\frac32}T\m_{\lf tm^{2/3}\rf}: t\ge 0\Big\} \xrightarrow{n\to\infty} \big\{\tfrac16 t: t\ge 0\big\} \quad\text{in }\mathbb D^{\ast}(\R_{+}, \R).
\]
\end{prop}

\begin{proof}
Conditional on $U_{k}$, $X\m_{k}$ is distributed as a Binom$(U_{k}, p)$ variable, so that $\mathbb E[(X\m_{k})_{3}]=\mathbb E[(U_{k})_{3}]p^{3}$. It then follows from~\eqref{bd: Uk-heavy} that 
\[
\max_{k\le tm^{2/3}} \Big| n^{3}p^{3}-\mathbb E\big[(X\m_{k})_{3}\big]\Big| =o(n^3p^3), \; \max_{k\le tm^{2/3}}\Var\big[(X\m_{k})_{3}\big]  \le \max_{k\le tm^{2/3}}\mathbb E\big[(X\m_{k})^{6}\big]= \cO\big((np)^{6}\big). 
\]
Since $m^{2/3}n^{3}p^{3}\sim n^{3/2}m^{-5/6}\to\infty$, this implies that 
\[
\sup_{s\le t}\Big|m^{\frac56}n^{-\frac32}\mathbb E\Big[T\m_{\lf sm^{2/3}\rf}\Big] - \tfrac16 s\Big|\to 0 \quad \text{and} \quad m^{\frac53}n^{-3}\max_{k\le t m^{2/3}}\Var\big(T\n_{k}\big) \to 0.
\]
The conclusion then follows. 
\end{proof}

\begin{proof}[Proof of Theorem~\ref{thm: triangle-heavy}]
We denote by $\hat\rho_{m}: [n]\to \cU$ the bijection that maps the vertex $i$ of $G(n, m, p)$ to $w_{m+i}$ of $B(m, n, p)$. 
Let $\hat T\m_{k}$ be the number of triangles in $G(n, m, p)$ formed by the vertices from $\hat\rho_{m}^{-1}(\cup_{j\le k}\cM_{j})$. Let us show that for any $t>0$, 
\begin{equation}
\label{bd: diff-t'}
m^{\frac56}n^{-\frac32}\max_{k\le tm^{2/3}}\big|\hat T\m_{k}-T\m_{k}\big|\to 0, \quad \text{in probability}.
\end{equation}
To that end, we introduce the event
\[
\hat \cI_{t}: =\big\{\forall\, 1\le j< k\le tm^{2/3}: \#(\cB(v_{k})\cap \cM_{j})< 2\big\}.
\]
We deduce from~\eqref{bd: binom-tail} that 
\begin{align}\notag
1-\mathbb P(\hat\cI_{t}) & \le tm^{\frac23}\sum_{j\le tm^{2/3}}\mathbb P\big(\text{Binom}(\#\cM_{j}, p)\ge 2\big) \le tm^{\frac23}\sum_{j\le tm^{2/3}}(\mathbb E[\#\cM_{j}] p^{2}+\mathbb E[(\#\cM_{j})^{2}]p^{2})\\ \label{bd: goodevent'}
& \le t^{2}m^{\frac43}(2np^{3}+n^{2}p^{4})  \to 0. 
\end{align}
where we have used that $\#\cM_{j}$ is stochastically bounded by a Binom$(n, p)$ variable. 
In the same manner as in the proof of Lemma~\ref{lem: diff-t}, we divide triangles into two types: type I refers to those whose vertices belong to a common community, type II refers to the opposite. If a type I triangle is missing from the count $T\n_{k}$, then letting $v_{k}$ be the smallest community that contains all its vertices, we must have at least one of its vertices from $\cB(v_{k})\setminus \cM_{k}$. We denote $Y\m_{k}=\sum_{j<k}\mathbf 1_{\{\cM_{j}\cap\cB(v_{k})\ne\varnothing\}}$. 
On the event $\hat\cI_{t}$, depending on how many of the vertices come form $\cB(v_{k})\setminus \cM_{k}$, the number of missing type I triangles generated by $v_{k}$ is bounded by
\[
\eta (m, n; k):= Y\m_{k}\cdot (X\m_{k})^{2}+(Y\m_{k})_{2}\cdot X\m_{k}+(Y\m_{k})_{3}.
\]
Let us show that $\sum_{k\le tm^{2/3}}\mathbb E[\eta(m, n; k)]=o(n^{\frac32}m^{-\frac56})$. We start with the first term. Let $\hat X_k=\#\cB(v_k, p)$, which has a Binom$(n-1, p)$ by the definition of $B(m, n, p)$. Given $\hat X_k$, we have $\hat X\m_k=\#\cM_k\le \hat X_k$ a.s.~and 
\[
\mathbb E\big[Y\m_k\,\big|\,\hat X_k\big] \le k\cdot \mathbb P(\text{Binom}(\hat X_k, p)\ge 1)\le k \hat X_k(p+p^2)
\]
by~\eqref{bd: binom-tail}. Taking the expectation yields 
\[
\max_{k\le tm^{2/3}}\mathbb E\big[Y\m_k\cdot (X\m_k)^2\big] \le k(p+p^2) \mathbb E[(\hat X_k)^3] = \cO(m^{\frac23}n^3p^4),
\]
since $\mathbb E[(\hat X_k)^3]\sim (np)^3$; see Appendix~\ref{sec: binom}. 
Similarly, we have
\[
\mathbb E\big[(Y\m_k)_2\,\big|\,\hat X_k\big]\le 2\sum_{j<j'<k} \mathbb P(\text{Binom}(\hat X_k, p)\ge 1)^2\le k^2 \hat X_k^2 (p+p^2)^2,
\]
so that $\max_{k\le tm^{2/3}}\mathbb E[(Y\m_k)_2 X\m_k]=\cO(m^{\frac43}n^3p^5)$. We also have  $\max_{k\le tm^{2/3}}\mathbb E[(Y\m_k)_3]=\cO(m^2n^3p^6)$. 
This leads to
\[
\sum_{k\le tm^{2/3}}\mathbb E\big[\eta (m, n; k)\big]=\cO(m^{\frac43}n^{3}p^{4}+m^{2}n^{3}p^{5}+m^{\frac83}n^{3}p^{6}) = o(n^{\frac23}m^{-\frac65}). 
\]
For type II triangles, the same arguments used in the proof of Lemma~\ref{lem: diff-t} show that
\[
\mathbb E[\#\{\text{type II triangles formed by vertices from }\hat\rho_{m}^{-1}(\cup_{j\le k}\cM_{j})\} ]\le  8k^{3}(np^{2})^{3}. 
\]
Combining all previous bounds yields
\[
\mathbb E\Big[\max_{k\le tm^{2/3}}\big|\hat T\m_{k}-T\m_{k}\big|\,;\hat\cI_{t}\Big]\le \cO(nm^{-\frac23}+m^{2}n^{3}p^{6}).
\]
Together with~\eqref{bd: goodevent'}, this implies~\eqref{bd: diff-t'}. The rest of the proof is similar to that of Theorem~\ref{thm: triangle-moderate}. 
\end{proof}

\appendix
\section{An application of Girsanov's Theorem} 
\label{sec: Girsanov}

\begin{proof}[Proof of \eqref{id: abs-cont}]
We use the following version of Girsanov's Theorem. Let $(\mathbf W_{t})_{t\ge 0}$ be a bi-dimensional standard Brownian motion defined on the probability space $(\Omega, \mathcal F, \mathbf P)$ and denote by $(\cG_{t})_{t\ge 0}$ is its natural  filtration. Suppose that $(\boldsymbol \xi_{t})_{t\ge 0}$ is a {\it deterministic} continuous function taking values in $\R^{2}$. Then 
\[
Z_{t}:= \exp\Big(\int_{0}^{t} \boldsymbol \xi_{s}\cdot d\mathbf W_{s} - \frac12 \int_{0}^{t} \|\boldsymbol \xi_{s}\|^{2} ds\Big), \quad t\ge 0
\]
is a positive martingale with respect to $(\cG_{t})_{t\ge 0}$ with unit mean. The dot notation above denotes the inner product between vectors. Define a probability measure $\mathbf Q$ on $(\Omega, \mathcal F)$ by
\[
\frac{d\mathbf Q}{d\mathbf P}|_{\cG_{t}} = Z_{t}, \quad t\ge0.
\]
Denote $\widetilde{\mathbf W}_{t} = \mathbf W_{t} - \int_{0}^{t} \boldsymbol\xi_{s} ds$. 
Then under $\mathbf Q$, $(\widetilde{\mathbf W}_{t})_{t\ge 0}$ is a bi-dimensional standard Brownian motion. 
Namely, for all $t\ge 0$ and measurable functional $F: \mathbb C([0, t], \R^{2})\to \R_+$, 
\begin{align*}
\mathbf P\Big[F\Big(\big(\mathbf W_{s})_{s\le t}\Big)\Big] &=
\mathbf P\Big[F\Big(\big(\widetilde{\mathbf W}_{s})_{s\le t}\Big)\, \exp\Big(\int_{0}^{t}\boldsymbol\xi_{s} \cdot d\mathbf W_{s} -\frac12\int_{0}^{t} \|\boldsymbol\xi_{s}\|^{2}ds\Big)\Big] \\
& = \mathbf P\Big[F\Big(\big(\widetilde{\mathbf W}_{s})_{s\le t}\Big)\, \exp\Big(\int_{0}^{t}\boldsymbol\xi_{s}\cdot d\widetilde{\mathbf W}_{s} +\frac12\int_{0}^{t} \|\boldsymbol\xi_{s}\|^{2}ds\Big)\Big]. 
\end{align*}
Here, we apply it to 
\[
\mathbf W_t = 
\begin{pmatrix}
W^{\ast}_{t}\\
W_{t}
\end{pmatrix} \quad\text{and}\quad 
\begin{pmatrix}
\mathcal U_{t}\\
\mathcal V_{t}
\end{pmatrix}
=\begin{pmatrix}
W^{\ast}_{t}+\theta^{\frac14} \lambda t -\frac12 \theta^{-\frac14} t^{2}\\
W_{t}+\lambda t-\frac12 t^{2}
\end{pmatrix}
\]
with 
\[
\boldsymbol\xi_{t} = 
\begin{pmatrix}
\theta^{-\frac14}t -\theta^{\frac14}\lambda, \\
t-\lambda
\end{pmatrix}
\]
and we find that  
\[
\mathbb E\big[F\big((\mathcal U_{s}, \mathcal V_{s})_{s\le t}\big)\big] = \mathbb E\Big[F\big((W^{\ast}_{s}, W_{s})_{s\le t}\big) \exp\Big(-\int_{0}^{t} \boldsymbol\xi_{s} \cdot d\mathbf W_{s}-\frac12\int_{0}^{t} \|\boldsymbol\xi_{s}\|^{2} ds\Big)\Big]. 
\]
 Also define
\[
\begin{pmatrix}
\hat{\mathcal U}_{t}\\
\hat{\mathcal V}_{t}
\end{pmatrix}
=\begin{pmatrix}
W^{\ast}_{t}+\theta^{\frac14} (\lambda-\lambda_{\ast}) t \\
W_{t}+(\lambda-\lambda_{\ast}) t
\end{pmatrix}
\]
A similar argument but this time with
\[
\hat{\boldsymbol\xi}_{t} = 
\begin{pmatrix}
\theta^{\frac14}(\lambda_{\ast}-\lambda)\\
\lambda_{\ast}-\lambda
\end{pmatrix}
\]
implies that 
\[
\mathbb E\Big[F\big((\hat{\mathcal U}_{s}, \hat{\mathcal V}_{s})_{s\le t}\big)\exp\Big(\int_{0}^{t} \hat{\boldsymbol\xi}_{s}\cdot d
\begin{pmatrix}
\hat{\mathcal U}_{s}\\
\hat{\mathcal V}_{s}
\end{pmatrix}
+\frac12\int_{0}^{t} \|\hat{\boldsymbol\xi}_{s}\|^{2} ds\Big)\Big]=\mathbb E\Big[F\big((W^{\ast}_{s}, W_{s})_{s\le t}\big)\Big]. 
\]
Combining the two, we arrive at
\[
\mathbb E\big[F\big((\mathcal U_{s}, \mathcal V_{s})_{s\le t}\big)\big] = \mathbb E\Big[F\big((\hat{\mathcal U}_{s}, \hat{\mathcal V}_{s})_{s\le t}\big)\exp\Big(\int_{0}^{t} (\hat{\boldsymbol\xi}_{s}-{\boldsymbol\xi}_{s}) \cdot d\begin{pmatrix}
\hat{\mathcal U}_{s}\\
\hat{\mathcal V}_{s}
\end{pmatrix}
+\frac12\int_{0}^{t} (\|\hat{\boldsymbol\xi}_{s}\|^{2}-\|{\boldsymbol\xi}_{s}\|^{2}) ds\Big)\Big]
\]
Let us denote
\begin{align*}
\log \hat{\cE}_t:&=\int_{0}^{t} (\hat{\boldsymbol\xi}_{s}-{\boldsymbol\xi}_{s}) \cdot d\begin{pmatrix}
\hat{\mathcal U}_{s}\\
\hat{\mathcal V}_{s}
\end{pmatrix}
+\frac12\int_{0}^{t} (\|\hat{\boldsymbol\xi}_{s}\|^{2}-\|{\boldsymbol\xi}_{s}\|^{2}) ds\\
&=\int_0^t \Big\{(\theta^{\frac14}\lambda_{\ast}-\theta^{-\frac14}s)d\,\hat{\mathcal U}_s+(\lambda_{\ast}-s)d\hat{\mathcal V}_s\Big\}+\frac{t}{2}(1+\theta^{\frac12})\big(\lambda_{\ast}^2-2\lambda\lambda_{\ast}\big)+\lambda t^2-\frac16(1+\theta^{-\frac12})t^3
\end{align*}
Recall from Proposition~\ref{prop: cvrw} the definitions of $\hat\cR, \hat\cS$. We have $\hat{\mathcal U}_t=\theta^{-\frac14}\hat\cR_t$ and $\hat{\mathcal V}_t=\hat\cS_t-\theta^{-\frac12}\hat\cR_t$. Substituting this into the previous display, we find that 
\[
\log \hat{\cE}_t=\int_0^t (1-\theta^{-\frac12})\lambda_{\ast}d\hat\cR_s+\int_0^t(\lambda_{\ast}-s)d\hat\cS_s + \frac{t}{2}(1+\theta^{\frac12})\big(\lambda_{\ast}^2-2\lambda\lambda_{\ast}\big)+\lambda t^2-\frac16(1+\theta^{-\frac12})t^3.
\]
We have shown that
\[
\mathbb E\big[F\big((\mathcal U_{s}, \mathcal V_{s})_{s\le t}\big)\big] = \mathbb E\Big[F\big((\hat{\mathcal U}_{s}, \hat{\mathcal V}_{s})_{s\le t}\big)\cdot \hat\cE_t\Big],
\]
which implies \eqref{id: abs-cont} and the display below it, as we have
\[
\begin{pmatrix}
\cR^{\lambda,\theta}_{t}\\
\cS^{\lambda,\theta}_{t}
\end{pmatrix}=
\begin{pmatrix}
\theta^{\frac14}\mathcal U_{t}\\
\mathcal V_{t} + \theta^{-\frac14}\mathcal U_{t}
\end{pmatrix}\quad\text{and}\quad
\begin{pmatrix}
\hat\cR_{t}\\
\hat\cS_{t}
\end{pmatrix}=
\begin{pmatrix}
\theta^{\frac14}\hat{\mathcal U}_{t}\\
\hat{\mathcal V}_{t} + \theta^{-\frac14}\hat{\mathcal U}_{t}
\end{pmatrix},
\]
namely, $(\cR^{\lambda,\theta}_{t}, \cS^{\lambda,\theta}_{t}
)$ is the same function of $(\mathcal U_t, \mathcal V_t)$ as $(\hat\cR_{t}, \hat\cS_{t})$ is of $(\hat{\mathcal U}_t, \hat{\mathcal V}_t)$. 
\end{proof}

\section{A martingale central limit theorem}
\label{sec: mart}

We follow closely Section 7.1 in~\cite{EKbook}. Let $d\in \N$ and $C=(C_{i, j})_{1\le i, j\le d}$ be a $d\times d$ symmetric real-valued matrix which is non negative definite. 
Standard arguments from linear algebra show that we can find some real-valued matrix $K$ so that $C=K^T K$, where $K^T$ stands for the transpose of $K$. 
Let $\mathbf W=(W_i)_{1\le i\le d}$ be a standard $d$-dimensional Brownian motion, i.e.~$W_i=(W_i(t))_{t\ge 0}, 1\le i\le d$ are independent standard linear Brownian  motions. Define
\begin{equation}
\label{def: Gaus}
X(t)=K^T \mathbf W(t) K, \quad t\ge 0.
\end{equation}
Then $\mathbf X=(X(t))_{t\ge 0}$ is $d$-dimensional Gaussian process with continuous sample paths, $\mathbb E[X(t)]=0$ and $\Var(X(t))=tC$ for any $t\ge 0$. 

For each $n\ge 1$, suppose that $\mathcal F^n$ is a filtration and $\mathbf M^n=(M^n(t))_{t\ge 0}$ is an $\R^d$-valued local martingale with respect to $\mathcal F^n$ that has c\`adl\`ag sample paths and zero means. We write $M^n(t)$ as a column vector. Assume that we can find a matrix-valued random process $\mathbf A^n=(A^n(t))_{t\ge 0}$ satisfying:
\begin{enumerate}[(i)]
\item
for any $0\le s\le t$, $A^n(t)-A^n(s)$ is a $d\times d$ symmetric non negative definite matrix;
\item
$t\mapsto A^n(t)$ is c\`adl\`ag;
\item
$M^n(t)(M^n(t))^T-A^n(t), t\ge 0$ is a local martingale with respect to $\mathcal F^n$.
\end{enumerate}
Let $\|\cdot\|$ stand for the Euclidean norms for $d$-dimensional vectors as well as $d\times d$ matrices. Assume that the following conditions are satisfied for all $t\ge 0$: 
\begin{align}\label{cvcon: var}
A^n(t) \xrightarrow[]{n\to\infty} tC \quad\text{in probability}. \\ \label{cvcon: An}
\mathbb E\Big[\sup_{s\le t}\big\|A^n(s)-A^n(s-)\big\|\Big] \xrightarrow[]{n\to\infty} 0.\\ \label{cvcon: jump}
\mathbb E\Big[\sup_{s\le t}\big\|M^n(s)-M^n(s-)\big\|^2\Big] \xrightarrow[]{n\to\infty} 0.
\end{align}
Then Theorem 1.4 in~\cite{EKbook}, Chapter 7 says that 
\[
\mathbf M^n \Longrightarrow \mathbf X \quad\text{ in }\mathbb D^{\ast}(\R_+, \R^d), 
\]
where $\mathbf X=(X(t))_{t\ge 0}$ is defined in~\eqref{def: Gaus}. 

\section{Some properties of the Binomial and Poisson distributions}
\label{sec: binom}

Let $N\in \N$, $p\in [0, 1]$ and $c\in (0, \infty)$. Throughout this part, $X$ denotes a Binom$(N, p)$ variable and $Y$ a Poisson$(c)$ variable. 

\paragraph{Moments of Binomial and Poisson distributions. }
For all $k\in \N$ and $k\le N$, we have
\begin{equation}
\label{id: binom-moment}
\mathbb E\big[(X)_{k}\big] = \mathbb E\big[X(X-1)(X-2)\cdots(X-k+1)] = (N)_{k}\cdot p^{k}\le (Np)^{k}.  
\end{equation}
Similarly, for all $k\in \N$, we have
\begin{equation}
\label{id: poisson-moment'}
\mathbb E\big[(Y)_{k}\big] = c^{k}. 
\end{equation}
From \eqref{id: binom-moment} (resp.~\eqref{id: poisson-moment'}), we can deduce an expression for $\mathbb E[X^{k}]$ (resp.~for $\mathbb E[Y^{k}]$), since for any $\Z_+$-valued random variable $Z$, we have
\begin{equation}
\label{id: poisson-moment}
\mathbb E\big[Z^{k}\big] = \sum_{j=1}^{k} {k\brace j}\mathbb E\big[(Z)_{j}\big],
\end{equation}
where ${k \brace j}$ is the Stirling number of the second kind. 
This implies the following asymptotics for $\mathbb E[X^{k}]$. Assume that $N\to\infty$ and $p\to 0$. 
\begin{itemize}
\item
If $N p\to c\in (0, \infty)$, then $\mathbb E[X^{k}]\to \mathbb E[Y^{k}]$ for all $k\ge 1$. 
\item
If $Np\to \infty$, then $\mathbb E[X^{k}]\sim (Np)^{k}$ for all $k\ge 1$.
\item
If $Np\to 0$, then $\mathbb E[X^{k}] \sim Np$ for all $k\ge 1$ and for all $k\ge 2$, $\mathbb E[X^{k}]-Np\sim {k\brace 2}(Np)^{2}$. Note that ${k\brace 2}=2^{k-1}-1$. 
\end{itemize}
As a consequence, supposing $Np\to 0$, we have
\begin{equation}
\label{bd: binom2}
\mathbb E\Big[\Big((1+X)_{3}\Big)^{2}\Big]=\mathbb E[(X+1)^{2}X^{2}(X-1)^{2}] = \mathbb E[X^{6}-2X^{4}+X^{2}] \sim 18 (Np)^{2}.
\end{equation}

\paragraph{Tail estimates for Poisson distributions.}
Let $k\in \N$. We observe that
\begin{equation}
\label{bd: pois-tail}
\mathbb P(Y\ge k) =e^{-c} \sum_{j\ge k}\frac{c^{j}}{j!} = e^{-c}\frac{c^{k}}{k!}\sum_{j\ge k}\frac{k!\,c^{j-k}}{j!}\le \frac{c^{k}}{k!}. 
\end{equation}
We have
\begin{equation}
\label{bd: binom7'}
\mathbb E[Y|Y\ge 2] \le  (\mathbb P(Y=2))^{-1}\big(\mathbb E[Y]-\mathbb P(Y=1)\big)=\frac{c(1-e^{-c})}{e^{-c}c^2/2}\le 2e^{c}, 
\end{equation}
where we have used $1-e^{-x}\le x$ for $x\ge 0$. 
We also have
\begin{equation}
\label{bd: binom7''}
\mathbb E[Y|Y< 2] =\frac{ \mathbb P(Y=1)}{\mathbb P(Y\le 1)}=\frac{c}{1+c}\le c=\mathbb E[Y]. 
\end{equation}

\paragraph{Tail estimates for Binomial distributions.}
Combining~\eqref{bd: pois-tail} with \eqref{bd: binom-tv}, we find the following universal tail bound for Binomial distributions that is valid for all $N\ge 1$ and $p\in [0, 1]$: 
\begin{equation}
\label{bd: binom-tail}
\mathbb P(X\ge k)\le Np^2 + \mathbb P(\text{Poisson}(Np)\ge k)\le Np^2 + \tfrac{1}{k!}(Np)^{k}, \quad k\in \N.  
\end{equation}
Let us recall Chernoff's bounds for $X$ (\cite{concen}): for any $\lambda>0$: 
\begin{equation}
\label{bd: cc-binom}
\mathbb P(X\le \mathbb E[X]-\lambda)\le\exp\Big(-\frac{\lambda^2}{2\mathbb E[X]}\Big), \quad 
\mathbb P(X\ge \mathbb E[X]+\lambda)\le \exp\Big(-\frac{\lambda^2}{2\mathbb E[X]+2\lambda/3}\Big).
\end{equation}

\providecommand{\bysame}{\leavevmode\hbox to3em{\hrulefill}\thinspace}
\providecommand{\MR}{\relax\ifhmode\unskip\space\fi MR }
\providecommand{\MRhref}[2]{%
  \href{http://www.ams.org/mathscinet-getitem?mr=#1}{#2}
}
\providecommand{\href}[2]{#2}

%%%%%%%%%%%%%%%%%%%%%%%%%%%%%%%%%%%%%%%%%%%%%%%%%%%%%%%%%%%%%%%%%%%
%%                                                               %%
%% You may add acknowledgments (optional).                       %%
%%                                                               %%
%%%%%%%%%%%%%%%%%%%%%%%%%%%%%%%%%%%%%%%%%%%%%%%%%%%%%%%%%%%%%%%%%%%
\begin{acks}
%We are grateful to Martin Hairer who provided a nice \texttt{MR} macro and to S\'ebastien Gou\"ezel for his useful comments on the internals of the class file.
I am very grateful to the referee and the editor for their insightful and detailed suggestions. The author’s research is partially funded by EPSRC, United Kingdom EP/W033585/1 grant and receives support from the Dr Perry James (Jim) Browne Research Centre. No new data were created or analysed during this study. Data sharing is not applicable to this article. 
\end{acks}

%%%%%%%%%%%%%%%%%%%%%%%%%%%%%%%%%%%%%%%%%%%%%%%%%%%%%%%%%%%%%%%%%%%
%%                                                               %%
%% You have reached the end of your document.                    %%
%%                                                               %%
%%%%%%%%%%%%%%%%%%%%%%%%%%%%%%%%%%%%%%%%%%%%%%%%%%%%%%%%%%%%%%%%%%%

\end{document}